\documentclass[a4paper,3p,10pt]{elsarticle}

\usepackage{amsmath}
\usepackage{cases}
\usepackage{amssymb}
\usepackage{amsthm}
\usepackage{bm}
\usepackage{mathtools}
\usepackage{physics}
\usepackage{booktabs}

\usepackage[noabbrev]{cleveref}

\biboptions{sort&compress}

\theoremstyle{plain}
\newtheorem{example}{Example}[section]
\AddToHook{env/example/begin}{\pushQED{\hfill$\diamond$}}
\AddToHook{env/example/end}{\popQED}

\newcommand\SCALAR[1]{\mathrm{#1}}
\newcommand\hSCALAR[1]{\hat{\SCALAR{#1}}}

\newcommand\VEC[1]{\vb*{#1}}
\newcommand\hVEC[1]{\hat{\VEC{#1}}}

\newcommand\MAT[1]{\mathrm{#1}}

\newcommand\eye{\MAT{I}}

\newcommand\ncp{\text{n}_{\text{sf}}}
\newcommand\strain{\bm{\varepsilon}}

\usepackage[dvipsnames]{xcolor}
\usepackage{subcaption}
\usepackage{tikz}
\usetikzlibrary{calc,positioning,spy}
\usetikzlibrary{patterns}
\usetikzlibrary{patterns.meta}
\usetikzlibrary{decorations.markings,decorations.shapes,decorations.pathmorphing,decorations.pathreplacing}
\usetikzlibrary{decorations.text,matrix}
\usetikzlibrary{external}
\usepackage{pgfplots}
\usepgfplotslibrary{units,groupplots}
\usepackage{pgfplotstable}

\usepackage{algorithm2e}
\RestyleAlgo{ruled}
\SetKwInput{KwInput}{Input}
\SetKwInput{KwOutput}{Output}
\SetKwComment{Comment}{\# }{}

\pgfdeclarelayer{background}
\pgfsetlayers{background,main}

\pgfplotsset{
    grid=both,
    grid style={dashed, gray!30},
    ylabel near ticks,
    xlabel near ticks,
    enlargelimits,
	legend style={fill=white, fill opacity=0.6, draw opacity=1,draw=black,text opacity=1, font=\normalsize},
    legend cell align={left},
    legend image with text/.style={
    legend image code/.code={%
        \node[anchor=center] at (0.3cm,0cm) {#1};
    }
},
legend pos = north west,
tick label style={font=\normalsize},
label style={font=\normalsize},
every axis plot/.append style={line width=0.5pt},
every axis plot post/.append style={
    every mark/.append style={solid,scale=0.5,line width=0.5pt}
                                    }
}

\definecolor{col1}{HTML}{b3274a}
\definecolor{col2}{HTML}{2746b3}
\definecolor{col3}{HTML}{27b33b}
\definecolor{col4}{HTML}{b38e27}

\definecolor{red1}{HTML}{b3274a}
\definecolor{red2}{HTML}{d2748b}
\definecolor{green1}{HTML}{2ed043}
\definecolor{green2}{HTML}{4e6c52}
\definecolor{blue1}{HTML}{2746b3}
\definecolor{yellow1}{HTML}{b38e27}
\definecolor{yellow2}{HTML}{b29f6b}

\pgfplotsset{
TP_AT1_Order2_coarse/.style={thick,col1,densely dashed,mark=square,mark options={solid}},
TP_AT1_Order4_coarse/.style={thick,col2!90!white,densely dashed,mark=square,mark options={solid}},
TP_AT2_Order2_coarse/.style={thick,col3,densely dashed,mark=square,mark options={solid}},
TP_AT2_Order4_coarse/.style={thick,col4,densely dashed,mark=square,mark options={solid}},
TP_AT1_Order2_fine/.style={thick,col1,densely dotted,mark=o,mark options={solid}},
TP_AT1_Order4_fine/.style={thick,col2!90!white,densely dotted,mark=o,mark options={solid}},
TP_AT2_Order2_fine/.style={thick,col3,densely dotted,mark=o,mark options={solid}},
TP_AT2_Order4_fine/.style={thick,col4,densely dotted,mark=o,mark options={solid}},
THB_AT1_Order2_coarse_implicit/.style={thick,col1,densely dashed,mark=square,mark options={solid}},
THB_AT1_Order4_coarse_implicit/.style={thick,col2!90!white,densely dashed,mark=square,mark options={solid}},
THB_AT2_Order2_coarse_implicit/.style={thick,col3,densely dashed,mark=square,mark options={solid}},
THB_AT2_Order4_coarse_implicit/.style={thick,col4,densely dashed,mark=square,mark options={solid}},
THB_AT1_Order2_fine_implicit/.style={thick,col1,densely dotted,mark=o,mark options={solid}},
THB_AT1_Order4_fine_implicit/.style={thick,col2!90!white,densely dotted,mark=o,mark options={solid}},
THB_AT2_Order2_fine_implicit/.style={thick,col3,densely dotted,mark=o,mark options={solid}},
THB_AT2_Order4_fine_implicit/.style={thick,col4,densely dotted,mark=o,mark options={solid}},
THB_AT1_Order2_coarse_explicit/.style={thick,col1!50!white,densely dashed,mark=square,mark options={solid}},
THB_AT1_Order4_coarse_explicit/.style={thick,col2!50!white,densely dashed,mark=square,mark options={solid}},
THB_AT2_Order2_coarse_explicit/.style={thick,col3!50!white,densely dashed,mark=square,mark options={solid}},
THB_AT2_Order4_coarse_explicit/.style={thick,col4!50!white,densely dashed,mark=square,mark options={solid}},
THB_AT1_Order2_fine_explicit/.style={thick,col1!50!white,densely dotted,mark=o,mark options={solid}},
THB_AT1_Order4_fine_explicit/.style={thick,col2!50!white,densely dotted,mark=o,mark options={solid}},
THB_AT2_Order2_fine_explicit/.style={thick,col3!50!white,densely dotted,mark=o,mark options={solid}},
THB_AT2_Order4_fine_explicit/.style={thick,col4!50!white,densely dotted,mark=o,mark options={solid}},
THB_AT1_Order2_coarse_hybrid/.style={thick,col1!50!black,densely dashed,mark=square,mark options={solid}},
THB_AT1_Order4_coarse_hybrid/.style={thick,col2!50!black,densely dashed,mark=square,mark options={solid}},
THB_AT2_Order2_coarse_hybrid/.style={thick,col3!50!black,densely dashed,mark=square,mark options={solid}},
THB_AT2_Order4_coarse_hybrid/.style={thick,col4!50!black,densely dashed,mark=square,mark options={solid}},
THB_AT1_Order2_fine_hybrid/.style={thick,col1!50!black,densely dotted,mark=o,mark options={solid}},
THB_AT1_Order4_fine_hybrid/.style={thick,col2!50!black,densely dotted,mark=o,mark options={solid}},
THB_AT2_Order2_fine_hybrid/.style={thick,col3!50!black,densely dotted,mark=o,mark options={solid}},
THB_AT2_Order4_fine_hybrid/.style={thick,col4!50!black,densely dotted,mark=o,mark options={solid}},
PFSolver/.style={fill=red1!50,draw=red1},
PFAssembler/.style={fill=red2!50,draw=red2},
ElSolver/.style={fill=green1!50,draw=green1},
ElAssembler/.style={fill=green2!50,draw=green2},
Projection/.style={fill=blue1!50,draw=blue1},
Marking/.style={fill=yellow1!50,draw=yellow1},
Refinement/.style={fill=yellow2!50,draw=yellow2}
}

\pgfplotscreateplotcyclelist{MyCycleList}{%
{col1,mark options={solid}},
{col2,mark options={solid}},
{col3,mark options={solid}},
}

\newcommand\NEW{}

\crefname{alg}{alg.}{algs.}
\Crefname{alg}{Algorithm}{Algorithms}
\crefname{algocf}{algorithm}{algorithms}
\Crefname{algocf}{Algorithm}{Algorithms}

\usepackage[compact]{titlesec}

\begin{document}

\begin{frontmatter}

\title{Adaptive isogeometric analysis of high-order phase-field fracture based on THB-splines}

\author[label1]{\corref{cor1}H.M. Verhelst}
\ead{h.m.verhelst@tue.nl}
\cortext[cor1]{Corresponding Author. }
\author[label2]{L. Greco}
\author[label2]{A. Reali}
\address[label1]{Eindhoven University of Technology, Department of Mechanical Engineering, Traverse, PO Box 513, 5600 MB Eindhoven, Eindhoven, The Netherlands}
\address[label2]{University of Pavia, Department of Civil Engineering and Architecture, Via Adolfo Ferrata 3, 27100 PV, Pavia, Italy}

\begin{abstract}
In recent decades, the study of fracture propagation in solids has increasingly relied on phase-field models. Several recent contributions have highlighted the potential of this approach in both static and dynamic frameworks. However, a major limitation remains the high computational cost. Two main strategies have been identified to mitigate this issue: the use of locally refined meshes and the adoption of higher-order models. In this work, leveraging Truncated Hierarchical B-splines (THB-splines), we introduce adaptive simulations of higher-order phase-field formulations (AT1 and AT2), focusing primarily on two-dimensional fracture problems.
\end{abstract}


\end{frontmatter}


\section{Introduction}
\label{sec:introduction}
Fracture mechanics is a relatively young discipline whose foundations lie in the pioneering work of Griffith \cite{Griffith1921}, who introduced the first energy-based approach to the theory of sharp brittle fracture. This framework was subsequently extended and widely adopted by Irwin \cite{irwin1957analysis} within the field of engineering. The development of this field has been strongly motivated by the profound societal impact of fracture phenomena, both in terms of safety implications and associated economic costs \cite{reed1983economic}.\\ 

This pressing need, combined with the concurrent advancement of computational technologies, has led to the development of dedicated numerical methods for the analysis of fracture problems. Among these, the most widely adopted are the extended finite element method (XFEM) \cite{belytschko1999elastic} and, more recently, the phase-field approach \cite{BOURDIN2000797}. The latter has become one of the most widely used methods to date, owing to its ability to capture complex fracture phenomena such as branching, coalescence, and merging. However, its main drawback lies in the high computational cost, which has so far limited its competitiveness compared to alternative approaches.\\

The need to reduce the computational cost  of phase-field fracture simulations is a well-recognized issue in the literature. Two main strategies have been explored in this context: the use of higher-order formulations and the use of adaptive meshing.
Firstly, the use of higher-order phase-field formulations follows the pioneering work of \citeauthor{BORDEN2014100} \cite{BORDEN2014100}, which involves the use of higher-order energy dissipation functionals based on the Cahn-Hilliard phase-field model. As shown by \citeauthor{greco2024higher} \cite{greco2024higher}, this formulation yields significant reduction of the computational costs, since coarser meshes can be used to achieve the same accuracy. In addition, another factor influencing the computational costs is the regularization functional, which determines the shape of the phase-field. The choice of this functional, originally proposed by \citeauthor{Ambrosio1990999} \cite{Ambrosio1990999}, can greatly reduce computational costs, as shown by \citeauthor{greco2024higher} \cite{greco2024higher}.\\
The computational advantages of the higher-order phase-field fracture models are largely attributed to the use of the fourth-order Cahn-Hilliard phase-field model over the second-order Allen-Cahn model. Since the former is of higher order, classical discretization based on $C^0$ continuity across element interfaces requires to introduce an auxiliary field for the solution gradient, doubling the computational costs of the model. Alternatively, higher continuity discretizations, such as spline-based Isogeometric Analysis, naturally provide higher-order derivatives, making them an attractive alternative to classical Finite Element Analysis for phase-field modeling \cite{HUGHES20054135, CottrellCMAME2007}.\\
One problem related to phase-field fracture modeling is the imposition of the irreversibility of the damage field. In the literature, several approaches have been proposed to ensure this, such as Strain history variable approach \cite{miehe2010IJNME}, Penalty method \cite{GERASIMOV2019990} and Projected Successive Over-Relaxation (PSOR) \cite{MARENGO2021114137}.\\

Secondly, besides the use of higher-order formulations to increase the efficiency of phase-field brittle fracture simulations, a second strategy relies on the use of adaptive mesh refinement techniques. This approach enables mesh refinement to be localized in the vicinity of the crack, thereby reducing the number of degrees of freedom in regions far from the fracture. The key idea is to employ the phase-field variable -- its value, variation, or gradient -- as an indicator to drive the refinement strategy. Since higher-order phase-field models benefit from isogeometric analysis, as mentioned above, we limit ourselves to adaptive spline techniques in the sequel.\\
Over the last years, different adaptive spline constructions have been successfully employed for adaptive phase-field fracture simulations. Firstly, Locally Refined (LR)-splines have been used by \citeauthor{PROSERPIO2020113363} \cite{PROSERPIO2020113363} in the context of phase-field fracture simulations for shell structures. In addition, Li and co-authors \cite{li2023dynamic, li2023modeling, li2022adaptive} studied fracture in rock-like materials -- both in static and dynamic -- using LR-splines. In a three-dimensional setting, however, LR-splines have hardly been used for isogeometric analysis, let alone to simulate phase-field fracture. Using Polynomial Splines over Hierarchical T-meshes (PHT-splines), phase-field fracture induced by thermal and electrical loading in composite materials has been simulated in three-dimensional space by \cite{kiran2022adaptive, kiran2023adaptive, kiran2024phase, xu2024adaptive}. In addition, mesh-free methods have been employed by \citeauthor{nguyen2020adaptive} \cite{nguyen2020adaptive} and \citeauthor{li2020phase} \cite{li2020phase} for adaptive fracture analysis in polycrystalline materials. Lastly, Truncated-Hierarchical B-splines \cite{giannelli2012thb} are known for their simple construction and refinement algorithms, generalized for domains of any dimension. Although this spline construction has been applied for adaptive simulation of tumor growth \cite{lorenzo2017hierarchically} and fluid immiscibility problems \cite{bracco2023adaptive} based on phase-field models, application of these splines to phase-field fracture is yet undiscovered.\\

{\NEW In this work, we present a numerical framework for the simulation of brittle fracture problems through higher-order phase-field models and adaptive meshing. The novelty of the proposed framework lies in the combination of a marking and refinement strategy dedicated to phase-field fracture simulations, based on the observations of fracture suddenness, irreversibility and cross-talk, as well as load stepping strategies to handle the suddenness of brittle fracture. In combination with the higher-order phase-field formulations from \cite{greco2024higher, greco_at1}, the presented framework leads to a significant reduction of the computational costs with respect to non-adaptive tensor-product B-splines, as will be shown with numerical experiments.}\\ 

The outline of this paper is as follows: \cref{sec:phase-field_fracture} provides a background on phase-field fracture models, including the higher-order formulations from \cite{greco2024higher, greco_at1}. Thereafter, \cref{sec:THB} provides a background on THB-splines and admissible refinement strategies. The novelty of this paper is presented in \cref{sec:adaptive_IGA}, where we define adaptive meshing strategies dedicated to phase-field fracture problems. \Cref{sec:benchmarks} presents numerical benchmarks to assess the performance of the proposed adaptive phase-field fracture model, comparing the results to tensor-product-based ones for different phase-field formulations. Finally, \cref{sec:conclusions} draws the conclusions of this work and provides perspectives for future research. Lastly, an appendix provides an algorithmic overview of the proposed adaptive phase-field fracture model.

\section{Isogeometric phase-field brittle Fracture}
\label{sec:phase-field_fracture}
In this section, we present the isogeometric model for brittle phase-field fracture used in the present paper. The section starts with the introduction of the variational formulation used for brittle phase-field fracture in \cref{subsec:PF_variationalForm}. This section provides energy functionals and their variations with respect to the displacements and the phase-field, needed for discretization. \Cref{subsec:PF_discretization} elaborates on the temporal and spatial discretization of the variational formulations provided in \cref{subsec:PF_variationalForm}. Using the discrete operators, \cref{subsec:PF_solution_scheme} elaborates on the staggered solution scheme employed in this work. Lastly, \cref{subsec:PF_initialization} provides a background on the initialization of a phase-field based on an initial crack.

\subsection{Phase-field variational formulations}
\label{subsec:PF_variationalForm}
In this section, we briefly recall the phase-field formulations considered in this work. Following the framework introduced in \cite{BORDEN2014100}, let us denote by $\Omega \subset \mathbb{R}^N$ the reference configuration, where $N$ represents the spatial dimension of the problem. Dirichlet boundary conditions are prescribed on $\partial \Omega_D \subseteq \partial \Omega$, while Neumann conditions are imposed on $\partial \Omega_N \subseteq \partial \Omega$, with $\partial \Omega = \partial \Omega_D \cup \partial \Omega_N$ and $\partial \Omega_D \cap \partial \Omega_N = \emptyset$. The space of admissible displacements is thus defined as
\begin{equation}
\mathcal{U}_{\VEC{u}} := \{ \VEC{u} \in H^1(\Omega; \mathbb{R}^N) \, : \, \VEC{u} = \bar{\mathbf{u}} \text{ on } \partial \Omega_D \}\,,
\end{equation}
whereas the phase-field variable $\SCALAR{d}$, representing a smooth approximation of the sharp crack topology, belongs to
\begin{align*}
\mathcal{U}_{\SCALAR{d}}^{ii} &:= \{ \SCALAR{d} \in H^{1}(\Omega) \, : \, \SCALAR{d} \ge 0 \} \,, \\
\mathcal{U}_{\SCALAR{d}}^{iv} &:= \{ \SCALAR{d} \in H^{2}(\Omega) \, : \, \SCALAR{d} \ge 0 \} \,,
\end{align*}
for the second- and fourth-order formulations, respectively.

Strictly speaking, the phase-field variable $\SCALAR{d}$ should take values within $[0,1]$, corresponding to the intact ($\SCALAR{d}=0$) and fully broken ($\SCALAR{d}=1$) material states. In the present work, we relax this constraint and allow $\SCALAR{d} \ge 0$, which \textit{a priori} permits $\SCALAR{d} > 1$. However, during the staggered evolution scheme, the minimization process naturally restricts $\SCALAR{d}$ to values within $[0,1]$. This weaker constraint is therefore adopted, as it proves to be more convenient from both the mathematical and numerical standpoints.

\subsubsection{Energy functionals}\label{subsec:functionals}
We assume small strains, i.e., $\strain = \frac{1}{2}\qty(\bm{\nabla} \VEC{u} + \bm{\nabla}^\top \VEC{u})$. Accordingly, for quasi-static brittle fracture, the total energy functionals for the second- and fourth-order formulations are expressed as
\begin{equation}\label{eq:pi}
	\Pi^{ii,iv} (\VEC{u},\SCALAR{d}) := \mathcal{E}(\VEC{u},\SCALAR{d})	+  G_c\,\mathcal{D}^{ii,iv} (\SCALAR{d}) - \mathcal{W}^{\text{ext}}(\VEC{u})\,.
\end{equation}
Here, $\mathcal{W}^{\text{ext}}(\VEC{u}):=\int_{\Omega} \mathbf{b} \cdot \VEC{u} \, d\Omega + \int_{\partial \Omega_{\text{N}}} \mathbf{t} \cdot \VEC{u} \, d\partial \Omega$ denotes the external work functional, where $\mathbf{b}$ and $\mathbf{t}$ are the body and traction force densities, respectively. In \cref{eq:pi}, $\mathcal{E}(\VEC{u},\SCALAR{d})$ represents the elastic strain energy and is defined as
\begin{equation}\label{eq:elastic_strain_en}
	\mathcal{E}(\VEC{u},\SCALAR{d}) := \int_{\Omega} \psi(\strain(\VEC{u}),\SCALAR{d}) \, d\Omega
	=\int_{\Omega} \big[\omega(\SCALAR{d})\, \psi^+_0(\strain(\VEC{u})) +  \psi^-_0(\strain(\VEC{u}))\big] \, d\Omega\,,
\end{equation}
where, following \cite{Comi2001IJSS, Amor2009JMPS}, the positive and negative parts of the strain energy density are given by
\begin{equation}\label{eq:free_en_split}
		\psi^+_0(\strain(\VEC{u})) := \tfrac{1}{2}\big[K \left(\varepsilon^+_v\right)^2 + \mu\,\vert\strain_d\vert^2\big]\,, \qquad
		\psi^-_0(\strain(\VEC{u})) := \tfrac{1}{2}\big[K \left(\varepsilon^-_v\right)^2\big]\,.
\end{equation}
The volumetric strain is $\varepsilon_v = \strain : \eye$, where $\eye$ is the identity tensor, and $\varepsilon^\pm_v := \langle \varepsilon_v \rangle_\pm$ denote its positive and negative parts. The deviatoric strain is defined as $\strain_d := \strain - \frac{1}{N} \varepsilon_v \eye$, while $K > 0$ and $\mu > 0$ are the bulk and shear moduli, respectively. The monotonically decreasing degradation function $\omega(\SCALAR{d}) = (1 - \SCALAR{d})^2 + \eta$ governs the reduction of stored elastic energy as damage evolves, and satisfies $\omega(0) = 1 + \eta$, $\omega(1) = \eta$, and $\omega'(1) = 0$. The small positive parameter $\eta \ll 1$ prevents full degradation of $\psi^+_0$ and ensures numerical stability by retaining a residual stiffness when $\SCALAR{d} \to 1$, as discussed in \cite{miehe2010IJNME}.

The fracture energy contributions $\mathcal{D}^{ii,iv}$ in \cref{eq:pi} describe the regularized crack surface densities. In this work, we consider two different families of functionals, namely the AT1 and AT2 models, originally proposed by Ambrosio and Tortorelli \cite{Ambrosio1990999}. The corresponding formulations are given as follows:
\begin{itemize}
	\item AT1 model:
	\begin{subequations}\label{eq:fracture_density_funct_AT1}
		\begin{align}
			\mathcal{D}_{\text{AT1}}^{ii} (\SCALAR{d}) &:= \int_\Omega \frac{3}{8l_0}\big(\SCALAR{d} + l_0^2 \vert \bm{\nabla}\SCALAR{d} \vert^2\big) \, d\Omega\,, \label{eq:2nd_order_AT1}\\[4pt]
			\mathcal{D}_{\text{AT1}}^{iv} (\SCALAR{d}) &:= \int_\Omega \frac{1}{c_\rho\, l_0}\big(\SCALAR{d} + \tfrac{l_0^2}{2}\vert \bm{\nabla}\SCALAR{d} \vert^2 + \rho (\nabla^2 \SCALAR{d})^2\big) \, d\Omega\,. \label{eq:4th_order_AT1}
		\end{align}
	\end{subequations}
	\item AT2 model:
	\begin{subequations}\label{eq:fracture_density_funct_AT2}
		\begin{align}
			\mathcal{D}_{\text{AT2}}^{ii} (\SCALAR{d}) &:= \int_\Omega \frac{1}{2l_0}\big(\SCALAR{d}^2 + l_0^2 \vert \bm{\nabla}\SCALAR{d} \vert^2\big) \, d\Omega\,, \label{eq:2nd_order_AT2}\\[4pt]
			\mathcal{D}_{\text{AT2}}^{iv} (\SCALAR{d}) &:= \int_\Omega \frac{1}{2l_0}\big(\SCALAR{d}^2 + \tfrac{l_0^2}{2}\vert \bm{\nabla}\SCALAR{d} \vert^2 + \tfrac{l_0^4}{16}(\nabla^2 \SCALAR{d})^2\big) \, d\Omega\,, \label{eq:4th_order_AT2}
		\end{align}
	\end{subequations}
\end{itemize}
In these equations, $c_{\rho}$ is the normalization constant associated with the optimal 1D profile. In this study, we adopt $\rho = 1$ and $c_{\rho} = 4.4485$, as suggested in \cite{greco_at1}. The internal length scale $l_0$ controls the width of the diffused crack zone; as $l_0 \to 0$, the regularized phase-field formulation $\Gamma$-converges to the classical sharp crack representation (see, e.g., \cite{NEGRI2020112858}).
It is possible to summarize the functional as follows:
\begin{equation} \label{eq: D general}
\mathcal{D}(\SCALAR{d}) := \int_\Omega \frac{1}{c_{\rho}} \qty(\frac{\SCALAR{d}^\beta}{l_0} + \gamma l_0\vert \bm{\nabla}\SCALAR{d} \vert^2 + \delta l_0^3(\nabla^2 \SCALAR{d})^2) \, d\Omega
\end{equation}

where the coefficients $\beta, \, \gamma \, \text{and} \, \delta$ are selected based on the regularization model (AT1 or AT2) and functional grade (second- or fourth-order). 

\subsubsection{Evolution in terms of energy variations}
\label{subsubsec:variations}

In this section, the variational formulation of fracture is presented. It is worth noting that an in-depth analysis of the variational problem lies beyond the scope of this work. However, for the sake of completeness, and with reference to \cite{greco2024higher, greco_at1}, we report below the system of equations governing the evolution of the problem in the continuum setting:

\begin{equation} \label{eq:evol-sint}
	\begin{cases}
		\partial_{\VEC{u}} \Pi^{ii,iv} ( \VEC{u}, \SCALAR{d}) = \bm{0}\,,  \\[4pt]
		\partial_{\SCALAR{d}} \Pi^{ii,iv} ( \VEC{u} , \SCALAR{d})  \ge 0\,,  \\[4pt]
		\partial_{\SCALAR{d}} \Pi^{ii,iv} ( \VEC{u}, \SCALAR{d} ) [ \dot{\SCALAR{d}} ] = 0\,, \quad  \dot{\SCALAR{d}} \ge 0\,.
	\end{cases}
\end{equation}

The system of \cref{eq:evol-sint} governs the evolution through the Karush--Kuhn--Tucker (KKT) conditions, from which the irreversibility condition of the problem directly follows from the third equation of the system. This variational formulation provides a unified energetic framework for the description of quasi-static fracture evolution, ensuring both mechanical equilibrium and the thermodynamic consistency of the damage growth process.\\

{\NEW \Cref{eq:evol-sint} can be solved using a penalty method proposed by \cite{GERASIMOV2019990}, penalizing energy contributions corresponding to negative values in the phase-field variable. This method, however, introduces a dependency on a penalty parameter. Alternatively, when considering \cref{eq:evol-sint} as a Symmetric Linear Complementarity Problem (SLCP), the irreversibility condition can be solved using the Projected Successive Over-Relaxation (PSOR) algorithm proposed by \cite{MARENGO2021114137}. The advantage of using the PSOR is that it is parameter-free, in contrast to the penalty-based approach. However, the PSOR algorithm solves the SLCP in a serial fashion, which can become a bottleneck when solving phase-field fracture problems on many computational ranks. Nevertheless, this adopts the serial PSOR algorithm since the computations are distributed on a relatively low number of computational ranks. The development of a parallel PSOR algorithm is considered an area for future research.}

\subsection{Space and time numerical discretization}
\label{subsec:PF_discretization}

Due to the irreversibility of the fracture process, the problem becomes path-dependent, and its solution requires a time integration of the mechanical model while enforcing the irreversibility condition throughout the evolution. Again, we refer the interested reader to \cite{greco2024higher, greco_at1} for the discritized version of \cref{{eq:evol-sint}}. The third equation of \eqref{eq:evol-sint}, in the discrete form, represents the irreversibility condition, solved \textit{via} Projected Successive Over-Relaxation (PSOR) Algorithm \cite{MARENGO2021114137}.\\

For what concerns the spatial discretization, a Galerkin isogeometric approach based on $C^1$-continuous quadratic basis functions is adopted (see, e.g., \cite{HUGHES20054135, pieg1996nurbs, CottrellCMAME2007}). This choice allows a consistent numerical approximation of the Laplacian operator. As will be discussed in \cref{sec:THB,sec:adaptive_IGA}, the adopted isogeometric framework relies on Truncated Hierarchical B-Splines (THB-splines), thereby enabling the advantages of local adaptivity.\\

We denote by $R_i^{\VEC{u}} = R_i^{\SCALAR{d}} = R_i : \Omega \to \mathbb{R}$ the isogeometric basis functions used for the approximation of both the displacement field $\VEC{u}$ and the phase-field variable $\SCALAR{d}$. Consequently, the approximate displacement field $\VEC{u}^h$ and phase-field variable $\SCALAR{d}^h$ are expressed as linear combinations of the IGA basis functions and the corresponding control variables, $\hVEC{u}_i \in \mathbb{R}^3$ and $\hSCALAR{d}_i \in \mathbb{R}$, respectively:
\begin{equation}
	\begin{aligned}
		\VEC{u}(\VEC{x} , t)&\approx\VEC{u}^h(\VEC{x},t) = \sum^{\ncp}_{i}  R_i(\VEC{x}) \, \hVEC{u}_{i}(t)\,,\\
		\SCALAR{d}(\VEC{x} , t)&\approx\SCALAR{d}^h(\VEC{x},t)  = \sum^{\ncp}_{i} R_i (\VEC{x}) \, \hSCALAR{d}_{i}(t)\,,
	\end{aligned}
	\label{eq:approx_global_u_d}
\end{equation}
where $\ncp$ is the total number of basis functions associated with the spatial discretization.\\

The approximated \eqref{eq:approx_global_u_d} are substituted into the weak form of the elasto-static problem and into the Karush–Kuhn–Tucker conditions. In this framework, the time-discretized variational formulation of the coupled system can be rewritten in matrix form as:
\begin{subnumcases}
		\MAT{K}( \VEC{u}^h_{n+1} , \SCALAR{d}^h_{n+1} ) \,  \VEC{u}_{n+1} - \VEC{F}_{n}^{\text{ext}} = \mathbf{0}, \label{eq:matrix_problem_u}\\
		\MAT{Q}^{\text{II,IV}} ( \VEC{u}^h_{n+1} ) \, \Delta \SCALAR{d}^h - \VEC{Q}^{\text{II,IV}} ( \VEC{u}^h_{n+1} , \SCALAR{d}^h_{n}) \ge \mathbf{0},\label{eq:matrix_problem_d1}\\
		\big( \MAT{Q}^{\text{II,IV}} ( \VEC{u}^h_{n+1}  ) \, \Delta \SCALAR{d}^h   - \VEC{Q}^{\text{II,IV}} ( \VEC{u}^h_{n+1} , \SCALAR{d}^h_{n}) \big) \cdot \Delta \SCALAR{d}^h= 0,\label{eq:matrix_problem_d2}
		\quad
		\Delta \SCALAR{d}^h \ge 0.
        \label{eq:matrix_problem_n}
\end{subnumcases}

The phase-field matrix and vector are defined as follows:
\begin{equation}\label{eq:Q_q}
	\MAT{Q}^{\text{II,IV}} := \bm{\Psi} (\hat{\mathbf{u}}) + G_c \,\bm{\Phi}^{\text{II,IV}}
	\,, \quad
	{\VEC{Q}}^{\text{II,IV}} := \MAT{Q}^{\text{II,IV}} \hat{\mathbf{d}} - \bm{\psi} (\hat{\mathbf{u}}) - \bm{\phi}^{\text{II,IV}}\,,
\end{equation}
where:
\begin{equation} \label{eq:SLCP_free_energy_matrix_vector}
	(\bm{\Psi}(\hat{\mathbf{u}}))_{\mathbf{ij}} := \int_{\Omega} 2\,\psi_0^+(\hat{\mathbf{u}}) \, R_{\mathbf{i,p}} \, R_{\mathbf{j,p}} \; \text{d} \Omega
	\,, \quad
	(\bm{\psi}(\hat{\mathbf{u}}))_{\mathbf{i}} := \int_{\Omega} 2\,\psi_0^+(\hat{\mathbf{u}}) \, R_{\mathbf{i,p}} \; \text{d} \Omega\,.
\end{equation}
The $\bm{\Phi}$ matrix, is the Galerkin approximation of eq: \eqref{eq: D general}. If $\beta = 2$ the $\bm{\Phi}$ AT2 phase-field matrix reads:
\begin{equation} 
\label{eq:dissipation_matrix}
	\begin{aligned}
		(\bm{\Phi}^{\text{II,IV}})_{\mathbf{ij}} &:= \int_{\Omega} \left( \frac{1}{l_0} R_{\mathbf{i,p}}\,R_{\mathbf{j,p}} \, + 2{l_0}\gamma \, \,\bm{\text{B}}_{\text{i,p}} \, \bm{\text{B}}_{\text{j,p}} + 2l_0^3\delta \, \bm{\text{C}}_{\text{i,p}}\bm{\text{C}}_{\text{j,p}} \right) \text{d} \Omega\,,
	\end{aligned}
\end{equation}
where, if $\delta = 0$, the second-order model is recovered. Otherwise, if $\beta = 1$, the AT1 phase-field matrix is:
\begin{equation} 
\label{eq:dissipation_matrix_AT1}
	\begin{aligned}
		(\bm{\Phi}^{\text{II,IV}})_{\mathbf{ij}} &:= \int_{\Omega} \left( 2{l_0}\gamma \, \,\bm{\text{B}}_{\text{i,p}} \, \bm{\text{B}}_{\text{j,p}} + 2l_0^3\delta \, \bm{\text{C}}_{\text{i,p}}\bm{\text{C}}_{\text{j,p}} \right) \text{d} \Omega\,.
	\end{aligned}
\end{equation}
For this case, the dissipated phase-field vector is defined as:
\begin{equation} \label{eq:dissipation_vector_AT1}
	\begin{aligned}
		(\bm{\phi}^{\text{II,IV}})_{\mathbf{i}} &:= \int_{\Omega} \frac{R_{\mathbf{i,p}}}{l_0}   \text{d} \Omega\,.
	\end{aligned}
\end{equation}
Instead, the stiffness matrix is defined as:
\begin{equation}\label{eq: stiffness matrix}
    (\bm{{K}})_{\mathbf{ij}} = \int_{\Omega} \qty( [(1 - {R}_{\bm{i,p}}\ \mathbf{\hat{d}}_e)^2(\psi^+_v + \psi_d) + \psi_d^-] {R}_{\bm{i,p}} {R}_{\bm{j,p}}) \, \text{d} \Omega\,,
\end{equation}
and the external force vector as:
\begin{equation}\label{eq: external force vector}
	(\VEC{F}^{\text{ext}})_{\mathbf{i}} = \int_{\Omega} \mathbf{b} \, R_{\bm{i,p}} \, \text{d} \Omega + \int_{\partial \Omega_{\text{N}}} \mathbf{t} \, R_{\bm{i,p}} \, \text{d} \partial \Omega\,,
\end{equation}
where $\VEC{b}$ and $\VEC{t}$ are the body and traction forces, respectively. For more details regarding the definitions of the matrices we refer the interested reader to \cite{greco2024higher, greco_at1}. 

\subsection{Numerical solution scheme} 
\label{subsec:PF_solution_scheme}

In order to solve \cref{eq:evol-sint}, we employ a staggered solution scheme, which alternates between solving the mechanical equilibrium and the phase-field evolution subproblems. Provided the solutions $\VEC{u}_{n}$ and $\SCALAR{d}_{n}$ at load step $n$, and defining the solutions at staggered iteration $i$ and load step $n+1$ as $\VEC{u}^i_{n+1}$ and $\SCALAR{d}^i_{n+1}$, the staggered scheme iteratively solves the mechanical and phase-field problems until convergence is achieved. Firstly, the mechanical subproblem from \cref{eq:matrix_problem_u} is solved using Picard iterations denoted by index $j$:
\begin{equation}
	\MAT{K}( \VEC{u}^{i,j-1}_{n+1} , \SCALAR{d}^{i-1}_{n} ) \,  \VEC{u}^{i,j}_{n+1} = \VEC{F}_{n}^{\text{ext}}\,.
\end{equation}
The Picard iterations are converged if the residual $\text{Res}^{j}_{\text{Pic},\VEC{u}} = \Vert \MAT{K}( \VEC{u}^{i,j}_{n+1} , \SCALAR{d}^{i-1}_{n} ) \,  \VEC{u}^{i,j}_{n+1} - \VEC{F}_{n}^{\text{ext}} \Vert_{L^2}$ is below a specified tolerance $\mathtt{TOL}_{\text{Pic},\VEC{u}}$. Once converged, the displacement solution is updated as $\VEC{u}^i_{n+1} = \VEC{u}^{i,j}_{n+1}$ and the phase-field subproblems from \eqref{eq:matrix_problem_u} are solved using the PSOR algorithm with matrix $\MAT{Q}^{\text{II,IV}} ( \VEC{u}^{i}_{n+1} )$ and vector $\VEC{R} = \MAT{Q}^{\text{II,IV}} ( \VEC{u}^{i}_{n+1} , \SCALAR{d}_{n}) \SCALAR{d}_{n} - \VEC{\psi} (\VEC{u}^{i}_{n+1}) + \VEC{\phi}^{\text{II,IV}}$ using a tolerance $\mathtt{TOL}_{\text{PSOR},\Delta\SCALAR{d}}$ (see \cite{MARENGO2021114137, greco2024higher} for more details). The PSOR solver provides the solution increment $\Delta \SCALAR{d}^i_{n+1}$, which is then used to update the phase-field solution as $\SCALAR{d}^i_{n+1} = \SCALAR{d}_{n} + \Delta \SCALAR{d}^i_{n+1}$. The staggered iterations are repeated until the residual $\text{Res}_{\text{stag}} = \Vert \MAT{K}( \VEC{u}^{i}_{n+1} , \SCALAR{d}^{i}_{n+1} ) \,  \VEC{u}^{i}_{n+1} - \VEC{F}_{n}^{\text{ext}} \Vert_{L^2}$ is below a specified tolerance $\mathtt{TOL}_{\text{stag}}$. For an algorithmic summary of the staggered solution scheme, we refer to \cref{alg:solve_elasticity,alg:solve_phasefield,alg:load_step} in \ref{app:algorithms}.

\subsection{Phase-field initialization}
\label{subsec:PF_initialization}
From physics perspective, it is well known that cracks initiate where stress concentrations are highest \cite{Griffith1921, irwin1957analysis}. Commonly, stress concentrations appear at geometric discontinuities, for example in notches, pores, sharp corners, or at places with pre-existing damage \cite{BORDEN201277, greco2024higher}. Alternatively, cracks can initiated at material discontinuities, as shown in the works \cite{GERASIMOV2019990, miehe2010phase, sargado2018high}.\\

In phase-field fracture simulations, cracks can be instantiated through geometric representation of geometric discontinuities or by providing an initial representation of the damage field, simulating the crack by degrading the initial material locally.\\

Firstly, geometric damage instantiation defines initial cracks or pores as geometric voids. This can be done by meshing of the initial geometry around these geometric discontinuities, as commonly done in FEM \cite{GERASIMOV2019990}, or by geometric methods such as trimming or imersion to locally ``disable'' the geometric domain. In the isogeometric analysis framework, which is the scope of this paper, the former approach typically requires multi-patch simulation, which becomes non-trivial for higher-order basis functions \cite{Verhelst2024a}. The latter approach, on the other hand, provides geometric generality, but the price is paid in terms of complicated quadrature schemes and potential advanced preconditioning \cite{dePrenter2023}.\\

Secondly, damage instantiation by defining an initial phase-field is a common problem in phase-field fracture simulation. In this case, initial cracks or pores are represented by damaged material rather than material voids, which is done by initializing the phase-field. To this end, the phase-field can be initialized by imposition of the history variable \cite{BORDEN201277}, or by an interpolated phase-field variable (IPF, \cite{greco2024higher}).\\
{\NEW In this paper, we adopt the IPF approach. In this approach, a pre-field of constant value $0.9999$ is defined within a $\beta$-neighborhood of the initial crack, where $\beta=\frac{1}{3}\ell_0$ for mesh size $h=\frac{\ell_0}{4}$, and $\beta=\frac{3}{4}\ell_0$ for mesh size $h=\frac{\ell_0}{2}$, as suggested in \cite{greco2024higher}. The phase-field is then initialized by performing an $L_2$ projection of this neighborhood onto the phase-field basis functions. Since the neighborhood is defined as a band around the initial crack, the $L_2$ projection can be solved by solving a low-rank subproblem containing only the basis functions with support in the neighborhood, which is computationally more efficient. This approach is general for both the AT1 and AT2 models, as it does not rely on the definition of a history variable, which is specific to the AT2 model.}

Ultimately, {\NEW the geometric damage instantiation and the phase-field initialization approaches} can be combined. For example, porous media can be represented geometrically, while cracks in this media can be represented by an initial phase-field. In this paper, we focus on the phase-field initialization approach, leaving the geometric or combined representations of discontinuities in the presented framework future research.\\

\section{Adaptive refinement using THB-splines}
\label{sec:THB}
This section provides preliminary information regarding THB-splines and their refinement. The section relies on previous works on Hierarchical B-splines \cite{Vuong2011}, Truncated Hierarchical B-splines \cite{Giannelli2012,Giannelli2016}, and admissible refinement \cite{Buffa2016,Bracco2018}. While the referenced works primarily form the mathematical foundation behind (T)HB-splines and their (admissible) refinement, the present section aims to provide an explanation of the THB-spline definition and the admissible meshing strategies by means of examples. For fundamental details about the material in this section, the reader is referred to the aforementioned references. The section is outlined as follows: \cref{subsec:THBsplines} elaborates on the construction of a THB-spline basis, and \cref{subsec:admissible} elaborates on the concept of mesh admissibility.\\

\subsection{Truncated Hierarchical B-splines}
\label{subsec:THBsplines}
Among other spline constructions, THB-splines provide a locally refinable spline basis. Compared to Hierarchical B-splines, THB-splines form a partition of unity, making them particularly appealing for solving partial differential equations using IGA. While we refer for the mathematical details behind THB-splines to the work of \cite{Vuong2011,Giannelli2012}, this section is limited to a brief definition of THB splines, and instead aims to provide the conceptual idea behind the construction of THB-spline bases.\\

Since THB-splines are a truncated version of Hierarchical B-splines (HB-splines), we first define the latter. Provided a sequence of $N$ nested tensor B-spline spaces in different levels $l=0,...,N-1$, denoted by $V^0\subset V^1\subset,...,V^{N-1}$ with an associated basis $\mathcal{B}^\ell$ of degree $p$. Using a sequence of parametric domains defined as $\Omega=\Omega^0\supseteq\Omega^1\supseteq...\supseteq\Omega^{N-1}=\emptyset$, the set of Hierarchical B-spline basis functions is defined as follows:
\begin{equation}
\mathcal{H} = \qty{ \beta \in \mathcal{B}^\ell \: : \: \text{supp}\qty(\beta) \subseteq \Omega^\ell \wedge  \text{supp}\qty(\beta) \not\subseteq \Omega^{\ell+1},\: \forall \ell=0,\dots,N}.
\end{equation}
In other words, the active functions of level $\ell$ in $\mathcal{H}$ are the basis functions of $\mathcal{B}^{\ell}$ which are fully contained in the parametric domain $\Omega^{\ell}$ and which are not fully contained in level $\ell+1$.\\

Similar to HB-splines, THB-splines are defined by a selection mechanism based on a sequence of nested parametric domains $\Omega^k$. However, the functions in the THB-spline basis are truncated using a truncation operation. Provided any function $\tau\in V^\ell$, it can be represented in the finer basis $V^{\ell+1}$ by taking the linear combination of a set of coefficients $c^{\ell+1}_\beta(\tau)\in\mathbb{R}$ and the basis functions $\beta\in\mathcal{B}^{\ell+1}$, i.e.,
\begin{equation}
\text{repr}^{\ell+1}\qty(\tau) = \sum_{\beta\in\mathcal{B}^{\ell+1}} c^{\ell+1}_\beta(\tau) \beta.
\end{equation}
In case of truncation of a function, this representation is restricted to the basis functions $\beta\in\mathcal{B}^{\ell+1}$ which are not fully contained in $\Omega^{\ell+1}$, i.e.
\begin{equation}
\text{trunc}^{\ell+1}\qty(\tau) = \sum_{\beta\in\mathcal{B}^{\ell+1},\:\text{supp}\qty(\beta)\not\subseteq \Omega^{\ell+1}} c^{\ell+1}_\beta(\tau) \beta.
\end{equation}
Since the truncation is to be applied recursively on all levels overlapping with the support of the considered basis function, a recursive definition of the THB-spline basis is more common:
\begin{enumerate}
	\item Initialize $\mathcal{T}^0 = \qty{\tau\in\mathcal{B}^0\: : \: \text{supp}\qty(\beta) \neq \emptyset}$
	\item Recursively, let $\mathcal{T}^{\ell+1} = \mathcal{T}_A^{\ell+1} + \mathcal{T}_B^{\ell+1}$ for $\ell=0,\dots,N-2$, where
	\begin{align*}
	\mathcal{T}_A^{\ell+1} &= \qty{\text{trunc}^{\ell+1}\qty(\tau)\: : \: \tau \in \mathcal{T}^\ell \wedge \text{supp}\qty()\tau \not\subseteq\Omega^{\ell+1}},\\
	\mathcal{T}_B^{\ell+1} &= \qty{\tau\in\beta^{\ell+1}\: : \: \text{supp}\qty(\tau)\subseteq\mathcal{B}^{\ell+1}}.
	\end{align*}
	\item Finally, $\mathcal{T} = \mathcal{T}^{N-1}$.
\end{enumerate}

\begin{example}[THB-spline representation and truncation]
Consider a knot vector with knots $\Xi^0 = \qty{0,0,0,\frac{1}{8},\frac{1}{4},\frac{3}{8},\frac{1}{2},\frac{5}{8},\frac{3}{4},1,1,1}$, hence defining a B-spline basis $\mathcal{B}^0$ of degree 2. Consequently, let the level $\mathcal{B}^1$ be defined using the knot vector $\Xi^1=\qty{0,0,0,\frac{1}{16},\dots,1,1,1}$. Additionally, let us refine the interval $[\frac{3}{8},\frac{7}{8}] = \Omega^1$, while $\Omega^0=\Omega=[0,1]$. The final THB-spline basis defined in this setting is given on the top of figure \cref{fig:IGArefinement}, where the black functions are B-spline functions $\beta\in\mathcal{B}^0$ for which $\text{supp}\qty(\beta)\cap\Omega^1=\emptyset$, the yellow and green/blue functions are, respectively, non-truncated and truncated functions $\beta\in\mathcal{B}^0$ for which $\text{supp}\qty(\beta)\cap\Omega^1\neq\emptyset$ and $\text{supp}\qty(\beta)\not\subseteq\Omega^1$ and the orange functions $\beta\in\mathcal{B}^1$ are from the fine level, hence satisfying $\text{supp}\qty(\beta)\subseteq\Omega^1$.\\

In the second up to the fourth rows of \cref{fig:IGArefinement}, the truncation mechanism is illustrated for three different functions $\beta\in\mathcal{B}^0$. Firstly, the second row of \cref{fig:IGArefinement} shows three different functions $\beta\in\mathcal{B}^0$. The green and blue function are active since $\text{supp}\qty(\beta)\not\subseteq\Omega^1$, whereas the orange function satisfies $\text{supp}\qty(\beta)\subseteq\Omega^1$, hence is not activated. Nevertheless, it will be shown that its truncation $\text{trunc}\qty(\beta)$ would yield zero coefficients. In the third row of \cref{fig:IGArefinement}, $\text{repr}^{1}\qty(\beta)$ is plotted, with the functions $\gamma\in\mathcal{B}^{1}$ in the background. Finally, the bottom row of \cref{fig:IGArefinement} shows $\text{trunc}^{1}\qty(\beta)$ for each of the three functions, with the functions $\gamma\in\mathcal{B}^{1}$ in the background, and with the functions from $\text{repr}^{1}\qty(\beta)$ which are not fully contained in $\Omega^1$, i.e., the functions contributing to $\text{trunc}^{1}\qty(\beta)$. For the orange function, which satisfies $\beta\subseteq\Omega^1$, it can be seen that all representing functions from $\mathcal{B}^1$ are fully contained in $\Omega^1$, hence $\text{trunc}^{1}\qty(\beta)=0$.
\end{example}



\begin{figure}
	\centering
	\begin{tikzpicture}
		\begin{groupplot}
			[
			height=0.15\textheight,
			width=0.4\linewidth,
			ytick = {0,1},
			legend pos = outer north east,
			group style={
				group name=my plots,
				group size=3 by 4,
				xlabels at=edge bottom,
				x descriptions at=edge bottom,
				y descriptions at=edge left,
				vertical sep=5pt,
				horizontal sep=5pt},
			xlabel={$\xi$},
			xtick={0.00,0.125,0.250,0.375,0.4375,0.500,0.5625,0.625,0.6875,0.750,0.875,1.000},
			xticklabels={0,$\frac{1}{8}$,$\frac{1}{4}$,$\frac{3}{8}$,$\frac{7}{16}$,$\frac{1}{2}$,$\frac{9}{16}$,$\frac{5}{8}$,$\frac{11}{16}$,$\frac{3}{4}$,$\frac{7}{8}$,1},
   			grid style={solid, gray!40},
			]
			\nextgroupplot[group/empty plot]

			\nextgroupplot[no markers,ylabel={$\beta^{THB}_i$}]
			\begin{pgfonlayer}{background}
				\fill[black!05] (axis cs: 0.000,0.0) rectangle (axis cs: 0.375,1.0);
				\fill[black!30] (axis cs: 0.375,0.0) rectangle (axis cs: 0.750,1.0);
				\fill[black!05] (axis cs: 0.750,0.0) rectangle (axis cs: 1.000,1.0);
			\end{pgfonlayer}

			\foreach \k in {0,1,2,3,6,7,8}
			{
				\pgfmathtruncatemacro{\y}{2*\k+1};
				\pgfmathtruncatemacro{\x}{2*\k};
				\addplot+[black,thin,solid] table[col sep = comma, header=false, y index = {\y}, x index = {\x}] {Data/Refinement_illustration/THBspline.csv};
			}

			\foreach \k in {3,6}
			{
				\pgfmathtruncatemacro{\y}{2*\k+1};
				\pgfmathtruncatemacro{\x}{2*\k};
				\addplot+[col4,thick,solid] table[col sep = comma, header=false, y index = {\y}, x index = {\x}] {Data/Refinement_illustration/THBspline.csv};
			}

			\def\k{4}
			\pgfmathtruncatemacro{\y}{2*\k+1};
			\pgfmathtruncatemacro{\x}{2*\k};
			\addplot+[col1,thick,solid] table[col sep = comma, header=false, y index = {\y}, x index = {\x}] {Data/Refinement_illustration/THBspline.csv};

			\def\k{5}
			\pgfmathtruncatemacro{\y}{2*\k+1};
			\pgfmathtruncatemacro{\x}{2*\k};
			\addplot+[col3,thick,solid] table[col sep = comma, header=false, y index = {\y}, x index = {\x}] {Data/Refinement_illustration/THBspline.csv};

			\foreach \k in {9,...,12}
			{
				\pgfmathtruncatemacro{\y}{2*\k+1};
				\pgfmathtruncatemacro{\x}{2*\k};
				\addplot+[col2,thick,solid] table[col sep = comma, header=false, y index = {\y}, x index = {\x}] {Data/Refinement_illustration/THBspline.csv};
			}

			\nextgroupplot[group/empty plot]

			\nextgroupplot[no markers,ylabel={$\beta^0$}]
			\begin{pgfonlayer}{background}
			\fill[black!05] (axis cs: 0.000,0.0) rectangle (axis cs: 0.375,1.0);
			\fill[black!30] (axis cs: 0.375,0.0) rectangle (axis cs: 0.750,1.0);
			\fill[black!05] (axis cs: 0.750,0.0) rectangle (axis cs: 1.000,1.0);
			\end{pgfonlayer}

			\foreach \k in {0,...,3}
			{
				\pgfmathtruncatemacro{\y}{2*\k+1};
				\pgfmathtruncatemacro{\x}{2*\k};
				\addplot[black!40,thin,solid] table[col sep = comma, header=false, y index = {\y}, x index = {\x}] {Data/Refinement_illustration/init.csv};
			}
			\foreach \k in {5,...,9}
			{
				\pgfmathtruncatemacro{\y}{2*\k+1};
				\pgfmathtruncatemacro{\x}{2*\k};
				\addplot+[black!40,thin,solid] table[col sep = comma, header=false, y index = {\y}, x index = {\x}] {Data/Refinement_illustration/init.csv};
			}
			\def\k{4}
			\pgfmathtruncatemacro{\y}{2*\k+1};
			\pgfmathtruncatemacro{\x}{2*\k};
			\addplot+[col1,thick,solid] table[col sep = comma, header=false, y index = {\y}, x index = {\x}] {Data/Refinement_illustration/init.csv};

			\nextgroupplot[no markers]
			\begin{pgfonlayer}{background}
				\fill[black!05] (axis cs: 0.000,0.0) rectangle (axis cs: 0.375,1.0);
				\fill[black!30] (axis cs: 0.375,0.0) rectangle (axis cs: 0.750,1.0);
				\fill[black!05] (axis cs: 0.750,0.0) rectangle (axis cs: 1.000,1.0);
			\end{pgfonlayer}
			\foreach \k in {0,...,4}
			{
				\pgfmathtruncatemacro{\y}{2*\k+1};
				\pgfmathtruncatemacro{\x}{2*\k};
				\addplot+[black!40,thin,solid] table[col sep = comma, header=false, y index = {\y}, x index = {\x}] {Data/Refinement_illustration/init.csv};
			}
			\foreach \k in {6,...,9}
			{
				\pgfmathtruncatemacro{\y}{2*\k+1};
				\pgfmathtruncatemacro{\x}{2*\k};
				\addplot+[black!40,thin,solid] table[col sep = comma, header=false, y index = {\y}, x index = {\x}] {Data/Refinement_illustration/init.csv};
			}

			\def\k{5}
			\pgfmathtruncatemacro{\y}{2*\k+1};
			\pgfmathtruncatemacro{\x}{2*\k};
			\addplot+[col2,thick,solid] table[col sep = comma, header=false, y index = {\y}, x index = {\x}] {Data/Refinement_illustration/init.csv};

			\nextgroupplot[no markers]
			\begin{pgfonlayer}{background}
				\fill[black!05] (axis cs: 0.000,0.0) rectangle (axis cs: 0.375,1.0);
				\fill[black!30] (axis cs: 0.375,0.0) rectangle (axis cs: 0.750,1.0);
				\fill[black!05] (axis cs: 0.750,0.0) rectangle (axis cs: 1.000,1.0);
			\end{pgfonlayer}
			\foreach \k in {0,...,5}
			{
				\pgfmathtruncatemacro{\y}{2*\k+1};
				\pgfmathtruncatemacro{\x}{2*\k};
				\addplot+[black!40,thin,solid] table[col sep = comma, header=false, y index = {\y}, x index = {\x}] {Data/Refinement_illustration/init.csv};
			}
			\foreach \k in {7,...,9}
			{
				\pgfmathtruncatemacro{\y}{2*\k+1};
				\pgfmathtruncatemacro{\x}{2*\k};
				\addplot+[black!40,thin,solid] table[col sep = comma, header=false, y index = {\y}, x index = {\x}] {Data/Refinement_illustration/init.csv};
			}

			\def\k{6}
			\pgfmathtruncatemacro{\y}{2*\k+1};
			\pgfmathtruncatemacro{\x}{2*\k};
			\addplot+[col3,thick,solid] table[col sep = comma, header=false, y index = {\y}, x index = {\x}] {Data/Refinement_illustration/init.csv};

			\nextgroupplot[no markers,ylabel={$\text{repr}^{1}\qty(\beta^0)$}]
			\begin{pgfonlayer}{background}
				\fill[black!05] (axis cs: 0.000,0.0) rectangle (axis cs: 0.375,1.0);
				\fill[black!30] (axis cs: 0.375,0.0) rectangle (axis cs: 0.750,1.0);
				\fill[black!05] (axis cs: 0.750,0.0) rectangle (axis cs: 1.000,1.0);
			\end{pgfonlayer}
			\foreach \k in {0,...,17}
			{
				\pgfmathtruncatemacro{\y}{2*\k+1};
				\pgfmathtruncatemacro{\x}{2*\k};
				\addplot+[black!40,thin,solid] table[col sep = comma, header=false, y index = {\y}, x index = {\x}] {Data/Refinement_illustration/BSpline.csv};
			}
			\def\k{4}
			\pgfmathtruncatemacro{\y}{2*\k+1};
			\pgfmathtruncatemacro{\x}{2*\k};
			\addplot+[col1,thick,dotted] table[col sep = comma, header=false, y index = {\y}, x index = {\x}] {Data/Refinement_illustration/init.csv};

			\def\k{6}
			\pgfmathtruncatemacro{\y}{2*\k+1};
			\pgfmathtruncatemacro{\x}{2*\k};
			\addplot+[col1,thick,solid] table[col sep = comma, header=false, y expr = 0.25*\thisrowno{\y}, x index = {\x}] {Data/Refinement_illustration/BSpline.csv};
			\def\k{7}
			\pgfmathtruncatemacro{\y}{2*\k+1};
			\pgfmathtruncatemacro{\x}{2*\k};
			\addplot+[col1,thick,solid] table[col sep = comma, header=false, y expr = 0.75*\thisrowno{\y}, x index = {\x}] {Data/Refinement_illustration/BSpline.csv};
			\def\k{8}
			\pgfmathtruncatemacro{\y}{2*\k+1};
			\pgfmathtruncatemacro{\x}{2*\k};
			\addplot+[col1,thick,solid] table[col sep = comma, header=false, y expr = 0.75*\thisrowno{\y}, x index = {\x}] {Data/Refinement_illustration/BSpline.csv};
			\def\k{9}
			\pgfmathtruncatemacro{\y}{2*\k+1};
			\pgfmathtruncatemacro{\x}{2*\k};
			\addplot+[col1,thick,solid] table[col sep = comma, header=false, y expr = 0.25*\thisrowno{\y}, x index = {\x}] {Data/Refinement_illustration/BSpline.csv};

			\nextgroupplot[no markers]
			\begin{pgfonlayer}{background}
				\fill[black!05] (axis cs: 0.000,0.0) rectangle (axis cs: 0.375,1.0);
				\fill[black!30] (axis cs: 0.375,0.0) rectangle (axis cs: 0.750,1.0);
				\fill[black!05] (axis cs: 0.750,0.0) rectangle (axis cs: 1.000,1.0);
			\end{pgfonlayer}
			\foreach \k in {0,...,17}
			{
				\pgfmathtruncatemacro{\y}{2*\k+1};
				\pgfmathtruncatemacro{\x}{2*\k};
				\addplot+[black!40,thin,solid] table[col sep = comma, header=false, y index = {\y}, x index = {\x}] {Data/Refinement_illustration/BSpline.csv};
			}
			\def\k{5}
			\pgfmathtruncatemacro{\y}{2*\k+1};
			\pgfmathtruncatemacro{\x}{2*\k};
			\addplot+[col2,thick,dotted] table[col sep = comma, header=false, y index = {\y}, x index = {\x}] {Data/Refinement_illustration/init.csv};

			\def\k{8}
			\pgfmathtruncatemacro{\y}{2*\k+1};
			\pgfmathtruncatemacro{\x}{2*\k};
			\addplot+[col2,thick,solid] table[col sep = comma, header=false, y expr = 0.25*\thisrowno{\y}, x index = {\x}] {Data/Refinement_illustration/BSpline.csv};
			\def\k{9}
			\pgfmathtruncatemacro{\y}{2*\k+1};
			\pgfmathtruncatemacro{\x}{2*\k};
			\addplot+[col2,thick,solid] table[col sep = comma, header=false, y expr = 0.75*\thisrowno{\y}, x index = {\x}] {Data/Refinement_illustration/BSpline.csv};
			\def\k{10}
			\pgfmathtruncatemacro{\y}{2*\k+1};
			\pgfmathtruncatemacro{\x}{2*\k};
			\addplot+[col2,thick,solid] table[col sep = comma, header=false, y expr = 0.75*\thisrowno{\y}, x index = {\x}] {Data/Refinement_illustration/BSpline.csv};
			\def\k{11}
			\pgfmathtruncatemacro{\y}{2*\k+1};
			\pgfmathtruncatemacro{\x}{2*\k};
			\addplot+[col2,thick,solid] table[col sep = comma, header=false, y expr = 0.25*\thisrowno{\y}, x index = {\x}] {Data/Refinement_illustration/BSpline.csv};

			\nextgroupplot[no markers]
			\begin{pgfonlayer}{background}
				\fill[black!05] (axis cs: 0.000,0.0) rectangle (axis cs: 0.375,1.0);
				\fill[black!30] (axis cs: 0.375,0.0) rectangle (axis cs: 0.750,1.0);
				\fill[black!05] (axis cs: 0.750,0.0) rectangle (axis cs: 1.000,1.0);
			\end{pgfonlayer}
			\foreach \k in {0,...,17}
			{
				\pgfmathtruncatemacro{\y}{2*\k+1};
				\pgfmathtruncatemacro{\x}{2*\k};
				\addplot+[black!40,thin,solid] table[col sep = comma, header=false, y index = {\y}, x index = {\x}] {Data/Refinement_illustration/BSpline.csv};
			}
			\def\k{6}
			\pgfmathtruncatemacro{\y}{2*\k+1};
			\pgfmathtruncatemacro{\x}{2*\k};
			\addplot+[col3,thick,dotted] table[col sep = comma, header=false, y index = {\y}, x index = {\x}] {Data/Refinement_illustration/init.csv};

			\def\k{10}
			\pgfmathtruncatemacro{\y}{2*\k+1};
			\pgfmathtruncatemacro{\x}{2*\k};
			\addplot+[col3,thick,solid] table[col sep = comma, header=false, y expr = 0.25*\thisrowno{\y}, x index = {\x}] {Data/Refinement_illustration/BSpline.csv};
			\def\k{11}
			\pgfmathtruncatemacro{\y}{2*\k+1};
			\pgfmathtruncatemacro{\x}{2*\k};
			\addplot+[col3,thick,solid] table[col sep = comma, header=false, y expr = 0.75*\thisrowno{\y}, x index = {\x}] {Data/Refinement_illustration/BSpline.csv};
			\def\k{12}
			\pgfmathtruncatemacro{\y}{2*\k+1};
			\pgfmathtruncatemacro{\x}{2*\k};
			\addplot+[col3,thick,solid] table[col sep = comma, header=false, y expr = 0.75*\thisrowno{\y}, x index = {\x}] {Data/Refinement_illustration/BSpline.csv};
			\def\k{13}
			\pgfmathtruncatemacro{\y}{2*\k+1};
			\pgfmathtruncatemacro{\x}{2*\k};
			\addplot+[col3,thick,solid] table[col sep = comma, header=false, y expr = 0.25*\thisrowno{\y}, x index = {\x}] {Data/Refinement_illustration/BSpline.csv};

			\nextgroupplot[no markers,ylabel={$\text{trunc}^{1}\qty(\tau)$}]
			\begin{pgfonlayer}{background}
				\fill[black!05] (axis cs: 0.000,0.0) rectangle (axis cs: 0.375,1.0);
				\fill[black!30] (axis cs: 0.375,0.0) rectangle (axis cs: 0.750,1.0);
				\fill[black!05] (axis cs: 0.750,0.0) rectangle (axis cs: 1.000,1.0);
			\end{pgfonlayer}
			\foreach \k in {0,...,17}
			{
				\pgfmathtruncatemacro{\y}{2*\k+1};
				\pgfmathtruncatemacro{\x}{2*\k};
				\addplot+[black!40,thin,solid] table[col sep = comma, header=false, y index = {\y}, x index = {\x}] {Data/Refinement_illustration/BSpline.csv};
			}
			\def\k{4}
			\pgfmathtruncatemacro{\y}{2*\k+1};
			\pgfmathtruncatemacro{\x}{2*\k};
			\addplot+[col1,thick,dotted] table[col sep = comma, header=false, y index = {\y}, x index = {\x}] {Data/Refinement_illustration/init.csv};

			\def\k{6}
			\pgfmathtruncatemacro{\y}{2*\k+1};
			\pgfmathtruncatemacro{\x}{2*\k};
			\addplot+[col1,thick,solid] table[col sep = comma, header=false, y expr = 0.25*\thisrowno{\y}, x index = {\x}] {Data/Refinement_illustration/BSpline.csv};
			\def\k{7}
			\pgfmathtruncatemacro{\y}{2*\k+1};
			\pgfmathtruncatemacro{\x}{2*\k};
			\addplot+[col1,thick,solid] table[col sep = comma, header=false, y expr = 0.75*\thisrowno{\y}, x index = {\x}] {Data/Refinement_illustration/BSpline.csv};
			\def\k{8}
			\pgfmathtruncatemacro{\y}{2*\k+1};
			\pgfmathtruncatemacro{\x}{2*\k};
			\addplot+[col1,thick,solid] table[col sep = comma, header=false, y expr = 0.00*\thisrowno{\y}, x index = {\x}] {Data/Refinement_illustration/BSpline.csv};
			\def\k{9}
			\pgfmathtruncatemacro{\y}{2*\k+1};
			\pgfmathtruncatemacro{\x}{2*\k};
			\addplot+[col1,thick,solid] table[col sep = comma, header=false, y expr = 0.00*\thisrowno{\y}, x index = {\x}] {Data/Refinement_illustration/BSpline.csv};

			\def\k{4}
			\pgfmathtruncatemacro{\y}{2*\k+1};
			\pgfmathtruncatemacro{\x}{2*\k};
			\addplot+[col1,thick,solid] table[col sep = comma, header=false, y index = {\y}, x index = {\x}] {Data/Refinement_illustration/THBspline.csv};

			\nextgroupplot[no markers]
			\begin{pgfonlayer}{background}
				\fill[black!05] (axis cs: 0.000,0.0) rectangle (axis cs: 0.375,1.0);
				\fill[black!30] (axis cs: 0.375,0.0) rectangle (axis cs: 0.750,1.0);
				\fill[black!05] (axis cs: 0.750,0.0) rectangle (axis cs: 1.000,1.0);
			\end{pgfonlayer}
			\foreach \k in {0,...,7}
			{
				\pgfmathtruncatemacro{\y}{2*\k+1};
				\pgfmathtruncatemacro{\x}{2*\k};
				\addplot+[black!40,thin,solid] table[col sep = comma, header=false, y index = {\y}, x index = {\x}] {Data/Refinement_illustration/BSpline.csv};
			}
			\foreach \k in {8,...,11}
			{
				\pgfmathtruncatemacro{\y}{2*\k+1};
				\pgfmathtruncatemacro{\x}{2*\k};
				\addplot+[black!40,thin,solid] table[col sep = comma, header=false, y index = {\y}, x index = {\x}] {Data/Refinement_illustration/BSpline.csv};
			}
			\foreach \k in {12,...,17}
			{
				\pgfmathtruncatemacro{\y}{2*\k+1};
				\pgfmathtruncatemacro{\x}{2*\k};
				\addplot+[black!40,thin,solid] table[col sep = comma, header=false, y index = {\y}, x index = {\x}] {Data/Refinement_illustration/BSpline.csv};
			}
			\def\k{5}
			\pgfmathtruncatemacro{\y}{2*\k+1};
			\pgfmathtruncatemacro{\x}{2*\k};
			\addplot+[col2,thick,dotted] table[col sep = comma, header=false, y index = {\y}, x index = {\x}] {Data/Refinement_illustration/init.csv};

			\def\k{8}
			\pgfmathtruncatemacro{\y}{2*\k+1};
			\pgfmathtruncatemacro{\x}{2*\k};
			\addplot+[col2,thick,solid] table[col sep = comma, header=false, y expr = 0.00*\thisrowno{\y}, x index = {\x}] {Data/Refinement_illustration/BSpline.csv};
			\def\k{9}
			\pgfmathtruncatemacro{\y}{2*\k+1};
			\pgfmathtruncatemacro{\x}{2*\k};
			\addplot+[col2,thick,solid] table[col sep = comma, header=false, y expr = 0.00*\thisrowno{\y}, x index = {\x}] {Data/Refinement_illustration/BSpline.csv};
			\def\k{10}
			\pgfmathtruncatemacro{\y}{2*\k+1};
			\pgfmathtruncatemacro{\x}{2*\k};
			\addplot+[col2,thick,solid] table[col sep = comma, header=false, y expr = 0.00*\thisrowno{\y}, x index = {\x}] {Data/Refinement_illustration/BSpline.csv};
			\def\k{11}
			\pgfmathtruncatemacro{\y}{2*\k+1};
			\pgfmathtruncatemacro{\x}{2*\k};
			\addplot+[col2,thick,solid] table[col sep = comma, header=false, y expr = 0.00*\thisrowno{\y}, x index = {\x}] {Data/Refinement_illustration/BSpline.csv};

			\nextgroupplot[no markers]
			\begin{pgfonlayer}{background}
				\fill[black!05] (axis cs: 0.000,0.0) rectangle (axis cs: 0.375,1.0);
				\fill[black!30] (axis cs: 0.375,0.0) rectangle (axis cs: 0.750,1.0);
				\fill[black!05] (axis cs: 0.750,0.0) rectangle (axis cs: 1.000,1.0);
			\end{pgfonlayer}
			\foreach \k in {0,...,17}
			{
				\pgfmathtruncatemacro{\y}{2*\k+1};
				\pgfmathtruncatemacro{\x}{2*\k};
				\addplot+[black!40,thin,solid] table[col sep = comma, header=false, y index = {\y}, x index = {\x}] {Data/Refinement_illustration/BSpline.csv};
			}
			\def\k{6}
			\pgfmathtruncatemacro{\y}{2*\k+1};
			\pgfmathtruncatemacro{\x}{2*\k};
			\addplot+[col3,thick,dotted] table[col sep = comma, header=false, y index = {\y}, x index = {\x}] {Data/Refinement_illustration/init.csv};

			\def\k{10}
			\pgfmathtruncatemacro{\y}{2*\k+1};
			\pgfmathtruncatemacro{\x}{2*\k};
			\addplot+[col3,thick,solid] table[col sep = comma, header=false, y expr = 0.00*\thisrowno{\y}, x index = {\x}] {Data/Refinement_illustration/BSpline.csv};
			\def\k{11}
			\pgfmathtruncatemacro{\y}{2*\k+1};
			\pgfmathtruncatemacro{\x}{2*\k};
			\addplot+[col3,thick,solid] table[col sep = comma, header=false, y expr = 0.00*\thisrowno{\y}, x index = {\x}] {Data/Refinement_illustration/BSpline.csv};
			\def\k{12}
			\pgfmathtruncatemacro{\y}{2*\k+1};
			\pgfmathtruncatemacro{\x}{2*\k};
			\addplot+[col3,thick,solid] table[col sep = comma, header=false, y expr = 0.75*\thisrowno{\y}, x index = {\x}] {Data/Refinement_illustration/BSpline.csv};
			\def\k{13}
			\pgfmathtruncatemacro{\y}{2*\k+1};
			\pgfmathtruncatemacro{\x}{2*\k};
			\addplot+[col3,thick,solid] table[col sep = comma, header=false, y expr = 0.25*\thisrowno{\y}, x index = {\x}] {Data/Refinement_illustration/BSpline.csv};

			\def\k{5}
			\pgfmathtruncatemacro{\y}{2*\k+1};
			\pgfmathtruncatemacro{\x}{2*\k};
			\addplot+[col3,thick,solid] table[col sep = comma, header=false, y index = {\y}, x index = {\x}] {Data/Refinement_illustration/THBspline.csv};
		\end{groupplot}
	\end{tikzpicture}
	\caption{The concept of Truncated Hierarchical B-spline refinement in one dimension. The THB-spline basis in the top row is locally refined on the interval $[\frac{3}{8},\frac{3}{4}]$ (shaded area) and the functions influenced by the local refinement are color-coded as follows: {\NEW (gold) non-truncated functions overlapping the refined area; (red and green) truncated function of the coarse level, (blue) activated functions of the fine level. The second up to the last rows decompose three functions of the coarse level which are affected by the refinement: left and right the truncated coarse functions and in the middle the eliminated coarse function. In the second row, the original functions of the coarse level are highlighted. In the third row, the representation of the function in the finer level is provided, with the original function represented by a dotted line. In the last row, the representation coefficients of the functions which are fully contained in the marked (shaded) interval are set to zero, yielding the truncated functions (light) in the left and right columns. Since the function in the middle row is fully eliminated, the middle plot in the bottom row provides the functions from the fine level inserted instead.}}
	\label{fig:IGArefinement}
\end{figure}

\subsection{Admissible refinement}
\label{subsec:admissible}
In the context of (T)HB-splines, the concept of admissible meshes was first presented by \citeauthor{Buffa2016} \cite{Buffa2016}. Admissible meshes are meshes where the number of elements acting in any mesh point is bounded and that the level of all active functions in an element is independent of the number of overall levels in the hierarchy. In the work \cite{Carraturo2019} among others, it was shown that inadmissible meshes may lead to oscillations in the numerical solution when solving PDEs with (T)HB-splines, while in case of admissible meshes oscillations are avoided.\\

The works of \citeauthor{Gantner2017} \cite{Gantner2017} and \citeauthor{Bracco2018} \cite{Bracco2018} provide algorithms for the construction of admissible meshes for HB- and THB- splines. In the present work, we employ the admissible refinement algorithms for THB-splines based on the $\mathcal{T}$-neigborhood \cite{Bracco2018}. In \cref{ex:admissible_meshing}, we provide a simple example of admissible refinement of class $2$, while we refer to \cite{Buffa2016,Gantner2017,Buffa2022,Gantner2017} for the mathematical details behind admissible meshing, and specifically to \cite{Bracco2018} to the refinement algorithms, which are later used in \ref{app:algorithms}.

\begin{example}[Admissible meshing]\label{ex:admissible_meshing}
We consider a THB-spline basis of degree $p=2$ with unique knot vectors $\Xi_u=\Xi_v=\qty{0,\frac{1}{8},\dots,\frac{7}{8},1}$. The initial configuration is given in \cref{fig:admissibility1} and consists of mesh elements of levels $1,2,3$, along with two elements of level 3 to be refined with admissibility of class $m=2$. The first step of the admissible meshing algorithm is to find the support extension of the marked elements, see \cref{fig:admissibility2}, which consists of the union of the supports of all functions of the tensor level $\ell=3$ active in the marked elements. Consequently, all elements of level $\ell-m=1$ intersecting these support extensions need to be refined in order to preserve admissibility of class $m=2$. The final result after refinement is depicted in \cref{fig:admissibility4}.

\begin{figure}
\centering
\begin{subfigure}[t]{0.225\linewidth}
\centering
\resizebox{\linewidth}{!}
{
\begin{tikzpicture}[scale=4]
	\draw[color=red1!10,line width=20](0,0.5)--(0.25,0.5);
	\draw[color=red1!10,line width=20,line cap=round](0.25,0.5)--(0.5,0.5);
	\draw[color=red1!30,line width=15](0,0.5)--(0.25,0.5);
	\draw[color=red1!30,line width=15,line cap=round](0.25,0.5)--(0.5,0.5);
	\draw[color=red1!50,line width=10](0,0.5)--(0.25,0.5);
	\draw[color=red1!50,line width=10,line cap=round](0.25,0.5)--(0.5,0.5);
	\draw[color=red1!70,line width=05](0,0.5)--(0.25,0.5);
	\draw[color=red1!70,line width=05,line cap=round](0.25,0.5)--(0.5,0.5);

	\foreach \x in {0,0.125,...,1.0}{\draw[opacity=0.5](\x,0)--(\x,1);}
	\foreach \y in {0,0.125,...,1.0}{\draw[opacity=0.5](0,\y)--(1,\y);}
	\foreach \x in {0.0625,0.1875,...,0.625}{\draw[opacity=0.5](\x,0.25)--(\x,0.75);}
	\foreach \y in {0.3125,0.4375,...,0.6875}{\draw[opacity=0.5](0,\y)--(0.625,\y);}
	\foreach \x in {0.03125,0.09375,...,0.5}{\draw[opacity=0.5](\x,0.4375)--(\x,0.5625);}
	\foreach \y in {0.46875,0.53125}{\draw[opacity=0.5](0,\y)--(0.5,\y);}

	\fill[col2,opacity=0.5](0.5,0.5) rectangle (0.5625,0.5625);
	\fill[col3,opacity=0.5](0.5,0.4375) rectangle (0.5625,0.5);
\end{tikzpicture}
}
\caption{}
\label{fig:admissibility1}
\end{subfigure}
\hfill
\begin{subfigure}[t]{0.225\linewidth}
\centering
\resizebox{\linewidth}{!}
{
\begin{tikzpicture}[scale=4]
	\draw[color=red1!10,line width=20](0,0.5)--(0.25,0.5);
	\draw[color=red1!10,line width=20,line cap=round](0.25,0.5)--(0.5,0.5);
	\draw[color=red1!30,line width=15](0,0.5)--(0.25,0.5);
	\draw[color=red1!30,line width=15,line cap=round](0.25,0.5)--(0.5,0.5);
	\draw[color=red1!50,line width=10](0,0.5)--(0.25,0.5);
	\draw[color=red1!50,line width=10,line cap=round](0.25,0.5)--(0.5,0.5);
	\draw[color=red1!70,line width=05](0,0.5)--(0.25,0.5);
	\draw[color=red1!70,line width=05,line cap=round](0.25,0.5)--(0.5,0.5);

	\foreach \x in {0,0.125,...,1.0}{\draw[opacity=0.5](\x,0)--(\x,1);}
	\foreach \y in {0,0.125,...,1.0}{\draw[opacity=0.5](0,\y)--(1,\y);}
	\foreach \x in {0.0625,0.1875,...,0.625}{\draw[opacity=0.5](\x,0.25)--(\x,0.75);}
	\foreach \y in {0.3125,0.4375,...,0.6875}{\draw[opacity=0.5](0,\y)--(0.625,\y);}
	\foreach \x in {0.03125,0.09375,...,0.5}{\draw[opacity=0.5](\x,0.4375)--(\x,0.5625);}
	\foreach \y in {0.46875,0.53125}{\draw[opacity=0.5](0,\y)--(0.5,\y);}

	\fill[col2,opacity=0.5](0.5,0.5) rectangle (0.5625,0.5625);
	\fill[col3,opacity=0.5](0.5,0.4375) rectangle (0.5625,0.5);
	\fill[thick,draw=col2,col2,opacity=0.2](0.375,0.375) rectangle (0.6875,0.6875);
	\fill[thick,draw=col3,col3,opacity=0.2](0.375,0.3125) rectangle (0.6875,0.625);
\end{tikzpicture}
}
\caption{}
\label{fig:admissibility2}
\end{subfigure}
\hfill
\begin{subfigure}[t]{0.225\linewidth}
\centering
\resizebox{\linewidth}{!}
{
\begin{tikzpicture}[scale=4]
	\draw[color=red1!10,line width=20](0,0.5)--(0.25,0.5);
	\draw[color=red1!10,line width=20,line cap=round](0.25,0.5)--(0.5,0.5);
	\draw[color=red1!30,line width=15](0,0.5)--(0.25,0.5);
	\draw[color=red1!30,line width=15,line cap=round](0.25,0.5)--(0.5,0.5);
	\draw[color=red1!50,line width=10](0,0.5)--(0.25,0.5);
	\draw[color=red1!50,line width=10,line cap=round](0.25,0.5)--(0.5,0.5);
	\draw[color=red1!70,line width=05](0,0.5)--(0.25,0.5);
	\draw[color=red1!70,line width=05,line cap=round](0.25,0.5)--(0.5,0.5);

	\foreach \x in {0,0.125,...,1.0}{\draw[opacity=0.5](\x,0)--(\x,1);}
	\foreach \y in {0,0.125,...,1.0}{\draw[opacity=0.5](0,\y)--(1,\y);}
	\foreach \x in {0.0625,0.1875,...,0.625}{\draw[opacity=0.5](\x,0.25)--(\x,0.75);}
	\foreach \y in {0.3125,0.4375,...,0.6875}{\draw[opacity=0.5](0,\y)--(0.625,\y);}
	\foreach \x in {0.03125,0.09375,...,0.5}{\draw[opacity=0.5](\x,0.4375)--(\x,0.5625);}
	\foreach \y in {0.46875,0.53125}{\draw[opacity=0.5](0,\y)--(0.5,\y);}

	\draw[thick,col2,opacity=0.2](0.375,0.375) rectangle (0.6875,0.6875);
	\draw[thick,col3,opacity=0.2](0.375,0.3125) rectangle (0.6875,0.625);
	\fill[thick,col2,opacity=0.5](0.5,0.5) rectangle (0.5625,0.5625);
	\fill[thick,col3,opacity=0.5](0.5,0.4375) rectangle (0.5625,0.5);
	\fill[col4,opacity=0.3](0.625,0.25) rectangle (0.75,0.375);
	\fill[col4,opacity=0.3](0.625,0.375) rectangle (0.75,0.5);
	\fill[col4,opacity=0.3](0.625,0.5) rectangle (0.75,0.625);
	\fill[col4,opacity=0.3](0.625,0.625) rectangle (0.75,0.75);
\end{tikzpicture}
}
\caption{}
\label{fig:admissibility3}
\end{subfigure}
\hfill
\begin{subfigure}[t]{0.225\linewidth}
\centering
\resizebox{\linewidth}{!}
{
\begin{tikzpicture}[scale=4]
	\draw[color=red1!10,line width=20](0,0.5)--(0.25,0.5);
	\draw[color=red1!10,line width=20,line cap=round](0.25,0.5)--(0.5,0.5);
	\draw[color=red1!30,line width=15](0,0.5)--(0.25,0.5);
	\draw[color=red1!30,line width=15,line cap=round](0.25,0.5)--(0.5,0.5);
	\draw[color=red1!50,line width=10](0,0.5)--(0.25,0.5);
	\draw[color=red1!50,line width=10,line cap=round](0.25,0.5)--(0.5,0.5);
	\draw[color=red1!70,line width=05](0,0.5)--(0.25,0.5);
	\draw[color=red1!70,line width=05,line cap=round](0.25,0.5)--(0.5,0.5);

	\foreach \x in {0,0.125,...,1.0}{\draw[opacity=0.5](\x,0)--(\x,1);}
	\foreach \y in {0,0.125,...,1.0}{\draw[opacity=0.5](0,\y)--(1,\y);}
	\foreach \x in {0.0625,0.1875,...,0.625}{\draw[opacity=0.5](\x,0.25)--(\x,0.75);}
	\foreach \y in {0.3125,0.4375,...,0.6875}{\draw[opacity=0.5](0,\y)--(0.625,\y);}
	\foreach \x in {0.03125,0.09375,...,0.5}{\draw[opacity=0.5](\x,0.4375)--(\x,0.5625);}
	\foreach \y in {0.46875,0.53125}{\draw[opacity=0.5](0,\y)--(0.5,\y);}

	\fill[col2,opacity=0.5](0.5,0.5) rectangle (0.5625,0.5625);
	\fill[col2,opacity=0.5](0.5,0.4375) rectangle (0.5625,0.5);
	\fill[col2,opacity=0.3](0.625,0.25) rectangle (0.75,0.375);
	\fill[col2,opacity=0.3](0.625,0.375) rectangle (0.75,0.5);
	\fill[col2,opacity=0.3](0.625,0.5) rectangle (0.75,0.625);
	\fill[col2,opacity=0.3](0.625,0.625) rectangle (0.75,0.75);
	\draw[col4,opacity=0.5](0.6875,0.25)--(0.6875,0.75);
	\foreach \y in {0.3125,0.4375,...,0.6875}{\draw[col2,opacity=0.5](0.625,\y)--(0.75,\y);}
	\draw[col4,opacity=0.5](0.53125,0.4375)--(0.53125,0.5625);
	\foreach \y in {0.46875,0.53125}{\draw[col2,opacity=0.5](0.5,\y)--(0.5625,\y);}
\end{tikzpicture}
}
\caption{}
\label{fig:admissibility4}
\end{subfigure}
\caption{Step-wise illustration of admissible refinement on a mesh corresponding to a THB-spline basis of degree $p=2$. Firstly, (\subref{fig:admissibility1}) illustrates the marking of two elements (in blue and green) of the finest level $\ell$, based on the damage field (in red). Secondly, (\subref{fig:admissibility2}) shows the support extensions of the marked elements as shaded areas. Thirdly, (c) (\cref{fig:admissibility3}) shows the element of level $\ell-2$ intersecting with the support extensions of the marked elements in yellow. Lastly, (d) (\cref{fig:admissibility4}) shows the refined mesh.}
\label{fig:}
\end{figure}
\end{example}

\section{Adaptive Isogeometric Phase-Field Fracture}
\label{sec:adaptive_IGA}
This section presents the novelty of the paper: an adaptive meshing strategy dedicated to phase-field fracture simulations using Truncated Hierarchical B-splines (THB-splines). Following the preliminaries on phase-field fracture models for brittle fracture in \cref{sec:phase-field_fracture} and the preliminaries on THB-splines and admissible refinement in \cref{sec:THB}, this section combines these two topics to present an adaptive meshing strategy for phase-field fracture simulations. Firstly, \cref{subsec:refStrategies} combines the contents of \cref{subsec:THBsplines,subsec:admissible} and {\NEW presents a marking and refinement strategy dedicated to phase-field fracture simulations, based on the observations of fracture suddenness, irreversibility and cross-talk. Secondly, \cref{subsec:LoadStepping} presents load stepping strategies to handle the suddenness of brittle fracture.}

\subsection{Refinement strategies for phase-field fracture}
\label{subsec:refStrategies}
In order to derive adaptive algorithms for phase-field brittle fracture, we first outline two observations about brittle fracture:
\begin{enumerate}
\item \textbf{Brittle fracture is sudden}, when compared to ductile fracture. As a consequence, when performing a load- or displacement-incremented simulation, it is expected that fracture propagation occurs in a small number of load steps. As computational costs are expected to be highest in the fracture propagation regime, it is expected that the total computational costs of the simulation are coming from the load or displacement steps corresponding to fracture propagation.
\item \textbf{Brittle fracture is irreversible}. {\NEW This implies that the damage phase only expands}. Consequently, adaptive meshing algorithms for phase-field brittle fracture can be driven by refinement only, omitting coarsening. The advantage of this approach is that solution projections from old to new meshes are exact.
\end{enumerate}
Besides the aforementioned observations regarding brittle fracture impacting the phase-field fracture simulation, the phenomenon of \textit{cross-talk} as explained by \cite{lian2025cross} is relevant for immersed or trimmed geometries, as well as for phase-field fracture simulations since the damage phase implies zero stiffness locally. In brief, cross-talk appears when the support of a basis function overlaps two sides of a disconnected part of the domain, e.g. due to trimming or damage. In this case, the overlapping basis function contributes to the stiffness matrix, hence add stiffness between the disconnected parts of the domain. For more information, we refer to \cite{lian2025cross} and to \cref{ex:cross-talk} later in this section.

In the following sub-sections, meshing and solution stepping approaches for phase-field brittle fracture are presented, based on the fracture suddenness, irreversibility and cross-talk. Firstly, \cref{subsec:ElementMarking} provides an adaptive meshing algorithm incorporating damage irreversibility and cross-talk. Secondly, \cref{subsec:LoadStepping} provides load stepping algorithms to handle the suddenness of brittle fracture.

\subsection{Phase-field-based element marking}
\label{subsec:ElementMarking}
As highlight in previous sections, the aim of this paper is to provide an adaptive meshing strategy for phase-field brittle fracture simulations using THB-splines. By this means, the aim is to have a sufficiently fine mesh (typically with mesh size $h_{\text{min}}=\ell_0/2$ or $h_{\text{min}}=\ell_0/4$) in the damaged region in order to accurately resolve the damage profile and the material degradation, while keeping the mesh coarse in the undamaged region to reduce computational costs. Due to damage irreversibility, the adaptive meshing strategy only considers element refinement and no coarsening. For the marking criterion, the damage field is sampled in $3^d$ uniformly distributed points in the element, and an element is marked for refinement to the finest level with size $h_{\text{min}}$ if the damage in any of the sampled points exceeds a given threshold value {\NEW $\SCALAR{d}_{\text{min}}$}. {\NEW Instead of taking uniformly distributed points, one could also re-use computed values of the phase-field on the quadrature points which are used for assembly of the elasticity problem. In this paper, however, the values are not cached, thus re-computed, avoiding memory overhead. As will be clear from \cref{fig:implicitexplicit_cpu}, the computational overhead related to the sampling of the phase-field is negligible.}\\

A \emph{field-based} marking strategy as outlined above is motivated by the fact that the marking strategy to avoid cross-talk as outlined below will be based on the damage field as well. Note, however, that this strategy does not guarantee the reduction of the numerical approximation error, since a metric of this error is not incorporated in the marking criterion.\\ 

\begin{example}[Cross-talk]
\label{ex:cross-talk}

Consider a spline basis $\mathcal{B}$ with knot vector $\Xi=\qty{0,0,0,\frac{1}{5},\frac{2}{5},\frac{3}{5},\frac{4}{5},1,1,1}$ and a THB-spline basis $\mathcal{T}$ constructed from $\Xi$ with $\Omega^1=\qty[\frac{1}{5},\frac{4}{5}]$, see the top left and top right plots in \cref{fig:cross-talk}, respectively. Using this basis, the brittle phase-field fracture model (second-order AT2) from \cref{sec:phase-field_fracture} is solved, subject to boundary conditions $u(0) = -1$ and $u(1) = 1$ and with a damage field defined by
\begin{equation}
d(\xi) =
\begin{dcases}
1 & \xi\in\qty[\frac{2}{5},\frac{3}{5}],\\
0 & \text{elsewhere}.
\end{dcases}
\end{equation}
In the bottom row of \cref{fig:cross-talk}, the solution $u(\xi)$ is plotted for the B-spline (bottom left) and THB-spline basis (bottom right). As can be seen from these results, the solution obtained by the B-spline basis is not decoupled, whereas the one obtained by the THB-spline basis is. This is explained by the fact that the middle basis function of the B-spline basis is supported on both sides of the interval $\qty[\frac{2}{5},\frac{3}{5}]$. On the contrary, the THB-spline basis has no functions that are active on both sides of the interval $\qty[\frac{2}{5},\frac{3}{5}]$, resulting in a fully decoupled solution field $u(\xi)$.
\end{example}

\begin{figure}
	\centering
	\begin{tikzpicture}
		\begin{groupplot}
			[
			height=0.125\textheight,
			width=0.5\linewidth,
			legend pos = outer north east,
			group style={
				group name=my plots,
				group size=2 by 2,
				xlabels at=edge bottom,
				x descriptions at=edge bottom,
				y descriptions at=edge left,
				vertical sep=5pt,
				horizontal sep=5pt},
			xlabel={$\xi$},
			grid style={solid, gray!40},
			]
			\nextgroupplot[no markers,
							ylabel={$\beta_i$},
							ytick = {0,1},
							xtick={0.0,0.2,0.4,0.6,0.8,1.0},
							]

			\begin{pgfonlayer}{background}
				\fill[black!05] (axis cs: 0.0,0.0) rectangle (axis cs: 0.4,1.0);
				\fill[col1!50] (axis cs: 0.4,0.0) rectangle (axis cs: 0.6,1.0);
				\fill[black!05] (axis cs: 0.6,0.0) rectangle (axis cs: 1.0,1.0);
			\end{pgfonlayer}

			\foreach \k in {1,...,7}
			{
				\pgfmathtruncatemacro{\y}{\k};
				\pgfmathtruncatemacro{\x}{0};
				\addplot+[black,thin,solid] table[col sep = comma, header=false, y index = {\y}, x index = {\x}] {Data/CrossTalk/TB/basis.csv};
			}

			\nextgroupplot[no markers,
							ylabel={$\beta_i$},
							ytick = {0,1},
							xtick={0.0,0.2,0.3,0.4,0.5,0.6,0.7,0.8,1.0},
							]

			\begin{pgfonlayer}{background}
				\fill[black!05] (axis cs: 0.0,0.0) rectangle (axis cs: 0.2,1.0);
				\fill[col2!50] (axis cs: 0.2,0.0) rectangle (axis cs: 0.4,1.0);
				\fill[col1!50] (axis cs: 0.4,0.0) rectangle (axis cs: 0.6,1.0);
				\fill[col2!50] (axis cs: 0.6,0.0) rectangle (axis cs: 0.8,1.0);
				\fill[black!05] (axis cs: 0.8,0.0) rectangle (axis cs: 1.0,1.0);
			\end{pgfonlayer}

			\foreach \k in {1,...,10}
			{
				\pgfmathtruncatemacro{\y}{\k};
				\pgfmathtruncatemacro{\x}{0};
				\addplot+[black,thin,solid] table[col sep = comma, header=false, y index = {\y}, x index = {\x}] {Data/CrossTalk/THB/basis.csv};
			}

			\nextgroupplot[no markers,
							ylabel={$u$},
							ytick = {-1,0,1},
							xtick={0.0,0.2,0.4,0.6,0.8,1.0},
							xticklabels={0,$\frac{1}{5}$,$\frac{2}{5}$,$\frac{3}{5}$,$\frac{4}{5}$,1},]

			\begin{pgfonlayer}{background}
				\fill[black!05] (axis cs: 0.0,-1.0) rectangle (axis cs: 0.4,1.0);
				\fill[col1!50] (axis cs: 0.4,-1.0) rectangle (axis cs: 0.6,1.0);
				\fill[black!05] (axis cs: 0.6,-1.0) rectangle (axis cs: 1.0,1.0);
			\end{pgfonlayer}

			\addplot+[black,thin,solid] table[col sep = comma, header=false, y index = {1}, x index = {0}] {Data/CrossTalk/TB/solution.csv};

			\nextgroupplot[no markers,
							ylabel={$u$},
							ytick = {-1,0,1},
							xtick={0.0,0.2,0.3,0.4,0.5,0.6,0.7,0.8,1.0},
							xticklabels={0,$\frac{1}{5}$,$\frac{3}{10}$,$\frac{2}{5}$,$\frac{1}{2}$,$\frac{3}{5}$,$\frac{7}{10}$,$\frac{4}{5}$,1},]

			\begin{pgfonlayer}{background}
				\fill[black!05] (axis cs: 0.0,-1.0) rectangle (axis cs: 0.2,1.0);
				\fill[col2!50] (axis cs: 0.2,-1.0) rectangle (axis cs: 0.4,1.0);
				\fill[col1!50] (axis cs: 0.4,-1.0) rectangle (axis cs: 0.6,1.0);
				\fill[col2!50] (axis cs: 0.6,-1.0) rectangle (axis cs: 0.8,1.0);
				\fill[black!05] (axis cs: 0.8,-1.0) rectangle (axis cs: 1.0,1.0);
			\end{pgfonlayer}

			\addplot+[black,thin,solid] table[col sep = comma, header=false, y index = {1}, x index = {0}] {Data/CrossTalk/THB/solution.csv};
		\end{groupplot}
	\end{tikzpicture}
\caption{The effect of cross-talk and the remedy of local refinement on a 1D bar. The top row represents a 1D B-spline basis (top left) and a 1D THB-spline basis (top right) locally refined in the interval $[\frac{1}{5},\frac{4}{5}]$ (blue) such that the middle function of the B-spline basis is eliminated. The bottom row represents the solution of a linear elasticity problem with a degraded material according to \cref{eq:elastic_strain_en,eq:free_en_split} using the red area between $[\frac{2}{5},\frac{3}{5}]$ as damage field, and end displacements of $-1$ and $1$ at $\xi=0$ and $\xi=1$, respectively. On the B-spline basis (bottom left), a gradual distribution of the displacement over the length coordinate $\xi$ is observed, whereas the result obtained on the THB-spline basis shows a discontinuity between the left ($[0,\frac{2}{5}]$) and right ($[\frac{3}{5},1]$) part of the damaged region.}
\label{fig:cross-talk}
\end{figure}

To guarantee the elimination of cross-talk, we note that it appears when basis functions overlap both sides of a disconnected region, e.g., due to trimming or damage, as briefly illustrated in \cref{ex:cross-talk}. Consequently, the marking strategy should ensure that basis functions of coarse levels are eliminated over the damaged region. Considering the construction of THB-splines as presented in \cref{subsec:THBsplines} a function of a coarser level is eliminated if all elements in its support are refined to a finer level. In addition, on an admissible mesh (see \cref{subsec:admissible}) of class $m=2$, the support of a basis function of level $\ell$ only contains elements of level $\ell$ and $\ell+1$. Combining these two observations, it can be concluded that cross-talk can naturally be eliminated by refining the supports of the functions of level $\ell_{\text{max{}}}-1$ admissibly to the finest level $\ell_{\text{max}}$. Looking at \cref{ex:cross-talk}, it suffices to refine a region covering the support of the basis functions active on both sides of the damaged region. This region has a minimum width of $p+1$ elements of level $\ell_{\text{max}}-1$, where $p$ is the polynomial degree of the basis. It should be noted that the choice of the mesh size of the finest level, in this paper chosen to be $h=\frac{\ell_0}{2}$ or $h=\frac{\ell_0}{4}$, according to \cite{greco2024higher, greco_at1, GERASIMOV2019990} guarantees that there is a sufficient number of basis functions in the damaged region provided basis functions with sufficiently low degrees. In case of this paper, we restrict ourselves to quadratic basis functions, since ealier works \cite{greco2024higher,greco_at1} have shown accurate results for these mesh sizes for quadratic basis functions.\\

In practice, finding the support of all basis functions of level $\ell_{\text{max}}-1$ overlapping the damaged region requires an algorithm checking the connectivity of the elements in the support of each basis function of level $\ell_{\text{max}}-1$. As this can be computational challenging, we propose to mark the support extensions of the elements of level $\ell_{\text{max}}-1$ to eliminate cross-talk. This approach re-uses the routine from admissible meshing to find the support extension of an element and is slightly conservative in nature: it marks maximum $2p$ elements around the damaged region. An example of our approach is provided in \cref{ex:cross-talk-elimination} and its algorithm is presented in \cref{alg:refine_elements}.\\

\begin{example}[Cross-talk elimination by local refinement]
\label{ex:cross-talk-elimination}
In this example we consider a domain composed of $8\times 8$ tensor B-spline elements of degree 2. A fictitious damage field is prescribed as shown in red in \cref{fig:adaptivity_cross-talk_a}. We aim to refine the original mesh locally up to $\ell_{\text{max}}=2$ in the damaged region while eliminating cross-talk. The refinement procedure is outlined as follows. (\cref{fig:adaptivity_cross-talk_a}) shows marked elements of level $\ell=0$ in gray, corresponding to the elements overlapping the damaged region. (\cref{fig:adaptivity_cross-talk_b}) shows the mesh after refinement, with marked elements of level $\ell=1$ in gray, corresponding to the elements overlapping the damaged region as well as additional elements marked due to admissibility. (\cref{fig:adaptivity_cross-talk_c}) shows the mesh after refinement of the marked elements in (\cref{fig:adaptivity_cross-talk_b}), with additional elements of level $\ell=1$ from the support extension of the marked elements of level $\ell=1$ in (\cref{fig:adaptivity_cross-talk_b}). Finally, (\cref{fig:adaptivity_cross-talk_d}) shows the final mesh after refinement, where all basis functions of level $\ell=1$ are eliminated over the damaged region by refining the support extensions of the marked elements of level $\ell=1$ to level $\ell=2$.
\end{example}

\begin{figure}
\centering
\begin{subfigure}[t]{0.225\linewidth}
\centering
\resizebox{\linewidth}{!}
{
\begin{tikzpicture}[scale=4]

	\draw[color=red1!10,line width=20](0,0.5)--(0.25,0.5);
	\draw[color=red1!10,line width=20,line cap=round](0.25,0.5)--(0.41,0.5);
	\draw[color=red1!30,line width=15](0,0.5)--(0.25,0.5);
	\draw[color=red1!30,line width=15,line cap=round](0.25,0.5)--(0.41,0.5);
	\draw[color=red1!50,line width=10](0,0.5)--(0.25,0.5);
	\draw[color=red1!50,line width=10,line cap=round](0.25,0.5)--(0.41,0.5);
	\draw[color=red1!70,line width=05](0,0.5)--(0.25,0.5);
	\draw[color=red1!70,line width=05,line cap=round](0.25,0.5)--(0.41,0.5);

	\foreach \x in {0,0.125,...,1.0}{\draw[opacity=0.5](\x,0)--(\x,1);}
	\foreach \y in {0,0.125,...,1.0}{\draw[opacity=0.5](0,\y)--(1,\y);}

	\fill[opacity=0.25] (0,0.375) rectangle (0.5,0.625);

\end{tikzpicture}
}
\caption{}
\label{fig:adaptivity_cross-talk_a}
\end{subfigure}
\begin{subfigure}[t]{0.225\linewidth}
\centering
\resizebox{\linewidth}{!}
{
\begin{tikzpicture}[scale=4]

	\draw[color=red1!10,line width=20](0,0.5)--(0.25,0.5);
	\draw[color=red1!10,line width=20,line cap=round](0.25,0.5)--(0.41,0.5);
	\draw[color=red1!30,line width=15](0,0.5)--(0.25,0.5);
	\draw[color=red1!30,line width=15,line cap=round](0.25,0.5)--(0.41,0.5);
	\draw[color=red1!50,line width=10](0,0.5)--(0.25,0.5);
	\draw[color=red1!50,line width=10,line cap=round](0.25,0.5)--(0.41,0.5);
	\draw[color=red1!70,line width=05](0,0.5)--(0.25,0.5);
	\draw[color=red1!70,line width=05,line cap=round](0.25,0.5)--(0.41,0.5);

	\foreach \x in {0,0.125,...,1.0}{\draw[opacity=0.5](\x,0)--(\x,1);}
	\foreach \y in {0,0.125,...,1.0}{\draw[opacity=0.5](0,\y)--(1,\y);}
	\foreach \x in {0.0625,0.1875,...,0.5}{\draw[opacity=0.5](\x,0.375)--(\x,0.625);}
	\foreach \y in {0.4375,0.5,...,0.5625}{\draw[opacity=0.5](0,\y)--(0.5,\y);}

	\fill[opacity=0.25] (0,0.4375) rectangle (0.4375,0.5625);

	\fill[opacity=0.125] (0,0.25) rectangle (0.5,0.375);
	\fill[opacity=0.125] (0,0.625) rectangle (0.5,0.75);
	\fill[opacity=0.125] (0.5,0.25) rectangle (0.625,0.75);
\end{tikzpicture}
}
\caption{}
\label{fig:adaptivity_cross-talk_b}
\end{subfigure}
\begin{subfigure}[t]{0.225\linewidth}
\centering
\resizebox{\linewidth}{!}
{
\begin{tikzpicture}[scale=4]

	\draw[color=red1!10,line width=20](0,0.5)--(0.25,0.5);
	\draw[color=red1!10,line width=20,line cap=round](0.25,0.5)--(0.41,0.5);
	\draw[color=red1!30,line width=15](0,0.5)--(0.25,0.5);
	\draw[color=red1!30,line width=15,line cap=round](0.25,0.5)--(0.41,0.5);
	\draw[color=red1!50,line width=10](0,0.5)--(0.25,0.5);
	\draw[color=red1!50,line width=10,line cap=round](0.25,0.5)--(0.41,0.5);
	\draw[color=red1!70,line width=05](0,0.5)--(0.25,0.5);
	\draw[color=red1!70,line width=05,line cap=round](0.25,0.5)--(0.41,0.5);

	\foreach \x in {0,0.125,...,1.0}{\draw[opacity=0.5](\x,0)--(\x,1);}
	\foreach \y in {0,0.125,...,1.0}{\draw[opacity=0.5](0,\y)--(1,\y);}
	\foreach \x in {0.0625,0.1875,...,0.625}{\draw[opacity=0.5](\x,0.25)--(\x,0.75);}
	\foreach \y in {0.3125,0.4375,...,0.6875}{\draw[opacity=0.5](0,\y)--(0.625,\y);}
	\foreach \x in {0.03125,0.09375,...,0.4375}{\draw[opacity=0.5](\x,0.4375)--(\x,0.5625);}
	\foreach \y in {0.46875,0.53125}{\draw[opacity=0.5](0,\y)--(0.4375,\y);}

	\fill[opacity=0.375] (0,0.375) rectangle (0.4375,0.4375);
	\fill[opacity=0.375] (0,0.5625) rectangle (0.4375,0.625);
	\fill[opacity=0.375] (0.4375,0.375) rectangle (0.5,0.625);

\end{tikzpicture}
}
\caption{}
\label{fig:adaptivity_cross-talk_c}
\end{subfigure}
\begin{subfigure}[t]{0.225\linewidth}
\centering
\resizebox{\linewidth}{!}
{
\begin{tikzpicture}[scale=4]

	\draw[color=red1!10,line width=20](0,0.5)--(0.25,0.5);
	\draw[color=red1!10,line width=20,line cap=round](0.25,0.5)--(0.41,0.5);
	\draw[color=red1!30,line width=15](0,0.5)--(0.25,0.5);
	\draw[color=red1!30,line width=15,line cap=round](0.25,0.5)--(0.41,0.5);
	\draw[color=red1!50,line width=10](0,0.5)--(0.25,0.5);
	\draw[color=red1!50,line width=10,line cap=round](0.25,0.5)--(0.41,0.5);
	\draw[color=red1!70,line width=05](0,0.5)--(0.25,0.5);
	\draw[color=red1!70,line width=05,line cap=round](0.25,0.5)--(0.41,0.5);

	\foreach \x in {0,0.125,...,1.0}{\draw[opacity=0.5](\x,0)--(\x,1);}
	\foreach \y in {0,0.125,...,1.0}{\draw[opacity=0.5](0,\y)--(1,\y);}
	\foreach \x in {0.0625,0.1875,...,0.625}{\draw[opacity=0.5](\x,0.25)--(\x,0.75);}
	\foreach \y in {0.3125,0.4375,...,0.6875}{\draw[opacity=0.5](0,\y)--(0.625,\y);}	
	\foreach \x in {0.03125,0.09375,...,0.46875}{\draw[opacity=0.5](\x,0.375)--(\x,0.625);}
	\foreach \y in {0.40625,0.46875,...,0.59375}{\draw[opacity=0.5](0,\y)--(0.5,\y);}
\end{tikzpicture}
}
\caption{}
\label{fig:adaptivity_cross-talk_d}
\end{subfigure}
\caption{Illustration of cross-talk elimination through refinement of support extensions, provided a damage field (in red). The shaded cells represent elements marked for refinement. Firstly, (\subref{fig:adaptivity_cross-talk_a}) shows the damage on the coarsest (initial) mesh, where every element with sufficiently high damage value is marked for refinement. Secondly, (\subref{fig:adaptivity_cross-talk_b}) shows the mesh after refinement to level 1, where elements with sufficiently high damage value are again marked for refinement. The light shaded region of elements of level $\ell=0$ represents element marked by the admissibility algorithm. Thirdly, (\subref{fig:adaptivity_cross-talk_c}) shows the mesh after refinement to the finest level $\ell=2$, where now elements corresponding to the support extension of the previously marked elements are marked for refinement (in this case there is no need for additional refinement through admissibility). Finally, (\subref{fig:adaptivity_cross-talk_d}) shows the final mesh after refinement to the finest level $\ell=2$, where all basis functions of level $\ell=1$ have been eliminated over the damaged region, hence cross-talk is eliminated.}
\label{fig:adaptivity_cross-talk}
\end{figure}

\subsection{Load-stepping with adaptive meshing}
\label{subsec:LoadStepping}
Besides element marking strategies incorporating damage irreversibility and cross-talk elimination, another important aspect of adaptive meshing for phase-field brittle fracture simulations is the suddenness of fracture. When performing brittle fracture simulations using load- or displacement-stepping schemes, propagation of the damage field progresses slowly over a large number of steps, until a critical load or displacement has been reached after which the damage field evolves significantly. In case of tensor-product meshes which are fine over the whole domain, this suddenness does not pose any particular challenge, as the mesh provides sufficient resolution everywhere. However, in case of adaptive meshes, the suddenness of brittle fracture can cause the damage field to propagate significantly in a single load step, causing the damage field to be poorly resolved on coarse parts of the domain. To alleviate this issue, different strategies can be adopted when combining load stepping and adaptive meshing:
\begin{itemize}
	\item \textbf{Explicit mesh adaptivity}: After solution step $k$, the mesh will be refined and the solutions $\VEC{u}^k$ and $d^k$ will be projected onto the new mesh. Afterwards, the solution step $k+1$ is solved.
	\item \textbf{Implicit mesh adaptivity}: After solution step $k$, the mesh will be refined and the solutions $\VEC{u}^k$ and $d^k$ will be projected onto the new mesh. Afterwards, solution step $k$ is repeated with the new solutions and the process is repeated until new elements are refined in the mesh.
	\item \textbf{Quasi-implicit mesh adaptivity}: As in implicit mesh adaptivity, solution steps are repeated after mesh refinement. However, repetition is only done if the change in the mesh is `significant' compared to the previous refinement iteration. If not, the algorithm proceeds to the next solution step after projection of the solutions onto the mesh obtained in the latest refinement iteration.
\end{itemize}
An algorithm for adaptive load stepping used in this paper is presented in \cref{alg:adaptive_load_step} in \ref{app:algorithms}. In this algorithm, quasi-implicit mesh adaptivity is driven by the ratio in the number of elements before and after refinement. If this ratio is below a given tolerance, the change in the mesh is considered insignificant and the algorithm proceeds to the next load step. Otherwise, the current load step is repeated on the new mesh. Alternatively to the ratio in the number of elements, other metrics can be used to determine whether the change in the mesh is significant, e.g., the total area of the new elements added to the mesh.

\begin{example}[Adaptive phase-field refinement]
\label{ex:mesh_adaptivity}
In this example, the schematic propagation of the damage field and the computational mesh from \cref{fig:mesh_adaptivity} are studied. The panels in \cref{fig:mesh_adaptivity} are illustrative, hence not a result from computations. Initially, consider the damage field in load step $k-2$ (\cref{fig:propagation_step_km2}), with a two-level mesh refined around the non-zero part of the damage field. When arriving in load step $k-1$ (\cref{fig:propagation_step_km1}), the damage has propagated towards the right of the domain by a little. In this case, the two extra elements are not considered a significant change in the mesh, hence the mesh adaptivity is explicit. Now, let us assume that significant propagation of the damage field happens in load step $k$. Then, the damage field computed on the mesh of load-step $k-1$ results in the damage field depicted in \cref{fig:propagation_step_k_it_i}. Consequently, a relatively large number of elements is marked for refinement, leading to a refinement iteration for load step $k$ on the new mesh. As a result, the approximation of the damage field is improved, as can be seen in \cref{fig:propagation_step_k_it_j}, yielding no extra elements to be refined.
\end{example}

\begin{figure}
\centering
\begin{subfigure}[t]{0.225\linewidth}
\centering
\resizebox{\linewidth}{!}
{
\begin{tikzpicture}[scale=4]
	\draw[color=red1!10,line width=20](0,0.5)--(0.25,0.5);
	\draw[color=red1!10,line width=20,line cap=round](0.25,0.5)--(0.41,0.5);
	\draw[color=red1!30,line width=15](0,0.5)--(0.25,0.5);
	\draw[color=red1!30,line width=15,line cap=round](0.25,0.5)--(0.41,0.5);
	\draw[color=red1!50,line width=10](0,0.5)--(0.25,0.5);
	\draw[color=red1!50,line width=10,line cap=round](0.25,0.5)--(0.41,0.5);
	\draw[color=red1!70,line width=05](0,0.5)--(0.25,0.5);
	\draw[color=red1!70,line width=05,line cap=round](0.25,0.5)--(0.41,0.5);

	\foreach \x in {0,0.125,...,1.0}{\draw[opacity=0.5](\x,0)--(\x,1);}
	\foreach \y in {0,0.125,...,1.0}{\draw[opacity=0.5](0,\y)--(1,\y);}
	\foreach \x in {0.0625,0.1875,...,0.5}{\draw[opacity=0.5](\x,0.375)--(\x,0.625);}
	\foreach \y in {0.4375,0.5625}{\draw[opacity=0.5](0,\y)--(0.5,\y);}
\end{tikzpicture}
}
\caption{Load step $k-2$}
\label{fig:propagation_step_km2}
\end{subfigure}
\hfill
\begin{subfigure}[t]{0.225\linewidth}
\centering
\resizebox{\linewidth}{!}
{
\begin{tikzpicture}[scale=4]
	\draw[color=red1!10,line width=20](0,0.5)--(0.25,0.5);
	\draw[color=red1!10,line width=20,line cap=round](0.25,0.5)--(0.5,0.5);
	\draw[color=red1!30,line width=15](0,0.5)--(0.25,0.5);
	\draw[color=red1!30,line width=15,line cap=round](0.25,0.5)--(0.5,0.5);
	\draw[color=red1!50,line width=10](0,0.5)--(0.25,0.5);
	\draw[color=red1!50,line width=10,line cap=round](0.25,0.5)--(0.5,0.5);
	\draw[color=red1!70,line width=05](0,0.5)--(0.25,0.5);
	\draw[color=red1!70,line width=05,line cap=round](0.25,0.5)--(0.5,0.5);

	\foreach \x in {0,0.125,...,1.0}{\draw[opacity=0.5](\x,0)--(\x,1);}
	\foreach \y in {0,0.125,...,1.0}{\draw[opacity=0.5](0,\y)--(1,\y);}
	\foreach \x in {0.0625,0.1875,...,0.5}{\draw[opacity=0.5](\x,0.375)--(\x,0.625);}
	\foreach \y in {0.4375,0.5625}{\draw[opacity=0.5](0,\y)--(0.5,\y);}
	\fill[black,opacity=0.2] (0.5,0.375) rectangle (0.625,0.625);
\end{tikzpicture}
}
\caption{Load step $k-1$}
\label{fig:propagation_step_km1}
\end{subfigure}
\hfill
\begin{subfigure}[t]{0.225\linewidth}
\centering
\resizebox{\linewidth}{!}
{
\begin{tikzpicture}[scale=4]
	\draw[color=red1!10,line width=20](0,0.5)--(0.5,0.5);
	\draw[color=red1!10,line width=20,line cap=round,line join=round](0.5,0.5)--(0.625,0.375)--(0.75,0.3);
	\draw[color=red1!30,line width=15](0,0.5)--(0.5,0.5);
	\draw[color=red1!30,line width=15,line cap=round,line join=round](0.5,0.5)--(0.625,0.375)--(0.75,0.3);
	\draw[color=red1!50,line width=10](0,0.5)--(0.5,0.5);
	\draw[color=red1!50,line width=10,line cap=round,line join=round](0.5,0.5)--(0.625,0.375)--(0.75,0.3);
	\draw[color=red1!70,line width=05](0,0.5)--(0.5,0.5);
	\draw[color=red1!70,line width=05,line cap=round,line join=round](0.5,0.5)--(0.625,0.375)--(0.75,0.3);

	\foreach \x in {0,0.125,...,1.0}{\draw[opacity=0.5](\x,0)--(\x,1);}
	\foreach \y in {0,0.125,...,1.0}{\draw[opacity=0.5](0,\y)--(1,\y);}
	\foreach \x in {0.0625,0.1875,...,0.5625}{\draw[opacity=0.5](\x,0.375)--(\x,0.625);}
	\foreach \y in {0.4375,0.5625}{\draw[opacity=0.5](0,\y)--(0.625,\y);}

	\fill[black,opacity=0.2] (0.500,0.250) rectangle (0.625,0.375);
	\fill[black,opacity=0.2] (0.625,0.125) rectangle (0.750,0.500);
	\fill[black,opacity=0.2] (0.750,0.125) rectangle (0.875,0.375);
\end{tikzpicture}
}
\caption{Load step $k$, iteration $i$.}
\label{fig:propagation_step_k_it_i}
\end{subfigure}
\hfill
\begin{subfigure}[t]{0.225\linewidth}
\centering
\resizebox{\linewidth}{!}
{
\begin{tikzpicture}[scale=4]
	\draw[color=red1!10,line width=20](0,0.5)--(0.5,0.5);
	\draw[color=red1!10,line width=20,line cap=round,line join=round](0.5,0.5)--(0.75,0.25);
	\draw[color=red1!30,line width=15](0,0.5)--(0.5,0.5);
	\draw[color=red1!30,line width=15,line cap=round,line join=round](0.5,0.5)--(0.75,0.25);
	\draw[color=red1!50,line width=10](0,0.5)--(0.5,0.5);
	\draw[color=red1!50,line width=10,line cap=round,line join=round](0.5,0.5)--(0.75,0.25);
	\draw[color=red1!70,line width=05](0,0.5)--(0.5,0.5);
	\draw[color=red1!70,line width=05,line cap=round,line join=round](0.5,0.5)--(0.75,0.25);

	\foreach \x in {0,0.125,...,1.0}{\draw[opacity=0.5](\x,0)--(\x,1);}
	\foreach \y in {0,0.125,...,1.0}{\draw[opacity=0.5](0,\y)--(1,\y);}
	\foreach \x in {0.0625,0.1875,...,0.5625}{\draw[opacity=0.5](\x,0.375)--(\x,0.625);}
	\foreach \y in {0.4375,0.5625}{\draw[opacity=0.5](0,\y)--(0.625,\y);}

	\draw[opacity=0.5](0.5625,0.25)--(0.5625,0.375);
	\draw[opacity=0.5](0.6825,0.125)--(0.6825,0.5);
	\draw[opacity=0.5](0.8125,0.125)--(0.8125,0.375);

	\draw[opacity=0.5](0.625,0.1875)--(0.875,0.1875);
	\draw[opacity=0.5](0.5,0.3125)--(0.875,0.3125);
	\draw[opacity=0.5](0.5,0.4375)--(0.75,0.4375);
\end{tikzpicture}
}
\caption{Load step $k$, iteration $j>i$}
\label{fig:propagation_step_k_it_j}
\end{subfigure}
\caption{Schematic representation of sudden phase-field propagation combined with mesh adaptivity. From left to right: (\subref{fig:propagation_step_km2}) phase-field at step $k-2$, (\subref{fig:propagation_step_km1}) phase-field at step $k-1$; (\subref{fig:propagation_step_km1}) slightly propagated phase-field with elements marked for refinement (gray); (\subref{fig:propagation_step_k_it_i}) significantly propagated phase-field, partially distorted due to representation on coarse elements with elements marked for refinement in gray; (\subref{fig:propagation_step_k_it_j}) phase-field represented on a locally refined mesh according to the phase-field in iteration $i<j$ as presented in \subref{fig:propagation_step_k_it_i}. }
\label{fig:mesh_adaptivity}
\end{figure}

\subsection{Phase-field initialization on THB-meshes}
\label{subsec:PF_initialization_THB}
As discussed in \cref{subsec:PF_initialization}, the initial phase-field can be constructed using the IPF method introduced by \cite{greco2024higher}. {\NEW When using THB-splines as a basis for phase-field fracture, the mesh must be initialized so that the IPF yields a phase-field identical to one obtained on a uniform mesh with size $h_{\text{min}}$ -- the finest mesh size of the THB basis to be constructed. Since the IPF method performs a local $L_2$-projection onto all basis functions with support in a $\beta$-neighborhood around the crack, the initial mesh must contain all finest-level basis functions within this neighborhood.} For THB-splines, this is achieved by applying the same method used for cross-talk elimination (see \cref{subsec:ElementMarking}), where the $\beta$-neighborhood is used to mark elements for refinement.


\section{Benchmark examples}
\label{sec:benchmarks}
In this section, the proposed adaptive refinement framework for brittle phase-field fracture simulations is benchmarked using two basic examples: the Single Edge Notched (SEN) tensile and shear tests from \cite{miehe2010IJNME, GERASIMOV2019990, greco_at1}. The goal of the benchmarks is to assess the efficiency of the proposed framework in terms of computational costs, considering different phase-field formulations and meshing strategies. Therefore, \cref{subsec:benchmarks_definitions} introduces the benchmark problems, \cref{subsec:modelcomparison} assesses the effect of different phase-field formulations on the model performance and finally \cref{subsec:adaptivity,subsec:CPU} elaborate on the efficiency gains followed from different meshing settings. Unless stated otherwise, all simulations are performed with tolerances $\mathtt{TOL}_{\text{Pic},\VEC{u}}=10^{-5}$, $\mathtt{TOL}_{\text{PSOR},\Delta\SCALAR{d}}=10^{-9}$, $\mathtt{TOL}_{\text{stag}}=10^{-5}$, $\mathtt{TOL}_{ref}=0.005$. {\NEW Furthermore, we consider results obtained on tensor-product B-spline bases as the ground truth in this paper when assessing the results from our adaptive meshing framework using THB-splines.}.\\

The framework proposed in this paper is implemented in the Geometry + Simulation modules \cite{mantzaflarisGeometrySimulationModules2025,Juttler2014}. Assembly of the linear systems is performed element-by-element with $\qty(p+1)^d$ quadrature points per element, using shared-memory parallelization via OpenMP on 10 threads of an Intel${}^\text{\textregistered}$ Xeon${}^\text{\textregistered}$ Silver 4316 CPU.  Unless stated otherwise, the Pardiso \cite{schenk2004solving} solver from Intel's Math Kernel Library (MKL) is used to solve the linear systems related to the elasticity problem, while the (inherently serial) PSOR method \cite{MARENGO2021114137} is used to solve the phase-field equations, as stated in \cref{subsubsec:variations}.\\

It should be noted that the computational costs in this section are reported in terms of CPU walltime and in terms of degrees of freedom. Although the latter is more common, it implies that a lower number of degrees of freedom leads to lower computational costs. However, this implicit assumption omits any computational overhead imposed by changing the spline basis from B-splines to THB-splines. On the other hand, CPU walltime heavily depends on the efficiency of the implementation of different spline constructions, and optimization of assemblers and solvers in favor of one spline construction over the other. The aforementioned points motivate that reporting of only degrees of freedom or only CPU time is not fully objective in either way.\\

{\NEW In the supplementary material to this paper, we provide videos of the damage field evolution for the different phase-field models and meshing strategies considered in this paper. Videos 1-4 show the damage evolution on coarse (Video 1 and 3) and coarse fine (Video 2 and 4), respectively, comparing tensor-product and THB with hybrid load stepping for the shear (Video 1-2) and tensile tests (Video 3-4). Videos 5-8 show the damage evolution for the different load stepping schemes (explicit, implicit and hybrid) for the shear (Video 5-6) and tensile tests (Video 7-8) on coarse (Video 5 and 7) and fine (Video 6 and 8) meshes.}\\

\subsection{Problem definitions}
\label{subsec:benchmarks_definitions}
The benchmark tests provided in this section are based on the SEN tensile and shear tests from \cite{miehe2010IJNME}. These tests have been widely studied by \cite{greco2024higher, greco_at1} among others. \Cref{fig:SEN} provides the geometry, boundary conditions and parameters related to the SEN tensile (\cref{fig:SEN_tensile}) and shear (\cref{fig:SEN_shear}) tests. For the analyses, B-spline and THB-spline bases of degree $p=2$ are used, and in both cases the results are studied with respect to a coarse and a fine reference element size. For the B-spline basis, the reference element size is the size everywhere in the domain, whereas for the THB-spline basis it is the size of the finest level used. {\NEW Furthermore, both cases are simulated in a displacement-controlled quasi-static fashion, meaning that dynamic effects are omitted and displacements are applied incrementally.}\\

\begin{figure}[]
\centering
\begin{subfigure}{\linewidth}
\centering
\begin{minipage}{0.35\linewidth}
\centering
\resizebox{\linewidth}{!}
{
\usetikzlibrary{patterns.meta}
\begin{tikzpicture}[scale=3]
\fill[white] (-0.20,-0.20) rectangle (1.20,1.20);
\filldraw[fill=black!50, draw=black,top color=black!20,bottom color=black!50](0,0) rectangle (1,1);
\draw[ultra thick](0,0.5) -- (0.5,0.5);
\fill[pattern={Lines[angle=45,distance=4pt]},pattern color=black!70] (-0.15,-0.1) rectangle (1.15,0);
\draw[ultra thick](0,0) -- (1,0);
\node[above] at (0.5,0.0){$\mathbf{u}=[0,0]$};

\draw[thick](-0.01,1) -- (1.01,1);
\draw[thick](-0.01,0.9) -- (-0.01,1.10);
\draw[thick](-0.075,0.9) -- (-0.075,1.10);
\draw[thick](1.01,0.9) -- (1.01,1.10);
\draw[thick](1.075,0.9) -- (1.075,1.10);
\fill[pattern={Lines[angle=45,distance=4pt]},pattern color=black!70] (-0.075,0.9) rectangle (-0.15,1.10);
\fill[pattern={Lines[angle=45,distance=4pt]},pattern color=black!70] (1.075,0.9) rectangle (1.15,1.10);
\foreach \y in {0.925,0.975,...,1.076}
{
\draw (-0.045,\y) circle [radius=0.025];
\draw (1.045,\y) circle [radius=0.025];
}

\foreach \x in {0.04,0.08,...,1}
{
\draw[-latex] (\x,1.01) -- (\x,1.11);
}
\node[below] at (0.5,0.99){$\mathbf{u}=[0,u_y]$};

\draw[latex-latex](0,-0.1) -- (1,-0.1) node[midway,below]{$L$};
\draw[latex-latex](1.1,0) -- (1.1,1) node[midway,right]{$L$};
\draw[latex-latex](-0.1,0) -- (-0.1,0.5) node[midway,left]{$L/2$};
\draw[latex-latex](0,0.45) -- (0.5,0.45) node[midway,below]{$L/2$};
\end{tikzpicture}
}
\end{minipage}
\hfill
\begin{minipage}{0.6\linewidth}
\footnotesize
\centering
\begin{tabular}{p{0.15\linewidth}p{0.15\linewidth}p{0.1\linewidth}p{0.35\linewidth}}
\toprule
\textbf{Parameter} & \textbf{Value} & \textbf{Unit} & \textbf{Description}\\
\midrule
$L$ & $1.0$ & $\text{mm}$ & Length parameter\\
$u_x$ & $[0,6]\cdot 10^{-3}$ & $\text{mm}$ & Incremental horizontal displacement\\
$E$ & $210$ & $\text{kN}/\text{mm}^2$ & Young's modulus\\
$\nu$ & $0.3$ & $-$ & Poisson's ratio\\
$G_c$ & $2.7 \cdot 10^{-3}$ & $\text{kN}/\text{mm}$ & Toughness\\
$\ell_0$ & $0.015$ & $\text{mm}$ & Internal length\\
\bottomrule
\end{tabular}
\end{minipage}
\caption{SEN tensile test.}
\label{fig:SEN_tensile}
\end{subfigure}

\begin{subfigure}{\linewidth}
\centering
\begin{minipage}{0.35\linewidth}
\centering
\resizebox{\linewidth}{!}
{
\begin{tikzpicture}[scale=3]
\fill[white] (-0.20,-0.20) rectangle (1.20,1.20);
\filldraw[fill=black!50, draw=black,top color=black!20,bottom color=black!50](0,0) rectangle (1,1);
\draw[ultra thick](0,0.5) -- (0.5,0.5);
\fill[pattern={Lines[angle=45,distance=4pt]},pattern color=black!70] (-0.15,-0.1) rectangle (1.15,0);
\draw[ultra thick](0,0) -- (1,0);
\node[above] at (0.5,0.0){$\mathbf{u}=[0,0]$};

\draw[ultra thick](0,1) -- (1,1);
\draw[thick](0,1.06) -- (1,1.06);
\foreach \x in {0.04,0.1,...,0.99}
{
\draw (\x,1.03) circle [radius=0.025];
}

\foreach \x in {0,0.1,0.2,...,0.99}
{
\draw[-latex] (\x,0.95) -- (\x+0.09,0.95);
}
\node[below] at (0.5,0.95){$\mathbf{u}=[u_x,0]$};

\draw[latex-latex](0,-0.1) -- (1,-0.1) node[midway,below]{$L$};
\draw[latex-latex](1.1,0) -- (1.1,1) node[midway,right]{$L$};
\draw[latex-latex](-0.1,0) -- (-0.1,0.5) node[midway,left]{$L/2$};
\draw[latex-latex](0,0.45) -- (0.5,0.45) node[midway,below]{$L/2$};
\end{tikzpicture}
}
\end{minipage}
\hfill
\begin{minipage}{0.6\linewidth}
\centering
\footnotesize
\centering
\begin{tabular}{p{0.15\linewidth}p{0.15\linewidth}p{0.1\linewidth}p{0.35\linewidth}}
\toprule
\textbf{Parameter} & \textbf{Value} & \textbf{Unit} & \textbf{Description}\\
\midrule
$L$ & $1.0$ & $\text{mm}$ & Length parameter\\
$u_x$ & $[6, 12]\cdot 10^{-3}$ & $\text{mm}$ & Incremental horizontal displacement\\
$E$ & $210$ & $\text{kN}/\text{mm}^2$ & Young's modulus\\
$\nu$ & $0.3$ & $-$ & Poisson's ratio\\
$G_c$ & $2.7 \cdot10^{-3}$ & $\text{kN}/\text{mm}$ & Toughness\\
$\ell_0$ & $0.010$ & $\text{mm}$ & Internal length\\
\bottomrule
\end{tabular}
\end{minipage}
\caption{SEN shear test}
\label{fig:SEN_shear}
\end{subfigure}
\caption{Benchmark setup for the SEN (\subref{fig:SEN_shear}) shear and (\subref{fig:SEN_tensile}) tensile tests including geometry, boundary conditions and parameters.}
\label{fig:SEN}
\end{figure}

For the tensile test, the characteristic length $\ell_0=0.015\:[\text{mm}]$ is slightly higher than for the shear test, therefore a larger element size can be used. For the fine mesh, an element size of $h=0.003125\:[\text{mm}]$ is chosen, corresponding with $\frac{\ell_0}{h}=4.8$. For the THB-spline the hierarchy is composed of tensor B-spline levels with element sizes $h=0.05\:[\text{mm}]$ (level 0), $h=0.025\:[\text{mm}]$ (level 1), $h=0.0125\:[\text{mm}]$ (level 2), $h=0.00625\:[\text{mm}]$ (level 3), and $h=0.003125\:[\text{mm}]$ (level 4), the latter corresponding to the level of the B-spline basis. For the coarse mesh, the B-spline basis has mesh size $h=0.00625\:[\text{mm}]$ (hence $\frac{\ell_0}{h}=2.4$), whereas the THB-spline basis has the same levels up to level 3. Using the mesh initialization method discussed in \cref{subsec:PF_initialization_THB}, the fine THB-spline basis has 3316 elements and 3076 degrees of freedom and the coarse THB-spline basis has 1822 elements and 1748 degrees of freedom. The numbers of elements and degrees of freedom of these bases are significantly lower than the number of element in the fine and coarse tensor B-spline bases with, respectively, $320\times320$ elements and 103684 degrees of freedom and $160\times160$ elements and 26244 degrees of freedom for $p=2$.\\

For the shear test, the characteristic length is $\ell_0=0.010\:[\text{mm}]$, hence the mesh sizes will be slightly smaller compared to those obtained for the tensile test. The fine B-spline mesh has an element size of $h=0.0025\:[\text{mm}]$, corresponding to $\frac{\ell_0}{h}=4$. Similar to the tensile case, the THB-spline hierarchy is chosen such that its finest level is equal to the B-spline mesh. Therefore, for the fine THB-spline mesh, the hierarchy is composed of levels with mesh sizes $h=0.04\:[\text{mm}]$ (level 0), $h=0.02\:[\text{mm}]$ (level 1), $h=0.01\:[\text{mm}]$ (level 2), $h=0.005\:[\text{mm}]$ (level 3) and $h=0.0025\:[\text{mm}]$ (level 4). For the coarse B-spline basis, the mesh size is $h=0.005\:[\text{mm}]$, corresponding to $\frac{\ell_0}{h}=2$, and the coars THB-spline hierarchy is composed of levels 0 up to 3. Using the mesh initialization discussed in \cref{subsec:PF_initialization_THB}, the fine THB-spline basis has 4162 elements and 3868 degrees of freedom and the coarse THB-spline basis has 2308 elements and 2220 degrees of freedom. The number of elements and degrees of freedom of these bases are significantly lower than the number of element in the fine and coarse B-spline bases, with, respectively, $400\times400$ elements and 161604 degrees of freedom and $200\times200$ elements and 40804 degrees of freedom for $p=2$.\\

\subsection{Phase-field model comparison}
\label{subsec:modelcomparison}
Before elaborating on the adaptive refinement scheme in \cref{sec:adaptive_IGA}, we first elaborate on the choice of phase-field fracture model and its effect on the propagation of the damage field as well as on the reaction forces and the dissipated energy. {\NEW The purpose of this section primarily is to present the results on the tensor-product B-spline bases for different phase-field formulations.}

\begin{figure}[]
\centering
\begin{subfigure}{\linewidth}
\centering
\begin{tikzpicture}
\def\minIdx{15}
\def\maxIdx{20}
\def\minDispl{0.0042}
\def\stepDispl{0.0003}
\pgfmathsetmacro{\maxDispl}{\minDispl+\stepDispl*(\maxIdx-\minIdx)}
\pgfmathsetmacro{\minCol}{10}
\pgfmathsetmacro{\maxCol}{100}
\begin{groupplot}[
        group style={
            group name=myplot,
            group size=2 by 4,
            horizontal sep=5pt,
            vertical sep=5pt,
            xticklabels at=edge bottom,
            xlabels at=edge bottom,
            yticklabels at=edge left,
            ylabels at=edge left
        },
    width=0.45\linewidth,
    axis equal image,
    xlabel={$x$},
    x unit={mm},
    ylabel={$y$},
    y unit={mm},
    xmin = 0,xmax = 1,
    ymin = 0.4,ymax = 0.6,
    point meta min=\minDispl,
    point meta max=\maxDispl,
    every mark/.append style={solid},
    colormap={mycolormap}{color=(col1!\minCol!col2) color=(col1!\maxCol!col2)},
    cycle list={[samples of colormap={\maxIdx-\minIdx+1} of mycolormap]},
    colorbar style={axis equal image=false,width=0.2cm,x unit={}, y unit={}},
]
\nextgroupplot
\node[anchor=north west,scale=0.8] at (rel axis cs: 0.05,0.95){AT-1, Order 2, $h=\ell_0/2$};
\foreach \step[count=\i] in {\minIdx,...,\maxIdx}
{
    \pgfmathsetmacro{\col}{100-90*\i/(\maxIdx-\minIdx)}
    \pgfmathsetmacro{\disp}{\minDispl+\i*\stepDispl}
    \edef\temp
    {\noexpand\addplot+[only marks,mark size=0.3pt] table[header=true,col sep=comma, x index=1, y index=2,] {Data/TB/tensile_TB_0.50l0_AT-1_Order2_contours/contour_\step.csv};
    }
    \temp
}

\nextgroupplot[colorbar,colorbar style={axis equal image=false,width=0.2cm,x unit={}, y unit={mm},ylabel={Vertical top displacement},height=
        4*\pgfkeysvalueof{/pgfplots/parent axis height}+
        3*\pgfkeysvalueof{/pgfplots/group/vertical sep}}]
\node[anchor=north west,scale=0.8] at (rel axis cs: 0.05,0.95){AT-1, Order 2, $h=\ell_0/4$};
\foreach \step[count=\i] in {\minIdx,...,\maxIdx}
{
    \pgfmathsetmacro{\col}{100-90*\i/(\maxIdx-\minIdx)}
    \pgfmathsetmacro{\disp}{\minDispl+\i*\stepDispl}
    \edef\temp
    {\noexpand\addplot+[only marks,mark size=0.3pt] table[header=true,col sep=comma, x index=1, y index=2,] {Data/TB/tensile_TB_0.25l0_AT-1_Order2_contours/contour_\step.csv};
    }
    \temp
}
\nextgroupplot
\node[anchor=north west,scale=0.8] at (rel axis cs: 0.05,0.95){AT-1, Order 4, $h=\ell_0/2$};
\foreach \step[count=\i] in {\minIdx,...,\maxIdx}
{
    \pgfmathsetmacro{\col}{100-90*\i/(\maxIdx-\minIdx)}
    \pgfmathsetmacro{\disp}{\minDispl+\i*\stepDispl}
    \edef\temp
    {\noexpand\addplot+[only marks,mark size=0.3pt] table[header=true,col sep=comma, x index=1, y index=2,] {Data/TB/tensile_TB_0.50l0_AT-1_Order4_contours/contour_\step.csv};
    }
    \temp
}

\nextgroupplot[]
\node[anchor=north west,scale=0.8] at (rel axis cs: 0.05,0.95){AT-1, Order 4, $h=\ell_0/4$};
\foreach \step[count=\i] in {\minIdx,...,\maxIdx}
{
    \pgfmathsetmacro{\col}{100-90*\i/(\maxIdx-\minIdx)}
    \pgfmathsetmacro{\disp}{\minDispl+\i*\stepDispl}
    \edef\temp
    {\noexpand\addplot+[only marks,mark size=0.3pt] table[header=true,col sep=comma, x index=1, y index=2,] {Data/TB/tensile_TB_0.25l0_AT-1_Order4_contours/contour_\step.csv};
    }
    \temp
}

\nextgroupplot
\node[anchor=north west,scale=0.8] at (rel axis cs: 0.05,0.95){AT-2, Order 2, $h=\ell_0/2$};
\foreach \step[count=\i] in {\minIdx,...,\maxIdx}
{
    \pgfmathsetmacro{\col}{100-90*\i/(\maxIdx-\minIdx)}
    \pgfmathsetmacro{\disp}{\minDispl+\i*\stepDispl}
    \edef\temp
    {\noexpand\addplot+[only marks,mark size=0.3pt] table[header=true,col sep=comma, x index=1, y index=2,] {Data/TB/tensile_TB_0.50l0_AT-2_Order2_contours/contour_\step.csv};
    }
    \temp
}

\nextgroupplot[]
\node[anchor=north west,scale=0.8] at (rel axis cs: 0.05,0.95){AT-2, Order 2, $h=\ell_0/4$};
\foreach \step[count=\i] in {\minIdx,...,\maxIdx}
{
    \pgfmathsetmacro{\col}{100-90*\i/(\maxIdx-\minIdx)}
    \pgfmathsetmacro{\disp}{\minDispl+\i*\stepDispl}
    \edef\temp
    {\noexpand\addplot+[only marks,mark size=0.3pt] table[header=true,col sep=comma, x index=1, y index=2,] {Data/TB/tensile_TB_0.25l0_AT-2_Order2_contours/contour_\step.csv};
    }
    \temp
}

\nextgroupplot
\node[anchor=north west,scale=0.8] at (rel axis cs: 0.05,0.95){AT-2, Order 4, $h=\ell_0/2$};
\foreach \step[count=\i] in {\minIdx,...,\maxIdx}
{
    \pgfmathsetmacro{\col}{100-90*\i/(\maxIdx-\minIdx)}
    \pgfmathsetmacro{\disp}{\minDispl+\i*\stepDispl}
    \edef\temp
    {\noexpand\addplot+[only marks,mark size=0.3pt] table[header=true,col sep=comma, x index=1, y index=2,] {Data/TB/tensile_TB_0.50l0_AT-2_Order4_contours/contour_\step.csv};
    }
    \temp
}

\nextgroupplot[]
\node[anchor=north west,scale=0.8] at (rel axis cs: 0.05,0.95){AT-2, Order 4, $h=\ell_0/4$};
\foreach \step[count=\i] in {\minIdx,...,\maxIdx}
{
    \pgfmathsetmacro{\col}{100-90*\i/(\maxIdx-\minIdx)}
    \pgfmathsetmacro{\disp}{\minDispl+\i*\stepDispl}
    \edef\temp
    {\noexpand\addplot+[only marks,mark size=0.3pt] table[header=true,col sep=comma, x index=1, y index=2,] {Data/TB/tensile_TB_0.25l0_AT-2_Order4_contours/contour_\step.csv};
    }
    \temp
}
\end{groupplot}
\end{tikzpicture}
\caption{Tensile test.}
\label{fig:benchmarks_tensile_contour}
\end{subfigure}

\begin{subfigure}{\linewidth}
\centering
\begin{tikzpicture}
\def\minIdx{11}
\def\maxIdx{20}
\def\minDispl{0.0093}
\def\stepDispl{0.0003}
\pgfmathsetmacro{\maxDispl}{\minDispl+\stepDispl*(\maxIdx-\minIdx)}
\pgfmathsetmacro{\minCol}{10}
\pgfmathsetmacro{\maxCol}{100}
\begin{groupplot}[
        group style={
            group name=myplot,
            group size=4 by 2,
            horizontal sep=5pt,
            vertical sep=5pt,
            xticklabels at=edge bottom,
            xlabels at=edge bottom,
            yticklabels at=edge left,
            ylabels at=edge left
        },
    width=0.375\linewidth,
    axis equal image,
    xlabel={$x$},
    x unit={mm},
    ylabel={$y$},
    y unit={mm},
    xmin = 0.45,xmax = 0.75,
    ymin = 0.15,ymax = 0.55,
    point meta min=\minDispl,
    point meta max=\maxDispl,
    every mark/.append style={solid},
    colormap={mycolormap}{color=(col1!\minCol!col2) color=(col1!\maxCol!col2)},
    cycle list={[samples of colormap={\maxIdx-\minIdx+1} of mycolormap]},
]
\nextgroupplot[]
\node[anchor=north west,scale=0.8] at (rel axis cs: 0.05,0.95){AT-1, Order 2};
\node[anchor=south west,scale=0.8] at (rel axis cs: 0.05,0.05){$h=\ell_0/2$};
\foreach \step[count=\i] in {\minIdx,...,\maxIdx}
{
    \pgfmathsetmacro{\col}{100-90*\i/(\maxIdx-\minIdx)}
    \pgfmathsetmacro{\disp}{\minDispl+\i*\stepDispl}
    \edef\temp
    {\noexpand\addplot+[only marks,mark size=0.3pt] table[header=true,col sep=comma, x index=1, y index=2,]{Data/TB/shear_TB_0.50l0_AT-1_Order2_contours/contour_\step.csv};
    }
    \temp
}

\nextgroupplot[]
\node[anchor=north west,scale=0.8] at (rel axis cs: 0.05,0.95){AT-1, Order 4};
\node[anchor=south west,scale=0.8] at (rel axis cs: 0.05,0.05){$h=\ell_0/2$};
\foreach \step[count=\i] in {\minIdx,...,\maxIdx}
{
    \pgfmathsetmacro{\col}{100-90*\i/(\maxIdx-\minIdx)}
    \pgfmathsetmacro{\disp}{\minDispl+\i*\stepDispl}
    \edef\temp
    {\noexpand\addplot+[only marks,mark size=0.3pt] table[header=true,col sep=comma, x index=1, y index=2,]{Data/TB/shear_TB_0.50l0_AT-1_Order4_contours/contour_\step.csv};
    }
    \temp
}

\nextgroupplot[]
\node[anchor=north west,scale=0.8] at (rel axis cs: 0.05,0.95){AT-2, Order 2};
\node[anchor=south west,scale=0.8] at (rel axis cs: 0.05,0.05){$h=\ell_0/2$};
\foreach \step[count=\i] in {\minIdx,...,\maxIdx}
{
    \pgfmathsetmacro{\col}{100-90*\i/(\maxIdx-\minIdx)}
    \pgfmathsetmacro{\disp}{\minDispl+\i*\stepDispl}
    \edef\temp
    {\noexpand\addplot+[only marks,mark size=0.3pt] table[header=true,col sep=comma, x index=1, y index=2,] {Data/TB/shear_TB_0.50l0_AT-2_Order2_contours/contour_\step.csv};
    }
    \temp
}

\nextgroupplot[colorbar,colorbar style={axis equal image=false,width=0.2cm,x unit={}, y unit={mm},ylabel={Horizontal top displacement},height=
        2*\pgfkeysvalueof{/pgfplots/parent axis height}+
        \pgfkeysvalueof{/pgfplots/group/vertical sep}}]
\node[anchor=north west,scale=0.8] at (rel axis cs: 0.05,0.95){AT-2, Order 4};
\node[anchor=south west,scale=0.8] at (rel axis cs: 0.05,0.05){$h=\ell_0/2$};
\foreach \step[count=\i] in {\minIdx,...,\maxIdx}
{
    \pgfmathsetmacro{\col}{100-90*\i/(\maxIdx-\minIdx)}
    \pgfmathsetmacro{\disp}{\minDispl+\i*\stepDispl}
    \edef\temp
    {\noexpand\addplot+[only marks,mark size=0.3pt] table[header=true,col sep=comma, x index=1, y index=2,] {Data/TB/shear_TB_0.50l0_AT-2_Order4_contours/contour_\step.csv};
    }
    \temp
}

\nextgroupplot[]
\node[anchor=north west,scale=0.8] at (rel axis cs: 0.05,0.95){AT-1, Order 2};
\node[anchor=south west,scale=0.8] at (rel axis cs: 0.05,0.05){$h=\ell_0/4$};
\foreach \step[count=\i] in {\minIdx,...,\maxIdx}
{
    \pgfmathsetmacro{\col}{100-90*\i/(\maxIdx-\minIdx)}
    \pgfmathsetmacro{\disp}{\minDispl+\i*\stepDispl}
    \edef\temp
    {\noexpand\addplot+[only marks,mark size=0.3pt] table[header=true,col sep=comma, x index=1, y index=2,]{Data/TB/shear_TB_0.25l0_AT-1_Order2_contours/contour_\step.csv};
    }
    \temp
}

\nextgroupplot[]
\node[anchor=north west,scale=0.8] at (rel axis cs: 0.05,0.95){AT-1, Order 4};
\node[anchor=south west,scale=0.8] at (rel axis cs: 0.05,0.05){$h=\ell_0/4$};
\foreach \step[count=\i] in {\minIdx,...,\maxIdx}
{
    \pgfmathsetmacro{\col}{100-90*\i/(\maxIdx-\minIdx)}
    \pgfmathsetmacro{\disp}{\minDispl+\i*\stepDispl}
    \edef\temp
    {\noexpand\addplot+[only marks,mark size=0.3pt] table[header=true,col sep=comma, x index=1, y index=2,]{Data/TB/shear_TB_0.25l0_AT-1_Order4_contours/contour_\step.csv};
    }
    \temp
}

\nextgroupplot[]
\node[anchor=north west,scale=0.8] at (rel axis cs: 0.05,0.95){AT-2, Order 2};
\node[anchor=south west,scale=0.8] at (rel axis cs: 0.05,0.05){$h=\ell_0/4$};
\foreach \step[count=\i] in {\minIdx,...,\maxIdx}
{
    \pgfmathsetmacro{\col}{100-90*\i/(\maxIdx-\minIdx)}
    \pgfmathsetmacro{\disp}{\minDispl+\i*\stepDispl}
    \edef\temp
    {\noexpand\addplot+[only marks,mark size=0.3pt] table[header=true,col sep=comma, x index=1, y index=2,] {Data/TB/shear_TB_0.25l0_AT-2_Order2_contours/contour_\step.csv};
    }
    \temp
}

\nextgroupplot[]
\node[anchor=north west,scale=0.8] at (rel axis cs: 0.05,0.95){AT-2, Order 4};
\node[anchor=south west,scale=0.8] at (rel axis cs: 0.05,0.05){$h=\ell_0/4$};
\foreach \step[count=\i] in {\minIdx,...,\maxIdx}
{
    \pgfmathsetmacro{\col}{100-90*\i/(\maxIdx-\minIdx)}
    \pgfmathsetmacro{\disp}{\minDispl+\i*\stepDispl}
    \edef\temp
    {\noexpand\addplot+[only marks,mark size=0.3pt] table[header=true,col sep=comma, x index=1, y index=2,] {Data/TB/shear_TB_0.25l0_AT-2_Order4_contours/contour_\step.csv};
    }
    \temp
}
\end{groupplot}
\end{tikzpicture}
\caption{Shear test.}
\label{fig:benchmarks_shear_contour}
\end{subfigure}
\caption{Contours of the damage field at level $d=0.5$ for the (\subref{fig:benchmarks_tensile_contour}) tensile and (\subref{fig:benchmarks_shear_contour}) tests. The contours are provided for load steps in the interval $u\in[0.93,1.2]\cdot 10^{-3}\:[mm]$ (tensile) and $u\in[0.45,0.60]\cdot 10^{-3}\:[mm]$ (shear) with step size $\Delta u=0.03\cdot 10^{-3}\:[mm]$ as described in \cref{fig:SEN_tensile,fig:SEN_shear}. The rows represent the results obtained on a coarse ($h=\ell_0/2$) and fine ($h=\ell_0/4$) mesh, while the columns represent different fracture models. The color bars are used to indicate solutions in different load steps.}
\label{fig:benchmarks_contour}
\end{figure}

Firstly, \cref{fig:benchmarks_tensile_contour,fig:benchmarks_shear_contour} provide contour plots of the damage field at level $d=0.5$ in all load steps $u\in[0.45,0.600]\cdot 10^{-3}\:[mm]$ for the tensile case, and $u\in[0.93,1.2]\cdot 10^{-3}\:[mm]$ for the shear case. In these figures, the distance between two consecutive contours represents the speed of the crack. Consequently, it can be observed that the crack in the tensile test forms suddenly, typically in 2-3 load steps, while the crack in case of the shear test propagates over a series of time steps towards the end of the simulation. In addition, for both tests it can be observed that the crack for the fourth order AT-1 model is wider in general and that this model generally propagates faster and slightly earlier than the other models.\\

\begin{figure}[h]
\centering
\begin{subfigure}{\linewidth}
\centering
\begin{tikzpicture}
\begin{groupplot}[
        group style={
			group name=myplot,
			group size=4 by 2,
			vertical sep=15pt,
			horizontal sep=5pt,
			xticklabels at=edge bottom,
			xlabels at=edge bottom,
			yticklabels at=edge left,
			ylabels at=edge left
        },
	width=0.3\linewidth,
	height=0.2\textheight,
    xlabel = {$u_y$},
    x label style={at={(axis description cs:0.5,-0.2)},anchor=north},
    x unit={mm},
    legend columns = 6,
]
\nextgroupplot[ymin = 0, ymax = 0.75, 
				ylabel={$F_y$}, y unit={kN/mm}
				]
\node[anchor=north west,scale=0.8] at (rel axis cs: 0.05,0.95){AT1, Order 2};
\addplot+[TP_AT1_Order2_fine	,opacity=1.0] table[header=true,col sep=comma, x=u, y expr=-\thisrow{Fy}] {Data/TB/tensile_TB_0.25l0_AT-1_Order2.csv};
\addplot+[TP_AT1_Order2_coarse	,opacity=1.0] table[header=true,col sep=comma, x=u, y expr=-\thisrow{Fy}] {Data/TB/tensile_TB_0.50l0_AT-1_Order2.csv};
\addplot+[TP_AT1_Order4_fine	,opacity=0.1] table[header=true,col sep=comma, x=u, y expr=-\thisrow{Fy}] {Data/TB/tensile_TB_0.25l0_AT-1_Order4.csv};
\addplot+[TP_AT1_Order4_coarse	,opacity=0.1] table[header=true,col sep=comma, x=u, y expr=-\thisrow{Fy}] {Data/TB/tensile_TB_0.50l0_AT-1_Order4.csv};
\addplot+[TP_AT2_Order2_fine	,opacity=0.1] table[header=true,col sep=comma, x=u, y expr=-\thisrow{Fy}] {Data/TB/tensile_TB_0.25l0_AT-2_Order2.csv};
\addplot+[TP_AT2_Order2_coarse	,opacity=0.1] table[header=true,col sep=comma, x=u, y expr=-\thisrow{Fy}] {Data/TB/tensile_TB_0.50l0_AT-2_Order2.csv};
\addplot+[TP_AT2_Order4_fine	,opacity=0.1] table[header=true,col sep=comma, x=u, y expr=-\thisrow{Fy}] {Data/TB/tensile_TB_0.25l0_AT-2_Order4.csv};
\addplot+[TP_AT2_Order4_coarse	,opacity=0.1] table[header=true,col sep=comma, x=u, y expr=-\thisrow{Fy}] {Data/TB/tensile_TB_0.50l0_AT-2_Order4.csv};

\nextgroupplot[ymin = 0, ymax = 0.75]
\node[anchor=north west,scale=0.8] at (rel axis cs: 0.05,0.95){AT1, Order 4};
\addplot+[TP_AT1_Order2_fine	,opacity=0.1] table[header=true,col sep=comma, x=u, y expr=-\thisrow{Fy}] {Data/TB/tensile_TB_0.25l0_AT-1_Order2.csv};
\addplot+[TP_AT1_Order2_coarse	,opacity=0.1] table[header=true,col sep=comma, x=u, y expr=-\thisrow{Fy}] {Data/TB/tensile_TB_0.50l0_AT-1_Order2.csv};
\addplot+[TP_AT1_Order4_fine	,opacity=1.0] table[header=true,col sep=comma, x=u, y expr=-\thisrow{Fy}] {Data/TB/tensile_TB_0.25l0_AT-1_Order4.csv};
\addplot+[TP_AT1_Order4_coarse	,opacity=1.0] table[header=true,col sep=comma, x=u, y expr=-\thisrow{Fy}] {Data/TB/tensile_TB_0.50l0_AT-1_Order4.csv};
\addplot+[TP_AT2_Order2_fine	,opacity=0.1] table[header=true,col sep=comma, x=u, y expr=-\thisrow{Fy}] {Data/TB/tensile_TB_0.25l0_AT-2_Order2.csv};
\addplot+[TP_AT2_Order2_coarse	,opacity=0.1] table[header=true,col sep=comma, x=u, y expr=-\thisrow{Fy}] {Data/TB/tensile_TB_0.50l0_AT-2_Order2.csv};
\addplot+[TP_AT2_Order4_fine	,opacity=0.1] table[header=true,col sep=comma, x=u, y expr=-\thisrow{Fy}] {Data/TB/tensile_TB_0.25l0_AT-2_Order4.csv};
\addplot+[TP_AT2_Order4_coarse	,opacity=0.1] table[header=true,col sep=comma, x=u, y expr=-\thisrow{Fy}] {Data/TB/tensile_TB_0.50l0_AT-2_Order4.csv};

\nextgroupplot[ymin = 0, ymax = 0.75]
\node[anchor=north west,scale=0.8] at (rel axis cs: 0.05,0.95){AT2, Order 2};
\addplot+[TP_AT1_Order2_fine	,opacity=0.1] table[header=true,col sep=comma, x=u, y expr=-\thisrow{Fy}] {Data/TB/tensile_TB_0.25l0_AT-1_Order2.csv};
\addplot+[TP_AT1_Order2_coarse	,opacity=0.1] table[header=true,col sep=comma, x=u, y expr=-\thisrow{Fy}] {Data/TB/tensile_TB_0.50l0_AT-1_Order2.csv};
\addplot+[TP_AT1_Order4_fine	,opacity=0.1] table[header=true,col sep=comma, x=u, y expr=-\thisrow{Fy}] {Data/TB/tensile_TB_0.25l0_AT-1_Order4.csv};
\addplot+[TP_AT1_Order4_coarse	,opacity=0.1] table[header=true,col sep=comma, x=u, y expr=-\thisrow{Fy}] {Data/TB/tensile_TB_0.50l0_AT-1_Order4.csv};
\addplot+[TP_AT2_Order2_fine	,opacity=1.0] table[header=true,col sep=comma, x=u, y expr=-\thisrow{Fy}] {Data/TB/tensile_TB_0.25l0_AT-2_Order2.csv};
\addplot+[TP_AT2_Order2_coarse	,opacity=1.0] table[header=true,col sep=comma, x=u, y expr=-\thisrow{Fy}] {Data/TB/tensile_TB_0.50l0_AT-2_Order2.csv};
\addplot+[TP_AT2_Order4_fine	,opacity=0.1] table[header=true,col sep=comma, x=u, y expr=-\thisrow{Fy}] {Data/TB/tensile_TB_0.25l0_AT-2_Order4.csv};
\addplot+[TP_AT2_Order4_coarse	,opacity=0.1] table[header=true,col sep=comma, x=u, y expr=-\thisrow{Fy}] {Data/TB/tensile_TB_0.50l0_AT-2_Order4.csv};

\nextgroupplot[ymin = 0, ymax = 0.75]
\node[anchor=north west,scale=0.8] at (rel axis cs: 0.05,0.95){AT2, Order 4};
\addplot+[TP_AT1_Order2_fine	,opacity=0.1] table[header=true,col sep=comma, x=u, y expr=-\thisrow{Fy}] {Data/TB/tensile_TB_0.25l0_AT-1_Order2.csv};
\addplot+[TP_AT1_Order2_coarse	,opacity=0.1] table[header=true,col sep=comma, x=u, y expr=-\thisrow{Fy}] {Data/TB/tensile_TB_0.50l0_AT-1_Order2.csv};
\addplot+[TP_AT1_Order4_fine	,opacity=0.1] table[header=true,col sep=comma, x=u, y expr=-\thisrow{Fy}] {Data/TB/tensile_TB_0.25l0_AT-1_Order4.csv};
\addplot+[TP_AT1_Order4_coarse	,opacity=0.1] table[header=true,col sep=comma, x=u, y expr=-\thisrow{Fy}] {Data/TB/tensile_TB_0.50l0_AT-1_Order4.csv};
\addplot+[TP_AT2_Order2_fine	,opacity=0.1] table[header=true,col sep=comma, x=u, y expr=-\thisrow{Fy}] {Data/TB/tensile_TB_0.25l0_AT-2_Order2.csv};
\addplot+[TP_AT2_Order2_coarse	,opacity=0.1] table[header=true,col sep=comma, x=u, y expr=-\thisrow{Fy}] {Data/TB/tensile_TB_0.50l0_AT-2_Order2.csv};
\addplot+[TP_AT2_Order4_fine	,opacity=1.0] table[header=true,col sep=comma, x=u, y expr=-\thisrow{Fy}] {Data/TB/tensile_TB_0.25l0_AT-2_Order4.csv};
\addplot+[TP_AT2_Order4_coarse	,opacity=1.0] table[header=true,col sep=comma, x=u, y expr=-\thisrow{Fy}] {Data/TB/tensile_TB_0.50l0_AT-2_Order4.csv};

\nextgroupplot[	ymin = 1.3e-3, ymax = 4.0e-3, 
				ylabel={$\mathcal{D}$}, y unit={kJ/mm},
				legend to name={CommonLegend},]
\addlegendimage{TP_AT1_Order2_coarse,black}\addlegendentry{Coarse} 
\addlegendimage{TP_AT2_Order2_fine,black}\addlegendentry{Fine} 

\addlegendimage{TP_AT1_Order2_coarse,solid,no markers}\addlegendentry{AT1, Order 2}
\addlegendimage{TP_AT1_Order4_coarse,solid,no markers}\addlegendentry{AT1, Order 4}
\addlegendimage{TP_AT2_Order2_coarse,solid,no markers}\addlegendentry{AT2, Order 2}
\addlegendimage{TP_AT2_Order4_coarse,solid,no markers}\addlegendentry{AT2, Order 4}

\draw[black!50] (axis cs: 0e-3,1.35e-3) -- (axis cs: 6e-3,1.35e-3);
\node[anchor=south east,scale=0.8] at (axis cs: 6e-3,1.35e-3) {$1.35$};
\node[anchor=north west,scale=0.8] at (rel axis cs: 0.05,0.95){AT1, Order 2};
\addplot+[TP_AT1_Order2_fine	,opacity=1.0] table[header=true,col sep=comma, x=u, y=E_d] {Data/TB/tensile_TB_0.25l0_AT-1_Order2.csv};
\addplot+[TP_AT1_Order2_coarse	,opacity=1.0] table[header=true,col sep=comma, x=u, y=E_d] {Data/TB/tensile_TB_0.50l0_AT-1_Order2.csv};
\addplot+[TP_AT1_Order4_fine	,opacity=0.2] table[header=true,col sep=comma, x=u, y=E_d] {Data/TB/tensile_TB_0.25l0_AT-1_Order4.csv};
\addplot+[TP_AT1_Order4_coarse	,opacity=0.2] table[header=true,col sep=comma, x=u, y=E_d] {Data/TB/tensile_TB_0.50l0_AT-1_Order4.csv};
\addplot+[TP_AT2_Order2_fine	,opacity=0.2] table[header=true,col sep=comma, x=u, y=E_d] {Data/TB/tensile_TB_0.25l0_AT-2_Order2.csv};
\addplot+[TP_AT2_Order2_coarse	,opacity=0.2] table[header=true,col sep=comma, x=u, y=E_d] {Data/TB/tensile_TB_0.50l0_AT-2_Order2.csv};
\addplot+[TP_AT2_Order4_fine	,opacity=0.2] table[header=true,col sep=comma, x=u, y=E_d] {Data/TB/tensile_TB_0.25l0_AT-2_Order4.csv};
\addplot+[TP_AT2_Order4_coarse	,opacity=0.2] table[header=true,col sep=comma, x=u, y=E_d] {Data/TB/tensile_TB_0.50l0_AT-2_Order4.csv};

\nextgroupplot[ymin = 1.3e-3, ymax = 4.0e-3]
\draw[black!50] (axis cs: 0e-3,1.35e-3) -- (axis cs: 6e-3,1.35e-3);
\node[anchor=south east,scale=0.8] at (axis cs: 6e-3,1.35e-3) {$1.35$};
\node[anchor=north west,scale=0.8] at (rel axis cs: 0.05,0.95){AT1, Order 4};
\addplot+[TP_AT1_Order2_fine	,opacity=0.2] table[header=true,col sep=comma, x=u, y=E_d] {Data/TB/tensile_TB_0.25l0_AT-1_Order2.csv};
\addplot+[TP_AT1_Order2_coarse	,opacity=0.2] table[header=true,col sep=comma, x=u, y=E_d] {Data/TB/tensile_TB_0.50l0_AT-1_Order2.csv};
\addplot+[TP_AT1_Order4_fine	,opacity=1.0] table[header=true,col sep=comma, x=u, y=E_d] {Data/TB/tensile_TB_0.25l0_AT-1_Order4.csv};
\addplot+[TP_AT1_Order4_coarse	,opacity=1.0] table[header=true,col sep=comma, x=u, y=E_d] {Data/TB/tensile_TB_0.50l0_AT-1_Order4.csv};
\addplot+[TP_AT2_Order2_fine	,opacity=0.2] table[header=true,col sep=comma, x=u, y=E_d] {Data/TB/tensile_TB_0.25l0_AT-2_Order2.csv};
\addplot+[TP_AT2_Order2_coarse	,opacity=0.2] table[header=true,col sep=comma, x=u, y=E_d] {Data/TB/tensile_TB_0.50l0_AT-2_Order2.csv};
\addplot+[TP_AT2_Order4_fine	,opacity=0.2] table[header=true,col sep=comma, x=u, y=E_d] {Data/TB/tensile_TB_0.25l0_AT-2_Order4.csv};
\addplot+[TP_AT2_Order4_coarse	,opacity=0.2] table[header=true,col sep=comma, x=u, y=E_d] {Data/TB/tensile_TB_0.50l0_AT-2_Order4.csv};

\nextgroupplot[ymin = 1.3e-3, ymax = 4.0e-3]
\draw[black!50] (axis cs: 0e-3,1.35e-3) -- (axis cs: 6e-3,1.35e-3);
\node[anchor=south east,scale=0.8] at (axis cs: 6e-3,1.35e-3) {$1.35$};
\node[anchor=north west,scale=0.8] at (rel axis cs: 0.05,0.95){AT2, Order 2};
\addplot+[TP_AT1_Order2_fine	,opacity=0.2] table[header=true,col sep=comma, x=u, y=E_d] {Data/TB/tensile_TB_0.25l0_AT-1_Order2.csv};
\addplot+[TP_AT1_Order2_coarse	,opacity=0.2] table[header=true,col sep=comma, x=u, y=E_d] {Data/TB/tensile_TB_0.50l0_AT-1_Order2.csv};
\addplot+[TP_AT1_Order4_fine	,opacity=0.2] table[header=true,col sep=comma, x=u, y=E_d] {Data/TB/tensile_TB_0.25l0_AT-1_Order4.csv};
\addplot+[TP_AT1_Order4_coarse	,opacity=0.2] table[header=true,col sep=comma, x=u, y=E_d] {Data/TB/tensile_TB_0.50l0_AT-1_Order4.csv};
\addplot+[TP_AT2_Order2_fine	,opacity=1.0] table[header=true,col sep=comma, x=u, y=E_d] {Data/TB/tensile_TB_0.25l0_AT-2_Order2.csv};
\addplot+[TP_AT2_Order2_coarse	,opacity=1.0] table[header=true,col sep=comma, x=u, y=E_d] {Data/TB/tensile_TB_0.50l0_AT-2_Order2.csv};
\addplot+[TP_AT2_Order4_fine	,opacity=0.2] table[header=true,col sep=comma, x=u, y=E_d] {Data/TB/tensile_TB_0.25l0_AT-2_Order4.csv};
\addplot+[TP_AT2_Order4_coarse	,opacity=0.2] table[header=true,col sep=comma, x=u, y=E_d] {Data/TB/tensile_TB_0.50l0_AT-2_Order4.csv};

\nextgroupplot[ymin = 1.3e-3, ymax = 4.0e-3]
\draw[black!50] (axis cs: 0e-3,1.35e-3) -- (axis cs: 6e-3,1.35e-3);
\node[anchor=south east,scale=0.8] at (axis cs: 6e-3,1.35e-3) {$1.35$};
\node[anchor=north west,scale=0.8] at (rel axis cs: 0.05,0.95){AT2, Order 4};
\addplot+[TP_AT1_Order2_fine	,opacity=0.2] table[header=true,col sep=comma, x=u, y=E_d] {Data/TB/tensile_TB_0.25l0_AT-1_Order2.csv};
\addplot+[TP_AT1_Order2_coarse	,opacity=0.2] table[header=true,col sep=comma, x=u, y=E_d] {Data/TB/tensile_TB_0.50l0_AT-1_Order2.csv};
\addplot+[TP_AT1_Order4_fine	,opacity=0.2] table[header=true,col sep=comma, x=u, y=E_d] {Data/TB/tensile_TB_0.25l0_AT-1_Order4.csv};
\addplot+[TP_AT1_Order4_coarse	,opacity=0.2] table[header=true,col sep=comma, x=u, y=E_d] {Data/TB/tensile_TB_0.50l0_AT-1_Order4.csv};
\addplot+[TP_AT2_Order2_fine	,opacity=0.2] table[header=true,col sep=comma, x=u, y=E_d] {Data/TB/tensile_TB_0.25l0_AT-2_Order2.csv};
\addplot+[TP_AT2_Order2_coarse	,opacity=0.2] table[header=true,col sep=comma, x=u, y=E_d] {Data/TB/tensile_TB_0.50l0_AT-2_Order2.csv};
\addplot+[TP_AT2_Order4_fine	,opacity=1.0] table[header=true,col sep=comma, x=u, y=E_d] {Data/TB/tensile_TB_0.25l0_AT-2_Order4.csv};
\addplot+[TP_AT2_Order4_coarse	,opacity=1.0] table[header=true,col sep=comma, x=u, y=E_d] {Data/TB/tensile_TB_0.50l0_AT-2_Order4.csv};
\end{groupplot}
\path (myplot c2r1.north east) -- node[above]{\pgfplotslegendfromname{CommonLegend}} (myplot c3r1.north west);

\end{tikzpicture}
\caption{Tensile test.}
\label{fig:benchmarks_tensile_modelcomparison}
\end{subfigure}

\begin{subfigure}{\linewidth}
\centering
\begin{tikzpicture}
\begin{groupplot}[
        group style={
			group name=myplot,
			group size=4 by 2,
			vertical sep=15pt,
			horizontal sep=5pt,
			xticklabels at=edge bottom,
			xlabels at=edge bottom,
			yticklabels at=edge left,
			ylabels at=edge left
        },
	width=0.3\linewidth,
	height=0.2\textheight,
    xlabel = {$u_x$},
    x label style={at={(axis description cs:0.5,-0.2)},anchor=north},
    x unit={mm},
    legend columns = 6,
]
\nextgroupplot[ymin = 0.25, ymax = 0.45, 
				ylabel={$F_x$}, y unit={kN/mm}
				]
\node[anchor=north west,scale=0.8] at (rel axis cs: 0.05,0.95) {AT1, Order 2};
\addplot+[TP_AT1_Order2_fine	,opacity=1.0] table[header=true,col sep=comma, x=u, y expr=-\thisrow{Fx}] {Data/TB/shear_TB_0.25l0_AT-1_Order2.csv};
\addplot+[TP_AT1_Order2_coarse	,opacity=1.0] table[header=true,col sep=comma, x=u, y expr=-\thisrow{Fx}] {Data/TB/shear_TB_0.50l0_AT-1_Order2.csv};
\addplot+[TP_AT1_Order4_fine	,opacity=0.1] table[header=true,col sep=comma, x=u, y expr=-\thisrow{Fx}] {Data/TB/shear_TB_0.25l0_AT-1_Order4.csv};
\addplot+[TP_AT1_Order4_coarse	,opacity=0.1] table[header=true,col sep=comma, x=u, y expr=-\thisrow{Fx}] {Data/TB/shear_TB_0.50l0_AT-1_Order4.csv};
\addplot+[TP_AT2_Order2_fine	,opacity=0.1] table[header=true,col sep=comma, x=u, y expr=-\thisrow{Fx}] {Data/TB/shear_TB_0.25l0_AT-2_Order2.csv};
\addplot+[TP_AT2_Order2_coarse	,opacity=0.1] table[header=true,col sep=comma, x=u, y expr=-\thisrow{Fx}] {Data/TB/shear_TB_0.50l0_AT-2_Order2.csv};
\addplot+[TP_AT2_Order4_fine	,opacity=0.1] table[header=true,col sep=comma, x=u, y expr=-\thisrow{Fx}] {Data/TB/shear_TB_0.25l0_AT-2_Order4.csv};
\addplot+[TP_AT2_Order4_coarse	,opacity=0.1] table[header=true,col sep=comma, x=u, y expr=-\thisrow{Fx}] {Data/TB/shear_TB_0.50l0_AT-2_Order4.csv};

\nextgroupplot[ymin = 0.25, ymax = 0.45]
\node[anchor=north west,scale=0.8] at (rel axis cs: 0.05,0.95){AT1, Order 4};
\addplot+[TP_AT1_Order2_fine	,opacity=0.1] table[header=true,col sep=comma, x=u, y expr=-\thisrow{Fx}] {Data/TB/shear_TB_0.25l0_AT-1_Order2.csv};
\addplot+[TP_AT1_Order2_coarse	,opacity=0.1] table[header=true,col sep=comma, x=u, y expr=-\thisrow{Fx}] {Data/TB/shear_TB_0.50l0_AT-1_Order2.csv};
\addplot+[TP_AT1_Order4_fine	,opacity=1.0] table[header=true,col sep=comma, x=u, y expr=-\thisrow{Fx}] {Data/TB/shear_TB_0.25l0_AT-1_Order4.csv};
\addplot+[TP_AT1_Order4_coarse	,opacity=1.0] table[header=true,col sep=comma, x=u, y expr=-\thisrow{Fx}] {Data/TB/shear_TB_0.50l0_AT-1_Order4.csv};
\addplot+[TP_AT2_Order2_fine	,opacity=0.1] table[header=true,col sep=comma, x=u, y expr=-\thisrow{Fx}] {Data/TB/shear_TB_0.25l0_AT-2_Order2.csv};
\addplot+[TP_AT2_Order2_coarse	,opacity=0.1] table[header=true,col sep=comma, x=u, y expr=-\thisrow{Fx}] {Data/TB/shear_TB_0.50l0_AT-2_Order2.csv};
\addplot+[TP_AT2_Order4_fine	,opacity=0.1] table[header=true,col sep=comma, x=u, y expr=-\thisrow{Fx}] {Data/TB/shear_TB_0.25l0_AT-2_Order4.csv};
\addplot+[TP_AT2_Order4_coarse	,opacity=0.1] table[header=true,col sep=comma, x=u, y expr=-\thisrow{Fx}] {Data/TB/shear_TB_0.50l0_AT-2_Order4.csv};

\nextgroupplot[ymin = 0.25, ymax = 0.45]
\node[anchor=north west,scale=0.8] at (rel axis cs: 0.05,0.95){AT2, Order 2};
\addplot+[TP_AT1_Order2_fine	,opacity=0.1] table[header=true,col sep=comma, x=u, y expr=-\thisrow{Fx}] {Data/TB/shear_TB_0.25l0_AT-1_Order2.csv};
\addplot+[TP_AT1_Order2_coarse	,opacity=0.1] table[header=true,col sep=comma, x=u, y expr=-\thisrow{Fx}] {Data/TB/shear_TB_0.50l0_AT-1_Order2.csv};
\addplot+[TP_AT1_Order4_fine	,opacity=0.1] table[header=true,col sep=comma, x=u, y expr=-\thisrow{Fx}] {Data/TB/shear_TB_0.25l0_AT-1_Order4.csv};
\addplot+[TP_AT1_Order4_coarse	,opacity=0.1] table[header=true,col sep=comma, x=u, y expr=-\thisrow{Fx}] {Data/TB/shear_TB_0.50l0_AT-1_Order4.csv};
\addplot+[TP_AT2_Order2_fine	,opacity=1.0] table[header=true,col sep=comma, x=u, y expr=-\thisrow{Fx}] {Data/TB/shear_TB_0.25l0_AT-2_Order2.csv};
\addplot+[TP_AT2_Order2_coarse	,opacity=1.0] table[header=true,col sep=comma, x=u, y expr=-\thisrow{Fx}] {Data/TB/shear_TB_0.50l0_AT-2_Order2.csv};
\addplot+[TP_AT2_Order4_fine	,opacity=0.1] table[header=true,col sep=comma, x=u, y expr=-\thisrow{Fx}] {Data/TB/shear_TB_0.25l0_AT-2_Order4.csv};
\addplot+[TP_AT2_Order4_coarse	,opacity=0.1] table[header=true,col sep=comma, x=u, y expr=-\thisrow{Fx}] {Data/TB/shear_TB_0.50l0_AT-2_Order4.csv};

\nextgroupplot[ymin = 0.25, ymax = 0.45]
\node[anchor=north west,scale=0.8] at (rel axis cs: 0.05,0.95){AT2, Order 4};
\addplot+[TP_AT1_Order2_fine	,opacity=0.1] table[header=true,col sep=comma, x=u, y expr=-\thisrow{Fx}] {Data/TB/shear_TB_0.25l0_AT-1_Order2.csv};
\addplot+[TP_AT1_Order2_coarse	,opacity=0.1] table[header=true,col sep=comma, x=u, y expr=-\thisrow{Fx}] {Data/TB/shear_TB_0.50l0_AT-1_Order2.csv};
\addplot+[TP_AT1_Order4_fine	,opacity=0.1] table[header=true,col sep=comma, x=u, y expr=-\thisrow{Fx}] {Data/TB/shear_TB_0.25l0_AT-1_Order4.csv};
\addplot+[TP_AT1_Order4_coarse	,opacity=0.1] table[header=true,col sep=comma, x=u, y expr=-\thisrow{Fx}] {Data/TB/shear_TB_0.50l0_AT-1_Order4.csv};
\addplot+[TP_AT2_Order2_fine	,opacity=0.1] table[header=true,col sep=comma, x=u, y expr=-\thisrow{Fx}] {Data/TB/shear_TB_0.25l0_AT-2_Order2.csv};
\addplot+[TP_AT2_Order2_coarse	,opacity=0.1] table[header=true,col sep=comma, x=u, y expr=-\thisrow{Fx}] {Data/TB/shear_TB_0.50l0_AT-2_Order2.csv};
\addplot+[TP_AT2_Order4_fine	,opacity=1.0] table[header=true,col sep=comma, x=u, y expr=-\thisrow{Fx}] {Data/TB/shear_TB_0.25l0_AT-2_Order4.csv};
\addplot+[TP_AT2_Order4_coarse	,opacity=1.0] table[header=true,col sep=comma, x=u, y expr=-\thisrow{Fx}] {Data/TB/shear_TB_0.50l0_AT-2_Order4.csv};

\nextgroupplot[	ymin = 1e-3, ymax = 3.0e-3, 
				ylabel={$\mathcal{D}$}, y unit={kJ/mm},
				legend to name={CommonLegend},]
\addlegendimage{TP_AT1_Order2_coarse,black}\addlegendentry{Coarse} 
\addlegendimage{TP_AT2_Order2_fine,black}\addlegendentry{Fine} 

\addlegendimage{TP_AT1_Order2_coarse,solid,no markers}\addlegendentry{AT1, Order 2}
\addlegendimage{TP_AT1_Order4_coarse,solid,no markers}\addlegendentry{AT1, Order 4}
\addlegendimage{TP_AT2_Order2_coarse,solid,no markers}\addlegendentry{AT2, Order 2}
\addlegendimage{TP_AT2_Order4_coarse,solid,no markers}\addlegendentry{AT2, Order 4}

\draw[black!50] (axis cs: 6e-3,1.35e-3) -- (axis cs: 12e-3,1.35e-3);
\node[anchor=south east,scale=0.8] at (axis cs: 12e-3,1.35e-3) {$1.35$};
\node[anchor=north west,scale=0.8] at (rel axis cs: 0.05,0.95){AT1, Order 2};
\addplot+[TP_AT1_Order2_fine	,opacity=1.0] table[header=true,col sep=comma, x=u, y=E_d] {Data/TB/shear_TB_0.25l0_AT-1_Order2.csv};
\addplot+[TP_AT1_Order2_coarse	,opacity=1.0] table[header=true,col sep=comma, x=u, y=E_d] {Data/TB/shear_TB_0.50l0_AT-1_Order2.csv};
\addplot+[TP_AT1_Order4_fine	,opacity=0.1] table[header=true,col sep=comma, x=u, y=E_d] {Data/TB/shear_TB_0.25l0_AT-1_Order4.csv};
\addplot+[TP_AT1_Order4_coarse	,opacity=0.1] table[header=true,col sep=comma, x=u, y=E_d] {Data/TB/shear_TB_0.50l0_AT-1_Order4.csv};
\addplot+[TP_AT2_Order2_fine	,opacity=0.1] table[header=true,col sep=comma, x=u, y=E_d] {Data/TB/shear_TB_0.25l0_AT-2_Order2.csv};
\addplot+[TP_AT2_Order2_coarse	,opacity=0.1] table[header=true,col sep=comma, x=u, y=E_d] {Data/TB/shear_TB_0.50l0_AT-2_Order2.csv};
\addplot+[TP_AT2_Order4_fine	,opacity=0.1] table[header=true,col sep=comma, x=u, y=E_d] {Data/TB/shear_TB_0.25l0_AT-2_Order4.csv};
\addplot+[TP_AT2_Order4_coarse	,opacity=0.1] table[header=true,col sep=comma, x=u, y=E_d] {Data/TB/shear_TB_0.50l0_AT-2_Order4.csv};

\nextgroupplot[ymin = 1e-3, ymax = 3.0e-3]
\draw[black!50] (axis cs: 6e-3,1.35e-3) -- (axis cs: 12e-3,1.35e-3);
\node[anchor=south east,scale=0.8] at (axis cs: 12e-3,1.35e-3) {$1.35$};
\node[anchor=north west,scale=0.8] at (rel axis cs: 0.05,0.95){AT1, Order 4};
\addplot+[TP_AT1_Order2_fine	,opacity=0.1] table[header=true,col sep=comma, x=u, y=E_d] {Data/TB/shear_TB_0.25l0_AT-1_Order2.csv};
\addplot+[TP_AT1_Order2_coarse	,opacity=0.1] table[header=true,col sep=comma, x=u, y=E_d] {Data/TB/shear_TB_0.50l0_AT-1_Order2.csv};
\addplot+[TP_AT1_Order4_fine	,opacity=1.0] table[header=true,col sep=comma, x=u, y=E_d] {Data/TB/shear_TB_0.25l0_AT-1_Order4.csv};
\addplot+[TP_AT1_Order4_coarse	,opacity=1.0] table[header=true,col sep=comma, x=u, y=E_d] {Data/TB/shear_TB_0.50l0_AT-1_Order4.csv};
\addplot+[TP_AT2_Order2_fine	,opacity=0.1] table[header=true,col sep=comma, x=u, y=E_d] {Data/TB/shear_TB_0.25l0_AT-2_Order2.csv};
\addplot+[TP_AT2_Order2_coarse	,opacity=0.1] table[header=true,col sep=comma, x=u, y=E_d] {Data/TB/shear_TB_0.50l0_AT-2_Order2.csv};
\addplot+[TP_AT2_Order4_fine	,opacity=0.1] table[header=true,col sep=comma, x=u, y=E_d] {Data/TB/shear_TB_0.25l0_AT-2_Order4.csv};
\addplot+[TP_AT2_Order4_coarse	,opacity=0.1] table[header=true,col sep=comma, x=u, y=E_d] {Data/TB/shear_TB_0.50l0_AT-2_Order4.csv};

\nextgroupplot[ymin = 1e-3, ymax = 3.0e-3]
\draw[black!50] (axis cs: 6e-3,1.35e-3) -- (axis cs: 12e-3,1.35e-3);
\node[anchor=south east,scale=0.8] at (axis cs: 12e-3,1.35e-3) {$1.35$};
\node[anchor=north west,scale=0.8] at (rel axis cs: 0.05,0.95){AT2, Order 2};
\addplot+[TP_AT1_Order2_fine	,opacity=0.1] table[header=true,col sep=comma, x=u, y=E_d] {Data/TB/shear_TB_0.25l0_AT-1_Order2.csv};
\addplot+[TP_AT1_Order2_coarse	,opacity=0.1] table[header=true,col sep=comma, x=u, y=E_d] {Data/TB/shear_TB_0.50l0_AT-1_Order2.csv};
\addplot+[TP_AT1_Order4_fine	,opacity=0.1] table[header=true,col sep=comma, x=u, y=E_d] {Data/TB/shear_TB_0.25l0_AT-1_Order4.csv};
\addplot+[TP_AT1_Order4_coarse	,opacity=0.1] table[header=true,col sep=comma, x=u, y=E_d] {Data/TB/shear_TB_0.50l0_AT-1_Order4.csv};
\addplot+[TP_AT2_Order2_fine	,opacity=1.0] table[header=true,col sep=comma, x=u, y=E_d] {Data/TB/shear_TB_0.25l0_AT-2_Order2.csv};
\addplot+[TP_AT2_Order2_coarse	,opacity=1.0] table[header=true,col sep=comma, x=u, y=E_d] {Data/TB/shear_TB_0.50l0_AT-2_Order2.csv};
\addplot+[TP_AT2_Order4_fine	,opacity=0.1] table[header=true,col sep=comma, x=u, y=E_d] {Data/TB/shear_TB_0.25l0_AT-2_Order4.csv};
\addplot+[TP_AT2_Order4_coarse	,opacity=0.1] table[header=true,col sep=comma, x=u, y=E_d] {Data/TB/shear_TB_0.50l0_AT-2_Order4.csv};

\nextgroupplot[ymin = 1e-3, ymax = 3.0e-3]
\draw[black!50] (axis cs: 6e-3,1.35e-3) -- (axis cs: 12e-3,1.35e-3);
\node[anchor=south east,scale=0.8] at (axis cs: 12e-3,1.35e-3) {$1.35$};
\node[anchor=north west,scale=0.8] at (rel axis cs: 0.05,0.95){AT2, Order 4};
\addplot+[TP_AT1_Order2_fine	,opacity=0.1] table[header=true,col sep=comma, x=u, y=E_d] {Data/TB/shear_TB_0.25l0_AT-1_Order2.csv};
\addplot+[TP_AT1_Order2_coarse	,opacity=0.1] table[header=true,col sep=comma, x=u, y=E_d] {Data/TB/shear_TB_0.50l0_AT-1_Order2.csv};
\addplot+[TP_AT1_Order4_fine	,opacity=0.1] table[header=true,col sep=comma, x=u, y=E_d] {Data/TB/shear_TB_0.25l0_AT-1_Order4.csv};
\addplot+[TP_AT1_Order4_coarse	,opacity=0.1] table[header=true,col sep=comma, x=u, y=E_d] {Data/TB/shear_TB_0.50l0_AT-1_Order4.csv};
\addplot+[TP_AT2_Order2_fine	,opacity=0.1] table[header=true,col sep=comma, x=u, y=E_d] {Data/TB/shear_TB_0.25l0_AT-2_Order2.csv};
\addplot+[TP_AT2_Order2_coarse	,opacity=0.1] table[header=true,col sep=comma, x=u, y=E_d] {Data/TB/shear_TB_0.50l0_AT-2_Order2.csv};
\addplot+[TP_AT2_Order4_fine	,opacity=1.0] table[header=true,col sep=comma, x=u, y=E_d] {Data/TB/shear_TB_0.25l0_AT-2_Order4.csv};
\addplot+[TP_AT2_Order4_coarse	,opacity=1.0] table[header=true,col sep=comma, x=u, y=E_d] {Data/TB/shear_TB_0.50l0_AT-2_Order4.csv};
\end{groupplot}
\path (myplot c2r1.north east) -- node[above]{\pgfplotslegendfromname{CommonLegend}} (myplot c3r1.north west);

\end{tikzpicture}
\caption{Shear test.}
\label{fig:benchmarks_shear_modelcomparison}
\end{subfigure}
\caption{Reaction force $F_x$, $F_y$ (top) and dissipated energy $\mathcal{D}$ (bottom) for the (\subref{fig:benchmarks_tensile_modelcomparison}) tensile and (\subref{fig:benchmarks_shear_modelcomparison}) shear tests with respect to the displacement at the top boundary. {\NEW Each colun denotes another phase-field model, with the results of the other phase-field models plotted lightly to provide a comparison across model. In addition, the line labeled with $1.35$ denotes the value of the initial dissipation $G_c \cdot L_{\text{initial crack}}$.}}
\label{fig:benchmarks_modelcomparison}
\end{figure}

Besides the contour plots in \cref{fig:benchmarks_tensile_contour,fig:benchmarks_tensile_contour}, \cref{fig:benchmarks_tensile_modelcomparison,fig:benchmarks_tensile_modelcomparison} provide the vertical and horizontal reaction forces, $F_y$ and $F_x$ for the tensile and shear cases, respectively, as well as the dissipated energy $\mathcal{D}$. These plots are provided for all combinations of the AT1 and AT2 models of second and fourth order, on mesh sizes $h=\ell_0/2$ and $h=\ell_0/4$.\\

Besides the reaction forces being a global metric from a mechanical perspective, the dissipated energy provides a global metric on the phase-field: the jump between the dissipated energy between the first and last load step provides a measure for the accuracy of the model. In addition, the accuracy of the phase-field initialization can be noted in the graphs of the dissipated energy: in the undeformed configuration the dissipated energy should start at
\begin{equation} 
\int_{\Gamma} G_c \, d\Gamma = G_c \cdot L_{\text{initial crack}} = 1.35\cdot 10^{-3} \, \text{[kJ/mm]}.
\end{equation} 
As observed from the graph, a higher-order model or a finer mesh size provides more accurate phase-field initialization, hence more accurate initial dissipated energy. In addition, as can be seen in \cref{fig:benchmarks_tensile_modelcomparison}, the fourth-order AT1 model is the only model providing consistent reaction forces between the coarse and the fine mesh, meaning that it is the only model for which the coarse mesh can be used to obtain accurate results. This is in line with the observations from literature \cite{greco2024higher}. Considering the results of the shear test (see \cref{fig:benchmarks_shear_modelcomparison}), a similar conclusion can be drawn: all models except the fourth-order AT1 show a significant difference between the dissipation on the coarse and fine meshes, indicating that only for the fourth-order AT1 model the coarse mesh suffices. A similar conclusion was drawn in \cite{greco2024higher}.\\

\subsection{The role of adaptive meshing}
\label{subsec:adaptivity}
Following upon the comparison of the phase-field models, this section provides an assessment of the role of the adaptive meshing algorithm presented in this paper, applied to the problems presented in \cref{fig:SEN}. To this end, the benchmark problems are solved by means of \cref{alg:adaptive_load_step} in \ref{app:algorithms}.\\

Firstly, the contours of the propagating damage field using the fourth-order AT1 model are provided in \cref{fig:implicitexplicit_tensile,fig:implicitexplicit_shear}. From these figures, it can be seen that the implicit and hybrid adaptive meshing simulations provide an accurate reproduction of the results obtained using the tensor-product B-spline basis. In addition, it can be seen that the results obtained using the explicit scheme lag behind the tensor-product results. The latter is explained by the fact that the mesh is not sufficiently refined to resolve the crack in the next load step. Through refinement iterations, this can be resolved, as shown by the implicit and hybrid schemes.\\

Secondly, \cref{fig:implicitexplicit_plots_tensile,fig:implicitexplicit_plots_shear} provide the dissipated energy $\mathcal{D}$, the CPU-time and the number of degrees of freedom per load step for all phase-field models considered, on the finest mesh with mesh size $h=\ell_0/4$. The coarse mesh size is omitted for the sake of compactness. As can be seen from the dissipated energy distribution and the number of degree of freedom over time, no significant difference between the tensor-product and the implicit and hybrid approaches are observed. The explicit refinement approach, however, shows that the damage progresses slower, which is visible in the plots of the dissipated energy $\mathcal{D}$ as well as the degrees of freedom plot. In addition, the results show that the adaptive approaches are significantly more efficient than the tensor-product simulation, even when the implicit and hybrid approaches re-compute a load step. In addition, the CPU-time plot corresponding to the hybrid scheme shows equivalent CPU times to the explicit scheme in the pre-crack regime, while in the crack propagation regime its CPU times follow the same pattern as the implicit simulations.

\begin{figure}
\centering
\begin{subfigure}{\linewidth}
	\centering
	\begin{tikzpicture}
		\def\minIdx{15}
		\def\maxIdx{20}
		\def\minDispl{0.45} 
		\def\stepDispl{0.03}
		\pgfmathsetmacro{\maxDispl}{\minDispl+\stepDispl*(\maxIdx-\minIdx)}
		\begin{groupplot}[
			group style={
				group name=myplot,
				group size=2 by 3,
				horizontal sep=5pt,
				vertical sep=5pt,
				xticklabels at=edge bottom,
				xlabels at=edge bottom,
				yticklabels at=edge left,
				ylabels at=edge left
			},
			width=0.55\linewidth,
			axis equal image,
			xlabel={$x$},
			x unit={mm},
			ylabel={$y$},
			y unit={mm},
			xmin = 0.45,xmax = 1.0,
			ymin = 0.45,ymax = 0.55,
			cycle list name=MyCycleList,
			every mark/.style={solid},
			legend to name=zelda,
			legend columns = 5,
			legend style={/tikz/every even column/.append style={column sep=0.25cm}}
			]

			\pgfplotsforeachungrouped \step in {\minIdx,...,\maxIdx}%
			{
				\pgfmathsetmacro{\disp}{\minDispl+\stepDispl*(\step-\minIdx)}
				\def\AT{1}
				\def\ORDER{4}
				\def\MESH{0.50}
				\def\MESHl{\MESH l0}
				\edef\tmp
				{
					\noexpand\nextgroupplot[
					]
					\noexpand\addlegendimage{black!50,opacity=0.5,mark=*,only marks,mark size=1pt,solid};
					\noexpand\addlegendentry{TP};
					\noexpand\addlegendimage{col1,mark=*,only marks,mark size=1pt,solid};
					\noexpand\addlegendentry{Implicit};
					\noexpand\addlegendimage{col2,mark=*,only marks,mark size=1pt,solid};
					\noexpand\addlegendentry{Explicit};
					\noexpand\addlegendimage{col3,mark=*,only marks,mark size=1pt,solid};
					\noexpand\addlegendentry{Hybrid};
					\noexpand\addplot+[black!50,opacity=0.5,mark=*,only marks,mark size=0.15pt,solid] table[header=true,col sep=comma, x index=1, y index=2,] {Data/TB/tensile_TB_0.50l0_AT-\AT_Order\ORDER_contours/contour_\step.csv};

					\noexpand\addplot+[col1,mark=*,only marks,mark size=0.3pt] table[header=true,col sep=comma, x index=1, y index=2,] {Data/implicit/tensile_THB_4l0_\MESHl_AT-\AT_Order\ORDER_contours/contour_\step.csv};
					\noexpand\addplot+[col2,mark=*,only marks,mark size=0.3pt] table[header=true,col sep=comma, x index=1, y index=2,] {Data/explicit/tensile_THB_4l0_\MESHl_AT-\AT_Order\ORDER_contours/contour_\step.csv};
					\noexpand\addplot+[col3,mark=*,only marks,mark size=0.3pt] table[header=true,col sep=comma, x index=1, y index=2,] {Data/hybrid/tensile_THB_4l0_\MESHl_AT-\AT_Order\ORDER_contours/contour_\step.csv};
					\noexpand\node[anchor=north west,scale=0.8] at (rel axis cs: 0,1) {$u=\pgfmathprintnumber[precision=3]{\disp}\cdot 10^{-3}$};
					}
				\tmp
			}
		\end{groupplot}
		\node[anchor=south] (legend) at ($(myplot c1r1.north east)!0.5!(myplot c2r1.north west)$)
		{\pgfplotslegendfromname{zelda}};
	\end{tikzpicture}
	\caption{Tensile test.}
	\label{fig:implicitexplicit_tensile}
\end{subfigure}

\begin{subfigure}{\linewidth}
	\centering
	\begin{tikzpicture}
		\def\minIdx{11}
		\def\maxIdx{20}
		\def\minDispl{0.93} 
		\def\stepDispl{0.03}
		\pgfmathsetmacro{\maxDispl}{\minDispl+\stepDispl*(\maxIdx-\minIdx)}
		\begin{groupplot}[
			group style={
				group name=myplot,
				group size=5 by 2,
				horizontal sep=5pt,
				vertical sep=5pt,
				xticklabels at=edge bottom,
				xlabels at=edge bottom,
				yticklabels at=edge left,
				ylabels at=edge left
			},
			width=0.4\linewidth,
			axis equal image,
			xlabel={$x$},
			x unit={mm},
			ylabel={$y$},
			y unit={mm},
			xmin = 0.5,xmax = 0.7,
			ymin = 0.2,ymax = 0.55,
			cycle list name=MyCycleList,
			every mark/.style={solid},
			legend to name=zelda,
			legend columns = 5,
			legend style={/tikz/every even column/.append style={column sep=0.25cm}}
			]

			\pgfplotsforeachungrouped \step in {\minIdx,...,\maxIdx}%
			{
				\pgfmathsetmacro{\disp}{\minDispl+\stepDispl*(\step-\minIdx)}
				\def\AT{1}
				\def\ORDER{4}
				\def\MESH{0.50}
				\def\MESHl{\MESH l0}
				\edef\tmp
				{
					\noexpand\nextgroupplot
					\noexpand\addlegendimage{black!50,opacity=1.0,mark=*,only marks,mark size=1pt,solid};
					\noexpand\addlegendentry{TP};
					\noexpand\addlegendimage{col1,mark=*,only marks,mark size=1pt,solid};
					\noexpand\addlegendentry{Implicit};
					\noexpand\addlegendimage{col2,mark=*,only marks,mark size=1pt,solid};
					\noexpand\addlegendentry{Explicit};
					\noexpand\addlegendimage{col3,mark=*,only marks,mark size=1pt,solid};
					\noexpand\addlegendentry{Hybrid};

					\noexpand\addplot+[black!50,opacity=0.5,mark=*,only marks,mark size=0.15pt,solid] table[header=true,col sep=comma, x index=1, y index=2,] {Data/TB/shear_TB_0.50l0_AT-\AT_Order\ORDER_contours/contour_\step.csv};

					\noexpand\addplot+[col1,mark=*,only marks,mark size=0.3pt] table[header=true,col sep=comma, x index=1, y index=2,] {Data/implicit/shear_THB_4l0_\MESHl_AT-\AT_Order\ORDER_contours/contour_\step.csv};
					\noexpand\addplot+[col2,mark=*,only marks,mark size=0.3pt] table[header=true,col sep=comma, x index=1, y index=2,] {Data/explicit/shear_THB_4l0_\MESHl_AT-\AT_Order\ORDER_contours/contour_\step.csv};
					\noexpand\addplot+[col3,mark=*,only marks,mark size=0.3pt] table[header=true,col sep=comma, x index=1, y index=2,] {Data/hybrid/shear_THB_4l0_\MESHl_AT-\AT_Order\ORDER_contours/contour_\step.csv};
					\noexpand\node[anchor=north east,scale=0.8] at (rel axis cs: 1,1) {$u=\pgfmathprintnumber[precision=3]{\disp}\cdot 10^{-3}$};
				}
				\tmp
			}
		\end{groupplot}
		\node[anchor=south] (legend) at ($(myplot c1r1.north east)!0.5!(myplot c5r1.north west)$)
		{\pgfplotslegendfromname{zelda}};
	\end{tikzpicture}
	\caption{Shear test.}
	\label{fig:implicitexplicit_shear}
\end{subfigure}
\caption{Contour plots of the damage field $d$ at $d=0.5$ at different horizontal or vertical displacements $u$ for the tensile (\subref{fig:implicitexplicit_tensile}) and shear (\subref{fig:implicitexplicit_shear}) tests, comparing the tensor-product B-spline basis with the adaptive implicit, explicit and hybrid approaches using the fourth-order AT1 model on a mesh with coarsest mesh size $h=\frac{\ell_0}{2}$.}
\label{fig:implicitexplicit}
\end{figure}

\begin{figure}
\centering
\begin{subfigure}{\linewidth}
\centering

\caption{Tensile test.}
\label{fig:implicitexplicit_plots_tensile}
\end{subfigure}

\begin{subfigure}{\linewidth}
\centering
%
\caption{Shear test.}
\label{fig:implicitexplicit_plots_shear}
\end{subfigure}
\caption{Dissipated Energy $\mathcal{D}$ (top row), total CPU-time (middle row) and number of degrees of freedom (\# DoFs, bottom row) comparison for the (\subref{fig:implicitexplicit_plots_tensile}) tensile and (\subref{fig:implicitexplicit_plots_shear}) shear tests using the second order AT1 (left), fourth order AT1 (middle left), second order AT2 (middle right) and fourth order AT2 (right) models, all plotted against the vertical displacement of the top boundary $u_y$. All results are provided for the fine mesh ($h=\ell_0/4$) except for the fourth order AT-1 model, as it provides the same accuracy on the coarse mesh ($h=\ell_0/2$). {\NEW The markers indicate the results obtained on the tensor-product B-spline basis, while the lines indicate the results obtained on the adaptive THB-spline basis with hybrid load stepping.}}
\label{fig:implicitexplicit_plots}
\end{figure}

\subsection{CPU time breakdown}
\label{subsec:CPU}

\Cref{fig:benchmarks_tensile_CPU,fig:benchmarks_shear_CPU} provide plots of the computational costs in terms of CPU wall time for the tensile and shear simulations, respectively, comparing results obtained on a tensor-product B-spline basis and on an adaptive THB-spline basis (with hybrid load stepping). In \cref{fig:implicitexplicit_cpu}, we compare the computational costs of the the explicit, implicit and hybrid approaches. The CPU times are subdivided into assembly and solver times for the elasticity and phase-field and optionally (only for adaptive meshes) a projection part. Based on these figures, several conclusions can be drawn.\\

Firstly, {\NEW some general observations, independent of the phase-field model, mesh adaptivity or the load-stepping strategy can be made based on }\cref{fig:benchmarks_tensile_CPU,fig:benchmarks_shear_CPU}. {\NEW That is, the figures} show that computational costs are dominated by the elasticity problem. This is explained by the fact that the number of degrees of freedom of the elasticity problem is higher due to the vector-valued solution space. In addition, within the staggered iterations the elasticity problem is solved using multiple iterations, whereas the phase-field problem is restricted to one. Moreover, the phase-field system of equations can partially be assembled once since it does not depend on solutions $\VEC{u}$ and $d$, hence only \cref{eq:SLCP_free_energy_matrix_vector} is assembled in each staggered iteration. Lastly, the phase-field solution typically converges faster within the staggered iteration, hence the PSOR solver is initialized almost perfectly, lowering the solver costs of the phase-field problem.\\

{\NEW Secondly, \cref{fig:benchmarks_tensile_CPU,fig:benchmarks_shear_CPU} also show that, for tensor-product B-spline meshes, the second-order AT1 model is generally the most expensive model, followed by the fourth-order AT1 model, the fourth-order AT2 model and the second-order AT2 model. On the adaptive meshes, on the other hand, no clear ranking can be observed. Moreover, comparing fine to coarse meshes, the CPU time reduction is typically a factor 4 for tensor-product meshes, which is in line with the expected reduction of the number of degrees of freedom. For adaptive meshes, a factor 3 (fourth-order models) up to a factor 5 (second-order models) is observed. Provided that the fourth-order AT1 model provides similar accuracy on the coarse mesh compared to the fine mesh, a CPU time reduction of a factor 4 for tensor-product meshes and a factor 3 for adaptive meshes is achieved by using coarser meshes. This is in line with the observations made in \cite{greco_at1}. Moreover, choosing the adaptive meshing algorithm of the present paper on a coarse mesh over uniform tensor B-spline bases on a fine mesh for the fourth-order AT1 model results in a CPU time reduction of a factor 10 for the tensile test and a factor 13 for the shear test.}\\

{\NEW Thirdly, explict load stepping is always faster than implicit and hybrid load stepping, as shown in \cref{fig:implicitexplicit_cpu}. However, as observed in \cref{fig:implicitexplicit}, explicit load stepping leads to inaccurate results in cases where mesh refinement lacks behind the crack evolution. Overall, comparing the computational costs of hybrid or implicit load shows that the differences are marginal, with hybrid load stepping being slightly more efficient. This is explained by the fact that most of the computational costs are occuring during crack propagation (see \cref{fig:implicitexplicit_plots}), where the implicit and hybrid load stepping strategies are identical.}

\begin{figure}
	\centering
	\begin{subfigure}{\linewidth}
		\centering
		\begin{tikzpicture}
			\begin{groupplot}[
				group style={
					group name=myplot,
					group size=4 by 5,
					vertical sep=0pt,
					horizontal sep=10pt,
					xticklabels at=all,
					xlabels at=all
				},
				width=0.3\linewidth,
				height=0.1\textheight,
				legend columns = 4,
				xbar stacked,
				table/col sep=comma,
				table/header=true,
				xbar stacked,
				enlarge x limits={0.01},
				enlarge y limits={0},
				ytick={0},
				yticklabels={},
				grid=none,
				axis x line*=middle,
				scaled x ticks=false,
				xtick align=center,
				major tick length=15pt,
				minor tick length=12pt,
				/pgf/bar width=8pt,
				every x tick scale label/.style={
					at={(xticklabel cs:0.9,0)},
					anchor=near xticklabel
				},
				]

                \pgfplotsforeachungrouped \AT / \ORDER in {1/2,1/4,2/2,2/4}
                {
                    \def\MESH{0.25}
                    \def\MESHl{\MESH l0}
                    \def\fileTP{Data/sum_data/TB/tensile_TB_\MESHl_AT-\AT_Order\ORDER_sum.csv}

                    \pgfplotstableread[col sep=comma]{\fileTP}\tempTableTP

                    \pgfplotstablegetelem{0}{pfSolverTime}\of\tempTableTP \pgfmathsetmacro{\tpA}{\pgfplotsretval}
                    \pgfplotstablegetelem{0}{pfAssemblyTime}\of\tempTableTP \pgfmathsetmacro{\tpB}{\pgfplotsretval}
                    \pgfplotstablegetelem{0}{elSolverTime}\of\tempTableTP \pgfmathsetmacro{\tpC}{\pgfplotsretval}
                    \pgfplotstablegetelem{0}{elAssemblyTime}\of\tempTableTP \pgfmathsetmacro{\tpD}{\pgfplotsretval}
                    \pgfmathsetmacro{\totalTP}{(\tpA + \tpB + \tpC + \tpD)/1e3}

                    \expandafter\xdef\csname totalFine\AT\ORDER\endcsname{\totalTP}

                    \edef\doPlotTP
                    {
                        \noexpand\nextgroupplot[
                            title={AT\AT, Order \ORDER},
                            xmin=0, xmax=7e3, xticklabels={},
                            every node near coord/.append style={anchor=west, font=\noexpand\tiny, xshift=2pt,fill=white,text=black,}
                        ]
                        \noexpand\addplot+[PFSolver] table[y expr={0}, x=pfSolverTime] {\fileTP};
                        \noexpand\addplot+[PFAssembler] table[y expr={0}, x=pfAssemblyTime] {\fileTP};
                        \noexpand\addplot+[ElSolver] table[y expr={0}, x=elSolverTime] {\fileTP};
                        \noexpand\addplot+[ElAssembler, 
                            nodes near coords={\noexpand\pgfmathprintnumber[fixed, precision=2]{\totalTP}}
                        ] table[y expr={0}, x=elAssemblyTime] {\fileTP};
                    }
                    \doPlotTP
                }

                \pgfplotsforeachungrouped \AT / \ORDER in {1/2,1/4,2/2,2/4}
                {
                    \def\MESH{0.25}
                    \def\MESHl{\MESH l0}
                    \def\fileTP{Data/sum_data/TB/tensile_TB_\MESHl_AT-\AT_Order\ORDER_sum.csv}
                    \def\fileTHB{Data/sum_data/implicit/tensile_THB_4l0_\MESHl_AT-\AT_Order\ORDER_sum.csv}

                    \pgfplotstableread[col sep=comma]{\fileTP}\tempTableTP
                    \pgfplotstableread[col sep=comma]{\fileTHB}\tempTableTHB

                    \pgfplotstablegetelem{0}{pfSolverTime}\of\tempTableTP \pgfmathsetmacro{\tpA}{\pgfplotsretval}
                    \pgfplotstablegetelem{0}{pfAssemblyTime}\of\tempTableTP \pgfmathsetmacro{\tpB}{\pgfplotsretval}
                    \pgfplotstablegetelem{0}{elSolverTime}\of\tempTableTP \pgfmathsetmacro{\tpC}{\pgfplotsretval}
                    \pgfplotstablegetelem{0}{elAssemblyTime}\of\tempTableTP \pgfmathsetmacro{\tpD}{\pgfplotsretval}
                    \pgfmathsetmacro{\totalTP}{(\tpA + \tpB + \tpC + \tpD)/1e3}

                    \pgfplotstablegetelem{0}{pfSolverTime}\of\tempTableTHB \pgfmathsetmacro{\thbA}{\pgfplotsretval}
                    \pgfplotstablegetelem{0}{pfAssemblyTime}\of\tempTableTHB \pgfmathsetmacro{\thbB}{\pgfplotsretval}
                    \pgfplotstablegetelem{0}{elSolverTime}\of\tempTableTHB \pgfmathsetmacro{\thbC}{\pgfplotsretval}
                    \pgfplotstablegetelem{0}{elAssemblyTime}\of\tempTableTHB \pgfmathsetmacro{\thbD}{\pgfplotsretval}
                    \pgfmathsetmacro{\totalTHB}{(\thbA + \thbB + \thbC + \thbD)/1e3}

                    \pgfmathsetmacro{\perc}{\totalTHB / \totalTP * 100}

                    \edef\doPlotTHB
                    {
                        \noexpand\nextgroupplot[
                            xmin=0, xmax=7e3, 
                            scaled x ticks=true,
                            scaled x ticks=base 10:-3,
						    xtick scale label code/.code={},
                            every node near coord/.append style={anchor=west, font=\noexpand\tiny, xshift=2pt,fill=white,text=black,}
                        ]
                        \noexpand\addplot+[PFSolver] table[y expr={0}, x=pfSolverTime] {\fileTHB};
                        \noexpand\addplot+[PFAssembler] table[y expr={0}, x=pfAssemblyTime] {\fileTHB};
                        \noexpand\addplot+[ElSolver] table[y expr={0}, x=elSolverTime] {\fileTHB};
                        \noexpand\addplot+[ElAssembler, 
                            nodes near coords={\noexpand\pgfmathprintnumber[fixed, precision=2]{\totalTHB}}
                        ] table[y expr={0}, x=elAssemblyTime] {\fileTHB};
                    }
                    \doPlotTHB
                }

				\pgfplotsforeachungrouped \AT / \ORDER in {1/2,1/4,2/2,2/4}
				{
					\edef\tmp
					{
						\noexpand\nextgroupplot[group/empty plot]
					}
					\tmp
				}



                \pgfplotsforeachungrouped \AT / \ORDER in {1/2,1/4,2/2,2/4}
                {
                    \def\MESH{0.50}
                    \def\MESHl{\MESH l0}
                    \def\fileTP{Data/sum_data/TB/tensile_TB_\MESHl_AT-\AT_Order\ORDER_sum.csv}

                    \pgfplotstableread[col sep=comma]{\fileTP}\tempTableTP

                    \pgfplotstablegetelem{0}{pfSolverTime}\of\tempTableTP \pgfmathsetmacro{\tpA}{\pgfplotsretval}
                    \pgfplotstablegetelem{0}{pfAssemblyTime}\of\tempTableTP \pgfmathsetmacro{\tpB}{\pgfplotsretval}
                    \pgfplotstablegetelem{0}{elSolverTime}\of\tempTableTP \pgfmathsetmacro{\tpC}{\pgfplotsretval}
                    \pgfplotstablegetelem{0}{elAssemblyTime}\of\tempTableTP \pgfmathsetmacro{\tpD}{\pgfplotsretval}
                    \pgfmathsetmacro{\totalTP}{(\tpA + \tpB + \tpC + \tpD)/1e3}

                    \edef\doPlotTP
                    {
                        \noexpand\nextgroupplot[
                            xmin=0, xmax=2e3, xticklabels={},
                            every node near coord/.append style={anchor=west, font=\noexpand\tiny, xshift=2pt,fill=white,text=black,}
                        ]
                        \noexpand\addplot+[PFSolver] table[y expr={0}, x=pfSolverTime] {\fileTP};
                        \noexpand\addplot+[PFAssembler] table[y expr={0}, x=pfAssemblyTime] {\fileTP};
                        \noexpand\addplot+[ElSolver] table[y expr={0}, x=elSolverTime] {\fileTP};
                        \noexpand\addplot+[ElAssembler, 
                            nodes near coords={\noexpand\pgfmathprintnumber[fixed, precision=2]{\totalTP}}
                        ] table[y expr={0}, x=elAssemblyTime] {\fileTP};
                    }
                    \doPlotTP
                }

                \pgfplotsforeachungrouped \AT / \ORDER in {1/2,1/4,2/2,2/4}
                {
                    \def\MESH{0.50}
                    \def\MESHl{\MESH l0}
                    \def\fileTP{Data/sum_data/TB/tensile_TB_\MESHl_AT-\AT_Order\ORDER_sum.csv}
                    \def\fileTHB{Data/sum_data/implicit/tensile_THB_4l0_\MESHl_AT-\AT_Order\ORDER_sum.csv}

                    \pgfplotstableread[col sep=comma]{\fileTP}\tempTableTP
                    \pgfplotstableread[col sep=comma]{\fileTHB}\tempTableTHB

                    \pgfplotstablegetelem{0}{pfSolverTime}\of\tempTableTP \pgfmathsetmacro{\tpA}{\pgfplotsretval}
                    \pgfplotstablegetelem{0}{pfAssemblyTime}\of\tempTableTP \pgfmathsetmacro{\tpB}{\pgfplotsretval}
                    \pgfplotstablegetelem{0}{elSolverTime}\of\tempTableTP \pgfmathsetmacro{\tpC}{\pgfplotsretval}
                    \pgfplotstablegetelem{0}{elAssemblyTime}\of\tempTableTP \pgfmathsetmacro{\tpD}{\pgfplotsretval}
                    \pgfmathsetmacro{\totalTP}{(\tpA + \tpB + \tpC + \tpD)/1e3}

                    \pgfplotstablegetelem{0}{pfSolverTime}\of\tempTableTHB \pgfmathsetmacro{\thbA}{\pgfplotsretval}
                    \pgfplotstablegetelem{0}{pfAssemblyTime}\of\tempTableTHB \pgfmathsetmacro{\thbB}{\pgfplotsretval}
                    \pgfplotstablegetelem{0}{elSolverTime}\of\tempTableTHB \pgfmathsetmacro{\thbC}{\pgfplotsretval}
                    \pgfplotstablegetelem{0}{elAssemblyTime}\of\tempTableTHB \pgfmathsetmacro{\thbD}{\pgfplotsretval}
                    \pgfmathsetmacro{\totalTHB}{(\thbA + \thbB + \thbC + \thbD)/1e3}

                    \pgfmathsetmacro{\perc}{\totalTHB / \totalTP * 100}

                    \edef\doPlotTHB
                    {
                        \noexpand\nextgroupplot[
                            xmin = 0,xmax = 2e3,
                            xlabel={CPU time},x unit={s},
                            scaled x ticks=true,
                            scaled x ticks=base 10:-3,
                            legend to name={CommonLegend},
                            every node near coord/.append style={anchor=west, font=\noexpand\tiny, xshift=2pt,fill=white,text=black,}
                        ]
                        \noexpand\addplot+[PFSolver] table[y expr={0}, x=pfSolverTime] {\fileTHB};
						\noexpand\addlegendentry{Phase-field Solver}
                        \noexpand\addplot+[PFAssembler] table[y expr={0}, x=pfAssemblyTime] {\fileTHB};
						\noexpand\addlegendentry{Phase-field Assembly}
                        \noexpand\addplot+[ElSolver] table[y expr={0}, x=elSolverTime] {\fileTHB};
						\noexpand\addlegendentry{Elasticity Solver}
                        \noexpand\addplot+[ElAssembler, 
                            nodes near coords={\noexpand\pgfmathprintnumber[fixed, precision=2]{\totalTHB}}
                        ] table[y expr={0}, x=elAssemblyTime] {\fileTHB};
						\noexpand\addlegendentry{Elasticity Assembly}
                    }
                    \doPlotTHB
                }
			\end{groupplot}
			\node[anchor=east,text width=0.7cm,align=right,font=\small] (TPfine) at (myplot c1r1.west) {TP};
			\node[anchor=east,text width=0.7cm,align=right,font=\small] (THBfine) at (myplot c1r2.west) {THB};
			\draw[decorate,decoration={brace,amplitude=5pt,mirror,raise=0}] (TPfine.north west) -- (THBfine.south west) node[midway,above,rotate=90,inner sep=10]{Fine};

			\node[anchor=east,text width=0.7cm,align=right,font=\small] (TPcoarse) at (myplot c1r4.west) {TP};
			\node[anchor=east,text width=0.7cm,align=right,font=\small] (THBcoarse) at (myplot c1r5.west) {THB};
			\draw[decorate,decoration={brace,amplitude=5pt,mirror,raise=0}] (TPcoarse.north west) -- (THBcoarse.south west) node[midway,above,rotate=90,inner sep=10]{Coarse};
			\path (myplot c1r5) -- node[anchor=north,below=1.1cm]{\pgfplotslegendfromname{CommonLegend}} (myplot c4r5);
		\end{tikzpicture}
		\caption{Tensile test}
		\label{fig:benchmarks_tensile_CPU}
	\end{subfigure}

	\begin{subfigure}{\linewidth}
		\centering
		\begin{tikzpicture}
			\begin{groupplot}[
				group style={
					group name=myplot,
					group size=4 by 5,
					vertical sep=0pt,
					horizontal sep=10pt,
					xticklabels at=all,
					xlabels at=all
				},
				width=0.3\linewidth,
				height=0.1\textheight,
				legend columns = 4,
				xbar stacked,
				table/col sep=comma,
				table/header=true,
				xbar stacked,
				enlarge x limits={0.01},
				enlarge y limits={0},
				ytick={0},
				yticklabels={},
				grid=none,
				axis x line*=middle,
				scaled x ticks=false,
				xtick align=center,
				major tick length=15pt,
				minor tick length=12pt,
				/pgf/bar width=8pt,
				every x tick scale label/.style={
					at={(xticklabel cs:0.9,0)},
					anchor=near xticklabel
				},
				]

                \pgfplotsforeachungrouped \AT / \ORDER in {1/2,1/4,2/2,2/4}
                {
                    \def\MESH{0.25}
                    \def\MESHl{\MESH l0}
                    \def\fileTP{Data/sum_data/TB/shear_TB_\MESHl_AT-\AT_Order\ORDER_sum.csv}

                    \pgfplotstableread[col sep=comma]{\fileTP}\tempTableTP

                    \pgfplotstablegetelem{0}{pfSolverTime}\of\tempTableTP \pgfmathsetmacro{\tpA}{\pgfplotsretval}
                    \pgfplotstablegetelem{0}{pfAssemblyTime}\of\tempTableTP \pgfmathsetmacro{\tpB}{\pgfplotsretval}
                    \pgfplotstablegetelem{0}{elSolverTime}\of\tempTableTP \pgfmathsetmacro{\tpC}{\pgfplotsretval}
                    \pgfplotstablegetelem{0}{elAssemblyTime}\of\tempTableTP \pgfmathsetmacro{\tpD}{\pgfplotsretval}
                    \pgfmathsetmacro{\totalTP}{(\tpA + \tpB + \tpC + \tpD)/1e3}

                    \expandafter\xdef\csname totalFine\AT\ORDER\endcsname{\totalTP}

                    \edef\doPlotTP
                    {
                        \noexpand\nextgroupplot[
                            title={AT\AT, Order \ORDER},
                            xmin=0, xmax=15e3, xticklabels={},
                            every node near coord/.append style={anchor=west, font=\noexpand\tiny, xshift=2pt,fill=white,text=black,}
                        ]
                        \noexpand\addplot+[PFSolver] table[y expr={0}, x=pfSolverTime] {\fileTP};
                        \noexpand\addplot+[PFAssembler] table[y expr={0}, x=pfAssemblyTime] {\fileTP};
                        \noexpand\addplot+[ElSolver] table[y expr={0}, x=elSolverTime] {\fileTP};
                        \noexpand\addplot+[ElAssembler, 
                            nodes near coords={\noexpand\pgfmathprintnumber[fixed, precision=2]{\totalTP}}
                        ] table[y expr={0}, x=elAssemblyTime] {\fileTP};
                    }
                    \doPlotTP
                }

                \pgfplotsforeachungrouped \AT / \ORDER in {1/2,1/4,2/2,2/4}
                {
                    \def\MESH{0.25}
                    \def\MESHl{\MESH l0}
                    \def\fileTP{Data/sum_data/TB/shear_TB_\MESHl_AT-\AT_Order\ORDER_sum.csv}
                    \def\fileTHB{Data/sum_data/implicit/shear_THB_4l0_\MESHl_AT-\AT_Order\ORDER_sum.csv}

                    \pgfplotstableread[col sep=comma]{\fileTP}\tempTableTP
                    \pgfplotstableread[col sep=comma]{\fileTHB}\tempTableTHB

                    \pgfplotstablegetelem{0}{pfSolverTime}\of\tempTableTP \pgfmathsetmacro{\tpA}{\pgfplotsretval}
                    \pgfplotstablegetelem{0}{pfAssemblyTime}\of\tempTableTP \pgfmathsetmacro{\tpB}{\pgfplotsretval}
                    \pgfplotstablegetelem{0}{elSolverTime}\of\tempTableTP \pgfmathsetmacro{\tpC}{\pgfplotsretval}
                    \pgfplotstablegetelem{0}{elAssemblyTime}\of\tempTableTP \pgfmathsetmacro{\tpD}{\pgfplotsretval}
                    \pgfmathsetmacro{\totalTP}{(\tpA + \tpB + \tpC + \tpD)/1e3}

                    \pgfplotstablegetelem{0}{pfSolverTime}\of\tempTableTHB \pgfmathsetmacro{\thbA}{\pgfplotsretval}
                    \pgfplotstablegetelem{0}{pfAssemblyTime}\of\tempTableTHB \pgfmathsetmacro{\thbB}{\pgfplotsretval}
                    \pgfplotstablegetelem{0}{elSolverTime}\of\tempTableTHB \pgfmathsetmacro{\thbC}{\pgfplotsretval}
                    \pgfplotstablegetelem{0}{elAssemblyTime}\of\tempTableTHB \pgfmathsetmacro{\thbD}{\pgfplotsretval}
                    \pgfmathsetmacro{\totalTHB}{(\thbA + \thbB + \thbC + \thbD)/1e3}

                    \pgfmathsetmacro{\perc}{\totalTHB / \totalTP * 100}

                    \edef\doPlotTHB
                    {
                        \noexpand\nextgroupplot[
                            xmin=0, xmax=15e3, 
                            scaled x ticks=true,
                            scaled x ticks=base 10:-3,
						    xtick scale label code/.code={},
                            every node near coord/.append style={anchor=west, font=\noexpand\tiny, xshift=2pt,fill=white,text=black,}
                        ]
                        \noexpand\addplot+[PFSolver] table[y expr={0}, x=pfSolverTime] {\fileTHB};
                        \noexpand\addplot+[PFAssembler] table[y expr={0}, x=pfAssemblyTime] {\fileTHB};
                        \noexpand\addplot+[ElSolver] table[y expr={0}, x=elSolverTime] {\fileTHB};
                        \noexpand\addplot+[ElAssembler, 
                            nodes near coords={\noexpand\pgfmathprintnumber[fixed, precision=2]{\totalTHB}}
                        ] table[y expr={0}, x=elAssemblyTime] {\fileTHB};
                    }
                    \doPlotTHB
                }

				\pgfplotsforeachungrouped \AT / \ORDER in {1/2,1/4,2/2,2/4}
				{
					\edef\tmp
					{
						\noexpand\nextgroupplot[group/empty plot]
					}
					\tmp
				}



                \pgfplotsforeachungrouped \AT / \ORDER in {1/2,1/4,2/2,2/4}
                {
                    \def\MESH{0.50}
                    \def\MESHl{\MESH l0}
                    \def\fileTP{Data/sum_data/TB/shear_TB_\MESHl_AT-\AT_Order\ORDER_sum.csv}

                    \pgfplotstableread[col sep=comma]{\fileTP}\tempTableTP

                    \pgfplotstablegetelem{0}{pfSolverTime}\of\tempTableTP \pgfmathsetmacro{\tpA}{\pgfplotsretval}
                    \pgfplotstablegetelem{0}{pfAssemblyTime}\of\tempTableTP \pgfmathsetmacro{\tpB}{\pgfplotsretval}
                    \pgfplotstablegetelem{0}{elSolverTime}\of\tempTableTP \pgfmathsetmacro{\tpC}{\pgfplotsretval}
                    \pgfplotstablegetelem{0}{elAssemblyTime}\of\tempTableTP \pgfmathsetmacro{\tpD}{\pgfplotsretval}
                    \pgfmathsetmacro{\totalTP}{(\tpA + \tpB + \tpC + \tpD)/1e3}

                    \edef\doPlotTP
                    {
                        \noexpand\nextgroupplot[
                            xmin=0, xmax=3e3, xticklabels={},
                            every node near coord/.append style={anchor=west, font=\noexpand\tiny, xshift=2pt,fill=white,text=black,}
                        ]
                        \noexpand\addplot+[PFSolver] table[y expr={0}, x=pfSolverTime] {\fileTP};
                        \noexpand\addplot+[PFAssembler] table[y expr={0}, x=pfAssemblyTime] {\fileTP};
                        \noexpand\addplot+[ElSolver] table[y expr={0}, x=elSolverTime] {\fileTP};
                        \noexpand\addplot+[ElAssembler, 
                            nodes near coords={\noexpand\pgfmathprintnumber[fixed, precision=2]{\totalTP}}
                        ] table[y expr={0}, x=elAssemblyTime] {\fileTP};
                    }
                    \doPlotTP
                }

                \pgfplotsforeachungrouped \AT / \ORDER in {1/2,1/4,2/2,2/4}
                {
                    \def\MESH{0.50}
                    \def\MESHl{\MESH l0}
                    \def\fileTP{Data/sum_data/TB/shear_TB_\MESHl_AT-\AT_Order\ORDER_sum.csv}
                    \def\fileTHB{Data/sum_data/implicit/shear_THB_4l0_\MESHl_AT-\AT_Order\ORDER_sum.csv}

                    \pgfplotstableread[col sep=comma]{\fileTP}\tempTableTP
                    \pgfplotstableread[col sep=comma]{\fileTHB}\tempTableTHB

                    \pgfplotstablegetelem{0}{pfSolverTime}\of\tempTableTP \pgfmathsetmacro{\tpA}{\pgfplotsretval}
                    \pgfplotstablegetelem{0}{pfAssemblyTime}\of\tempTableTP \pgfmathsetmacro{\tpB}{\pgfplotsretval}
                    \pgfplotstablegetelem{0}{elSolverTime}\of\tempTableTP \pgfmathsetmacro{\tpC}{\pgfplotsretval}
                    \pgfplotstablegetelem{0}{elAssemblyTime}\of\tempTableTP \pgfmathsetmacro{\tpD}{\pgfplotsretval}
                    \pgfmathsetmacro{\totalTP}{(\tpA + \tpB + \tpC + \tpD)/1e3}

                    \pgfplotstablegetelem{0}{pfSolverTime}\of\tempTableTHB \pgfmathsetmacro{\thbA}{\pgfplotsretval}
                    \pgfplotstablegetelem{0}{pfAssemblyTime}\of\tempTableTHB \pgfmathsetmacro{\thbB}{\pgfplotsretval}
                    \pgfplotstablegetelem{0}{elSolverTime}\of\tempTableTHB \pgfmathsetmacro{\thbC}{\pgfplotsretval}
                    \pgfplotstablegetelem{0}{elAssemblyTime}\of\tempTableTHB \pgfmathsetmacro{\thbD}{\pgfplotsretval}
                    \pgfmathsetmacro{\totalTHB}{(\thbA + \thbB + \thbC + \thbD)/1e3}

                    \pgfmathsetmacro{\perc}{\totalTHB / \totalTP * 100}

                    \edef\doPlotTHB
                    {
                        \noexpand\nextgroupplot[
                            xmin = 0,xmax = 3e3,
                            xlabel={CPU time},x unit={s},
                            scaled x ticks=true,
                            scaled x ticks=base 10:-3,
                            legend to name={CommonLegend},
                            every node near coord/.append style={anchor=west, font=\noexpand\tiny, xshift=2pt,fill=white,text=black,}
                        ]
                        \noexpand\addplot+[PFSolver] table[y expr={0}, x=pfSolverTime] {\fileTHB};
						\noexpand\addlegendentry{Phase-field Solver}
                        \noexpand\addplot+[PFAssembler] table[y expr={0}, x=pfAssemblyTime] {\fileTHB};
						\noexpand\addlegendentry{Phase-field Assembly}
                        \noexpand\addplot+[ElSolver] table[y expr={0}, x=elSolverTime] {\fileTHB};
						\noexpand\addlegendentry{Elasticity Solver}
                        \noexpand\addplot+[ElAssembler, 
                            nodes near coords={\noexpand\pgfmathprintnumber[fixed, precision=2]{\totalTHB}}
                        ] table[y expr={0}, x=elAssemblyTime] {\fileTHB};
						\noexpand\addlegendentry{Elasticity Assembly}
                    }
                    \doPlotTHB
                }
			\end{groupplot}
			\node[anchor=east,text width=0.7cm,align=right,font=\small] (TPfine) at (myplot c1r1.west) {TP};
			\node[anchor=east,text width=0.7cm,align=right,font=\small] (THBfine) at (myplot c1r2.west) {THB};
			\draw[decorate,decoration={brace,amplitude=5pt,mirror,raise=0}] (TPfine.north west) -- (THBfine.south west) node[midway,above,rotate=90,inner sep=10]{Fine};

			\node[anchor=east,text width=0.7cm,align=right,font=\small] (TPcoarse) at (myplot c1r4.west) {TP};
			\node[anchor=east,text width=0.7cm,align=right,font=\small] (THBcoarse) at (myplot c1r5.west) {THB};
			\draw[decorate,decoration={brace,amplitude=5pt,mirror,raise=0}] (TPcoarse.north west) -- (THBcoarse.south west) node[midway,above,rotate=90,inner sep=10]{Coarse};
			\path (myplot c1r5) -- node[anchor=north,below=1.1cm]{\pgfplotslegendfromname{CommonLegend}} (myplot c4r5);
		\end{tikzpicture}
		\caption{Shear test}
		\label{fig:benchmarks_shear_CPU}
	\end{subfigure}
	\caption{Computational cost comparison for different phase-field models using adaptive meshing for the (\subref{fig:benchmarks_tensile_CPU}) tensile and (\subref{fig:benchmarks_shear_CPU}) shear tests for the tensor-product B-spline basis (TP) and the adaptive THB-spline basis with hybrid load stepping (THB). The bars are stacked and colored based on the solver and assembly time of the phase-field and elasticity problems. {\NEW The projection, marking and refinement time for the adaptive THB-spline basis are omitted as they are negligible. The total CPU time for each case is shown as a number on top of the bars.}}
\end{figure}

\begin{figure}
\centering
\begin{subfigure}{\linewidth}
\centering
\begin{tikzpicture}
\begin{groupplot}[
        group style={
            group name=myplot,
            group size=4 by 7,
            vertical sep=0pt,
            horizontal sep=10pt,
            xticklabels at=all,
            xlabels at=all
        },
    width=0.3\linewidth,
    height=0.1\textheight,
    legend columns = 4,
    xbar stacked,
    table/col sep=comma,
    table/header=true,
    xbar stacked,
    enlarge x limits={0.01},
    enlarge y limits={0},
    ytick={0},
    yticklabels={},
    grid=none,
    axis x line*=middle,
    scaled x ticks=false,
    xtick align=center,
    major tick length=15pt,
    minor tick length=12pt,
    /pgf/bar width=8pt,
    every x tick scale label/.style={
            at={(xticklabel cs:0.9,0)},
            anchor=near xticklabel
        },
]

\pgfplotsforeachungrouped \AT / \ORDER in {1/2,1/4,2/2,2/4}
{
    \def\MESH{0.25}
    \def\MESHl{\MESH l0}
	\def\fileTHB{Data/sum_data/implicit/tensile_THB_4l0_\MESHl_AT-\AT_Order\ORDER_sum.csv}

	\pgfplotstableread[col sep=comma]{\fileTHB}\tempTableTHB

	\pgfplotstablegetelem{0}{pfSolverTime}\of\tempTableTHB \pgfmathsetmacro{\THBA}{\pgfplotsretval}
	\pgfplotstablegetelem{0}{pfAssemblyTime}\of\tempTableTHB \pgfmathsetmacro{\THBB}{\pgfplotsretval}
	\pgfplotstablegetelem{0}{elSolverTime}\of\tempTableTHB \pgfmathsetmacro{\THBC}{\pgfplotsretval}
	\pgfplotstablegetelem{0}{elAssemblyTime}\of\tempTableTHB \pgfmathsetmacro{\THBD}{\pgfplotsretval}
	\pgfplotstablegetelem{0}{projectionTime}\of\tempTableTHB \pgfmathsetmacro{\THBE}{\pgfplotsretval}
	\pgfplotstablegetelem{0}{markingTime}\of\tempTableTHB \pgfmathsetmacro{\THBF}{\pgfplotsretval}
	\pgfplotstablegetelem{0}{refinementTime}\of\tempTableTHB \pgfmathsetmacro{\THBG}{\pgfplotsretval}
	\pgfmathsetmacro{\totalTHB}{(\THBA + \THBB + \THBC + \THBD + \THBE + \THBF + \THBG)/1e3}

	\expandafter\xdef\csname totalFine\AT\ORDER\endcsname{\totalTHB}

    \edef\tmp
    {
        \noexpand\nextgroupplot[
			title={AT\AT, Order \ORDER},
			xmin = 0,xmax = 2e3, xticklabels={},
			every node near coord/.append style={anchor=west, font=\noexpand\tiny, xshift=2pt,fill=white,text=black,}
		]
		\noexpand\addplot+[PFSolver] table[y expr={0}, x=pfSolverTime] {\fileTHB};
		\noexpand\addplot+[PFAssembler] table[y expr={0}, x=pfAssemblyTime] {\fileTHB};
		\noexpand\addplot+[ElSolver] table[y expr={0}, x=elSolverTime] {\fileTHB};
		\noexpand\addplot+[ElAssembler] table[y expr={0}, x=elAssemblyTime] {\fileTHB};
		\noexpand\addplot+[Projection] table[y expr={0}, x=projectionTime] {\fileTHB};           
		\noexpand\addplot+[Marking] table[y expr={0}, x=markingTime] {\fileTHB};
		\noexpand\addplot+[Refinement, 
			nodes near coords={\noexpand\pgfmathprintnumber[fixed, precision=2]{\totalTHB}}
		] table[y expr={0}, x=refinementTime] {\fileTHB};
    }
    \tmp
}

\pgfplotsforeachungrouped \AT / \ORDER in {1/2,1/4,2/2,2/4}
{
    \def\MESH{0.25}
    \def\MESHl{\MESH l0}
	\def\fileTHB{Data/sum_data/hybrid/tensile_THB_4l0_\MESHl_AT-\AT_Order\ORDER_sum.csv}

	\pgfplotstableread[col sep=comma]{\fileTHB}\tempTableTHB

	\pgfplotstablegetelem{0}{pfSolverTime}\of\tempTableTHB \pgfmathsetmacro{\THBA}{\pgfplotsretval}
	\pgfplotstablegetelem{0}{pfAssemblyTime}\of\tempTableTHB \pgfmathsetmacro{\THBB}{\pgfplotsretval}
	\pgfplotstablegetelem{0}{elSolverTime}\of\tempTableTHB \pgfmathsetmacro{\THBC}{\pgfplotsretval}
	\pgfplotstablegetelem{0}{elAssemblyTime}\of\tempTableTHB \pgfmathsetmacro{\THBD}{\pgfplotsretval}
	\pgfplotstablegetelem{0}{projectionTime}\of\tempTableTHB \pgfmathsetmacro{\THBE}{\pgfplotsretval}
	\pgfplotstablegetelem{0}{markingTime}\of\tempTableTHB \pgfmathsetmacro{\THBF}{\pgfplotsretval}
	\pgfplotstablegetelem{0}{refinementTime}\of\tempTableTHB \pgfmathsetmacro{\THBG}{\pgfplotsretval}
	\pgfmathsetmacro{\totalTHB}{(\THBA + \THBB + \THBC + \THBD + \THBE + \THBF + \THBG)/1e3}

	\expandafter\xdef\csname totalFine\AT\ORDER\endcsname{\totalTHB}

    \edef\tmp
    {
        \noexpand\nextgroupplot[
			xmin = 0,xmax = 2e3, xticklabels={},
			every node near coord/.append style={anchor=west, font=\noexpand\tiny, xshift=2pt,fill=white,text=black,}
		]
		\noexpand\addplot+[PFSolver] table[y expr={0}, x=pfSolverTime] {\fileTHB};
		\noexpand\addplot+[PFAssembler] table[y expr={0}, x=pfAssemblyTime] {\fileTHB};
		\noexpand\addplot+[ElSolver] table[y expr={0}, x=elSolverTime] {\fileTHB};
		\noexpand\addplot+[ElAssembler] table[y expr={0}, x=elAssemblyTime] {\fileTHB};
		\noexpand\addplot+[Projection] table[y expr={0}, x=projectionTime] {\fileTHB};           
		\noexpand\addplot+[Marking] table[y expr={0}, x=markingTime] {\fileTHB};
		\noexpand\addplot+[Refinement, 
			nodes near coords={\noexpand\pgfmathprintnumber[fixed, precision=2]{\totalTHB}}
		] table[y expr={0}, x=refinementTime] {\fileTHB};
    }
    \tmp
}

\pgfplotsforeachungrouped \AT / \ORDER in {1/2,1/4,2/2,2/4}
{
    \def\MESH{0.25}
    \def\MESHl{\MESH l0}
	\def\fileTHB{Data/sum_data/explicit/tensile_THB_4l0_\MESHl_AT-\AT_Order\ORDER_sum.csv}

	\pgfplotstableread[col sep=comma]{\fileTHB}\tempTableTHB

	\pgfplotstablegetelem{0}{pfSolverTime}\of\tempTableTHB \pgfmathsetmacro{\THBA}{\pgfplotsretval}
	\pgfplotstablegetelem{0}{pfAssemblyTime}\of\tempTableTHB \pgfmathsetmacro{\THBB}{\pgfplotsretval}
	\pgfplotstablegetelem{0}{elSolverTime}\of\tempTableTHB \pgfmathsetmacro{\THBC}{\pgfplotsretval}
	\pgfplotstablegetelem{0}{elAssemblyTime}\of\tempTableTHB \pgfmathsetmacro{\THBD}{\pgfplotsretval}
	\pgfplotstablegetelem{0}{projectionTime}\of\tempTableTHB \pgfmathsetmacro{\THBE}{\pgfplotsretval}
	\pgfplotstablegetelem{0}{markingTime}\of\tempTableTHB \pgfmathsetmacro{\THBF}{\pgfplotsretval}
	\pgfplotstablegetelem{0}{refinementTime}\of\tempTableTHB \pgfmathsetmacro{\THBG}{\pgfplotsretval}
	\pgfmathsetmacro{\totalTHB}{(\THBA + \THBB + \THBC + \THBD + \THBE + \THBF + \THBG)/1e3}

	\expandafter\xdef\csname totalFine\AT\ORDER\endcsname{\totalTHB}

    \edef\tmp
    {
        \noexpand\nextgroupplot[
			xmin = 0,xmax = 2e3,
			scaled x ticks=true,
			scaled x ticks=base 10:-3,
			xtick scale label code/.code={},
			every node near coord/.append style={anchor=west, font=\noexpand\tiny, xshift=2pt,fill=white,text=black,}
		]
		\noexpand\addplot+[PFSolver] table[y expr={0}, x=pfSolverTime] {\fileTHB};
		\noexpand\addplot+[PFAssembler] table[y expr={0}, x=pfAssemblyTime] {\fileTHB};
		\noexpand\addplot+[ElSolver] table[y expr={0}, x=elSolverTime] {\fileTHB};
		\noexpand\addplot+[ElAssembler] table[y expr={0}, x=elAssemblyTime] {\fileTHB};
		\noexpand\addplot+[Projection] table[y expr={0}, x=projectionTime] {\fileTHB};           
		\noexpand\addplot+[Marking] table[y expr={0}, x=markingTime] {\fileTHB};
		\noexpand\addplot+[Refinement, 
			nodes near coords={\noexpand\pgfmathprintnumber[fixed, precision=2]{\totalTHB}}
		] table[y expr={0}, x=refinementTime] {\fileTHB};
    }
    \tmp
}


\pgfplotsforeachungrouped \AT / \ORDER in {1/2,1/4,2/2,2/4}
{
    \def\MESH{0.25}
    \def\MESHl{\MESH l0}
    \edef\tmp
    {
        \noexpand\nextgroupplot[group/empty plot]
    }
    \tmp
}

\pgfplotsforeachungrouped \AT / \ORDER in {1/2,1/4,2/2,2/4}
{
    \def\MESH{0.50}
    \def\MESHl{\MESH l0}
	\def\fileTHB{Data/sum_data/implicit/tensile_THB_4l0_\MESHl_AT-\AT_Order\ORDER_sum.csv}

	\pgfplotstableread[col sep=comma]{\fileTHB}\tempTableTHB

	\pgfplotstablegetelem{0}{pfSolverTime}\of\tempTableTHB \pgfmathsetmacro{\THBA}{\pgfplotsretval}
	\pgfplotstablegetelem{0}{pfAssemblyTime}\of\tempTableTHB \pgfmathsetmacro{\THBB}{\pgfplotsretval}
	\pgfplotstablegetelem{0}{elSolverTime}\of\tempTableTHB \pgfmathsetmacro{\THBC}{\pgfplotsretval}
	\pgfplotstablegetelem{0}{elAssemblyTime}\of\tempTableTHB \pgfmathsetmacro{\THBD}{\pgfplotsretval}
	\pgfplotstablegetelem{0}{projectionTime}\of\tempTableTHB \pgfmathsetmacro{\THBE}{\pgfplotsretval}
	\pgfplotstablegetelem{0}{markingTime}\of\tempTableTHB \pgfmathsetmacro{\THBF}{\pgfplotsretval}
	\pgfplotstablegetelem{0}{refinementTime}\of\tempTableTHB \pgfmathsetmacro{\THBG}{\pgfplotsretval}
	\pgfmathsetmacro{\totalTHB}{(\THBA + \THBB + \THBC + \THBD + \THBE + \THBF + \THBG)/1e3}

	\expandafter\xdef\csname totalFine\AT\ORDER\endcsname{\totalTHB}

    \edef\tmp
    {
        \noexpand\nextgroupplot[
			xmin = 0,xmax = 0.6e3, xticklabels={},
			every node near coord/.append style={anchor=west, font=\noexpand\tiny, xshift=2pt,fill=white,text=black,}
		]
		\noexpand\addplot+[PFSolver] table[y expr={0}, x=pfSolverTime] {\fileTHB};
		\noexpand\addplot+[PFAssembler] table[y expr={0}, x=pfAssemblyTime] {\fileTHB};
		\noexpand\addplot+[ElSolver] table[y expr={0}, x=elSolverTime] {\fileTHB};
		\noexpand\addplot+[ElAssembler] table[y expr={0}, x=elAssemblyTime] {\fileTHB};
		\noexpand\addplot+[Projection] table[y expr={0}, x=projectionTime] {\fileTHB};           
		\noexpand\addplot+[Marking] table[y expr={0}, x=markingTime] {\fileTHB};
		\noexpand\addplot+[Refinement, 
			nodes near coords={\noexpand\pgfmathprintnumber[fixed, precision=2]{\totalTHB}}
		] table[y expr={0}, x=refinementTime] {\fileTHB};
    }
    \tmp
}

\pgfplotsforeachungrouped \AT / \ORDER in {1/2,1/4,2/2,2/4}
{
    \def\MESH{0.50}
    \def\MESHl{\MESH l0}
	\def\fileTHB{Data/sum_data/hybrid/tensile_THB_4l0_\MESHl_AT-\AT_Order\ORDER_sum.csv}

	\pgfplotstableread[col sep=comma]{\fileTHB}\tempTableTHB

	\pgfplotstablegetelem{0}{pfSolverTime}\of\tempTableTHB \pgfmathsetmacro{\THBA}{\pgfplotsretval}
	\pgfplotstablegetelem{0}{pfAssemblyTime}\of\tempTableTHB \pgfmathsetmacro{\THBB}{\pgfplotsretval}
	\pgfplotstablegetelem{0}{elSolverTime}\of\tempTableTHB \pgfmathsetmacro{\THBC}{\pgfplotsretval}
	\pgfplotstablegetelem{0}{elAssemblyTime}\of\tempTableTHB \pgfmathsetmacro{\THBD}{\pgfplotsretval}
	\pgfplotstablegetelem{0}{projectionTime}\of\tempTableTHB \pgfmathsetmacro{\THBE}{\pgfplotsretval}
	\pgfplotstablegetelem{0}{markingTime}\of\tempTableTHB \pgfmathsetmacro{\THBF}{\pgfplotsretval}
	\pgfplotstablegetelem{0}{refinementTime}\of\tempTableTHB \pgfmathsetmacro{\THBG}{\pgfplotsretval}
	\pgfmathsetmacro{\totalTHB}{(\THBA + \THBB + \THBC + \THBD + \THBE + \THBF + \THBG)/1e3}

	\expandafter\xdef\csname totalFine\AT\ORDER\endcsname{\totalTHB}

    \edef\tmp
    {
        \noexpand\nextgroupplot[
			xmin = 0,xmax = 0.6e3, xticklabels={},
			every node near coord/.append style={anchor=west, font=\noexpand\tiny, xshift=2pt,fill=white,text=black,}
		]
		\noexpand\addplot+[PFSolver] table[y expr={0}, x=pfSolverTime] {\fileTHB};
		\noexpand\addplot+[PFAssembler] table[y expr={0}, x=pfAssemblyTime] {\fileTHB};
		\noexpand\addplot+[ElSolver] table[y expr={0}, x=elSolverTime] {\fileTHB};
		\noexpand\addplot+[ElAssembler] table[y expr={0}, x=elAssemblyTime] {\fileTHB};
		\noexpand\addplot+[Projection] table[y expr={0}, x=projectionTime] {\fileTHB};           
		\noexpand\addplot+[Marking] table[y expr={0}, x=markingTime] {\fileTHB};
		\noexpand\addplot+[Refinement, 
			nodes near coords={\noexpand\pgfmathprintnumber[fixed, precision=2]{\totalTHB}}
		] table[y expr={0}, x=refinementTime] {\fileTHB};
    }
    \tmp
}

\pgfplotsforeachungrouped \AT / \ORDER in {1/2,1/4,2/2,2/4}
{
    \def\MESH{0.50}
    \def\MESHl{\MESH l0}
	\def\fileTHB{Data/sum_data/explicit/tensile_THB_4l0_\MESHl_AT-\AT_Order\ORDER_sum.csv}

	\pgfplotstableread[col sep=comma]{\fileTHB}\tempTableTHB

	\pgfplotstablegetelem{0}{pfSolverTime}\of\tempTableTHB \pgfmathsetmacro{\THBA}{\pgfplotsretval}
	\pgfplotstablegetelem{0}{pfAssemblyTime}\of\tempTableTHB \pgfmathsetmacro{\THBB}{\pgfplotsretval}
	\pgfplotstablegetelem{0}{elSolverTime}\of\tempTableTHB \pgfmathsetmacro{\THBC}{\pgfplotsretval}
	\pgfplotstablegetelem{0}{elAssemblyTime}\of\tempTableTHB \pgfmathsetmacro{\THBD}{\pgfplotsretval}
	\pgfplotstablegetelem{0}{projectionTime}\of\tempTableTHB \pgfmathsetmacro{\THBE}{\pgfplotsretval}
	\pgfplotstablegetelem{0}{markingTime}\of\tempTableTHB \pgfmathsetmacro{\THBF}{\pgfplotsretval}
	\pgfplotstablegetelem{0}{refinementTime}\of\tempTableTHB \pgfmathsetmacro{\THBG}{\pgfplotsretval}
	\pgfmathsetmacro{\totalTHB}{(\THBA + \THBB + \THBC + \THBD + \THBE + \THBF + \THBG)/1e3}

	\expandafter\xdef\csname totalFine\AT\ORDER\endcsname{\totalTHB}

    \edef\tmp
    {
        \noexpand\nextgroupplot[
			xmin = 0,xmax = 0.6e3,
			xlabel={CPU time},x unit={s},
			scaled x ticks=true,
			scaled x ticks=base 10:-3,
			every node near coord/.append style={anchor=west, font=\noexpand\tiny, xshift=2pt,fill=white,text=black,},
			legend to name={CommonLegend}
		]
		\noexpand\addplot+[PFSolver] table[y expr={0}, x=pfSolverTime] {\fileTHB};
		\noexpand\addlegendentry{Phase-field Solver}
		\noexpand\addplot+[PFAssembler] table[y expr={0}, x=pfAssemblyTime] {\fileTHB};
		\noexpand\addlegendentry{Phase-field Assembly}
		\noexpand\addplot+[ElSolver] table[y expr={0}, x=elSolverTime] {\fileTHB};
		\noexpand\addlegendentry{Elasticity Solver}
		\noexpand\addplot+[ElAssembler] table[y expr={0}, x=elAssemblyTime] {\fileTHB};
		\noexpand\addlegendentry{Elasticity Assembly}
		\noexpand\addplot+[Projection] table[y expr={0}, x=projectionTime] {\fileTHB};           
		\noexpand\addlegendentry{Projection}
		\noexpand\addplot+[Marking] table[y expr={0}, x=markingTime] {\fileTHB};
		\noexpand\addlegendentry{Marking}
		\noexpand\addplot+[Refinement, 
			nodes near coords={\noexpand\pgfmathprintnumber[fixed, precision=2]{\totalTHB}}
		] table[y expr={0}, x=refinementTime] {\fileTHB};
		\noexpand\addlegendentry{Refinement}
    }
    \tmp
}
\end{groupplot}

\node[anchor=east,text width=1.2cm,align=right,font=\small] (ImplicitFine)at (myplot c1r1.west) {Implicit};
\node[anchor=east,text width=1.2cm,align=right,font=\small] (HybridFine)  at (myplot c1r2.west) {Hybrid};
\node[anchor=east,text width=1.2cm,align=right,font=\small] (ExplicitFine)at (myplot c1r3.west) {Explicit};
\draw[decorate,decoration={brace,amplitude=5pt,mirror,raise=0}] (ImplicitFine.north west) -- (ExplicitFine.south west) node[midway,above,rotate=90,inner sep=10]{Fine};

\node[anchor=east,text width=1.2cm,align=right,font=\small] (ImplicitCoarse)at (myplot c1r5.west) {Implicit};
\node[anchor=east,text width=1.2cm,align=right,font=\small] (HybridCoarse)  at (myplot c1r6.west) {Hybrid};
\node[anchor=east,text width=1.2cm,align=right,font=\small] (ExplicitCoarse)at (myplot c1r7.west) {Explicit};
\draw[decorate,decoration={brace,amplitude=5pt,mirror,raise=0}] (ImplicitCoarse.north west) -- (ExplicitCoarse.south west) node[midway,above,rotate=90,inner sep=10]{Coarse};

\path (myplot c1r7) -- node[anchor=north,below=1.1cm]{\pgfplotslegendfromname{CommonLegend}} (myplot c4r7);

\end{tikzpicture}
\caption{Tensile test.}
\label{fig:implicitexplicit_cpu_tensile}
\end{subfigure}

\begin{subfigure}{\linewidth}
\centering
\begin{tikzpicture}
\begin{groupplot}[
        group style={
            group name=myplot,
            group size=4 by 7,
            vertical sep=0pt,
            horizontal sep=10pt,
            xticklabels at=all,
            xlabels at=all
        },
    width=0.3\linewidth,
    height=0.1\textheight,
    legend columns = 4,
    xbar stacked,
    table/col sep=comma,
    table/header=true,
    xbar stacked,
    enlarge x limits={0.01},
    enlarge y limits={0},
    ytick={0},
    yticklabels={},
    grid=none,
    axis x line*=middle,
    scaled x ticks=false,
    xtick align=center,
    major tick length=15pt,
    minor tick length=12pt,
    /pgf/bar width=8pt,
    every x tick scale label/.style={
            at={(xticklabel cs:0.9,0)},
            anchor=near xticklabel
        },
]

\pgfplotsforeachungrouped \AT / \ORDER in {1/2,1/4,2/2,2/4}
{
    \def\MESH{0.25}
    \def\MESHl{\MESH l0}
	\def\fileTHB{Data/sum_data/implicit/shear_THB_4l0_\MESHl_AT-\AT_Order\ORDER_sum.csv}

	\pgfplotstableread[col sep=comma]{\fileTHB}\tempTableTHB

	\pgfplotstablegetelem{0}{pfSolverTime}\of\tempTableTHB \pgfmathsetmacro{\THBA}{\pgfplotsretval}
	\pgfplotstablegetelem{0}{pfAssemblyTime}\of\tempTableTHB \pgfmathsetmacro{\THBB}{\pgfplotsretval}
	\pgfplotstablegetelem{0}{elSolverTime}\of\tempTableTHB \pgfmathsetmacro{\THBC}{\pgfplotsretval}
	\pgfplotstablegetelem{0}{elAssemblyTime}\of\tempTableTHB \pgfmathsetmacro{\THBD}{\pgfplotsretval}
	\pgfplotstablegetelem{0}{projectionTime}\of\tempTableTHB \pgfmathsetmacro{\THBE}{\pgfplotsretval}
	\pgfplotstablegetelem{0}{markingTime}\of\tempTableTHB \pgfmathsetmacro{\THBF}{\pgfplotsretval}
	\pgfmathsetmacro{\totalTHB}{(\THBA + \THBB + \THBC + \THBD + \THBE + \THBF)/1e3}

	\expandafter\xdef\csname totalFine\AT\ORDER\endcsname{\totalTHB}

    \edef\tmp
    {
        \noexpand\nextgroupplot[
			title={AT\AT, Order \ORDER},
			xmin = 0,xmax = 3.5e3, xticklabels={},
			every node near coord/.append style={anchor=west, font=\noexpand\tiny, xshift=2pt,fill=white,text=black,}
		]
		\noexpand\addplot+[PFSolver] table[y expr={0}, x=pfSolverTime] {\fileTHB};
		\noexpand\addplot+[PFAssembler] table[y expr={0}, x=pfAssemblyTime] {\fileTHB};
		\noexpand\addplot+[ElSolver] table[y expr={0}, x=elSolverTime] {\fileTHB};
		\noexpand\addplot+[ElAssembler] table[y expr={0}, x=elAssemblyTime] {\fileTHB};
		\noexpand\addplot+[Projection] table[y expr={0}, x=projectionTime] {\fileTHB};           
		\noexpand\addplot+[Marking] table[y expr={0}, x=markingTime] {\fileTHB};
		\noexpand\addplot+[Refinement, 
			nodes near coords={\noexpand\pgfmathprintnumber[fixed, precision=2]{\totalTHB}}
		] table[y expr={0}, x=refinementTime] {\fileTHB};
    }
    \tmp
}

\pgfplotsforeachungrouped \AT / \ORDER in {1/2,1/4,2/2,2/4}
{
    \def\MESH{0.25}
    \def\MESHl{\MESH l0}
	\def\fileTHB{Data/sum_data/hybrid/shear_THB_4l0_\MESHl_AT-\AT_Order\ORDER_sum.csv}

	\pgfplotstableread[col sep=comma]{\fileTHB}\tempTableTHB

	\pgfplotstablegetelem{0}{pfSolverTime}\of\tempTableTHB \pgfmathsetmacro{\THBA}{\pgfplotsretval}
	\pgfplotstablegetelem{0}{pfAssemblyTime}\of\tempTableTHB \pgfmathsetmacro{\THBB}{\pgfplotsretval}
	\pgfplotstablegetelem{0}{elSolverTime}\of\tempTableTHB \pgfmathsetmacro{\THBC}{\pgfplotsretval}
	\pgfplotstablegetelem{0}{elAssemblyTime}\of\tempTableTHB \pgfmathsetmacro{\THBD}{\pgfplotsretval}
	\pgfplotstablegetelem{0}{projectionTime}\of\tempTableTHB \pgfmathsetmacro{\THBE}{\pgfplotsretval}
	\pgfplotstablegetelem{0}{markingTime}\of\tempTableTHB \pgfmathsetmacro{\THBF}{\pgfplotsretval}
	\pgfmathsetmacro{\totalTHB}{(\THBA + \THBB + \THBC + \THBD + \THBE + \THBF)/1e3}

	\expandafter\xdef\csname totalFine\AT\ORDER\endcsname{\totalTHB}

    \edef\tmp
    {
        \noexpand\nextgroupplot[
			xmin = 0,xmax = 3.5e3, xticklabels={},
			every node near coord/.append style={anchor=west, font=\noexpand\tiny, xshift=2pt,fill=white,text=black,}
		]
		\noexpand\addplot+[PFSolver] table[y expr={0}, x=pfSolverTime] {\fileTHB};
		\noexpand\addplot+[PFAssembler] table[y expr={0}, x=pfAssemblyTime] {\fileTHB};
		\noexpand\addplot+[ElSolver] table[y expr={0}, x=elSolverTime] {\fileTHB};
		\noexpand\addplot+[ElAssembler] table[y expr={0}, x=elAssemblyTime] {\fileTHB};
		\noexpand\addplot+[Projection] table[y expr={0}, x=projectionTime] {\fileTHB};           
		\noexpand\addplot+[Marking] table[y expr={0}, x=markingTime] {\fileTHB};
		\noexpand\addplot+[Refinement, 
			nodes near coords={\noexpand\pgfmathprintnumber[fixed, precision=2]{\totalTHB}}
		] table[y expr={0}, x=refinementTime] {\fileTHB};
    }
    \tmp
}

\pgfplotsforeachungrouped \AT / \ORDER in {1/2,1/4,2/2,2/4}
{
    \def\MESH{0.25}
    \def\MESHl{\MESH l0}
	\def\fileTHB{Data/sum_data/explicit/shear_THB_4l0_\MESHl_AT-\AT_Order\ORDER_sum.csv}

	\pgfplotstableread[col sep=comma]{\fileTHB}\tempTableTHB

	\pgfplotstablegetelem{0}{pfSolverTime}\of\tempTableTHB \pgfmathsetmacro{\THBA}{\pgfplotsretval}
	\pgfplotstablegetelem{0}{pfAssemblyTime}\of\tempTableTHB \pgfmathsetmacro{\THBB}{\pgfplotsretval}
	\pgfplotstablegetelem{0}{elSolverTime}\of\tempTableTHB \pgfmathsetmacro{\THBC}{\pgfplotsretval}
	\pgfplotstablegetelem{0}{elAssemblyTime}\of\tempTableTHB \pgfmathsetmacro{\THBD}{\pgfplotsretval}
	\pgfplotstablegetelem{0}{projectionTime}\of\tempTableTHB \pgfmathsetmacro{\THBE}{\pgfplotsretval}
	\pgfplotstablegetelem{0}{markingTime}\of\tempTableTHB \pgfmathsetmacro{\THBF}{\pgfplotsretval}
	\pgfmathsetmacro{\totalTHB}{(\THBA + \THBB + \THBC + \THBD + \THBE + \THBF)/1e3}

	\expandafter\xdef\csname totalFine\AT\ORDER\endcsname{\totalTHB}

    \edef\tmp
    {
        \noexpand\nextgroupplot[
			xmin = 0,xmax = 3.5e3,
			scaled x ticks=true,
			scaled x ticks=base 10:-3,
			xtick scale label code/.code={},
			every node near coord/.append style={anchor=west, font=\noexpand\tiny, xshift=2pt,fill=white,text=black,}
		]
		\noexpand\addplot+[PFSolver] table[y expr={0}, x=pfSolverTime] {\fileTHB};
		\noexpand\addplot+[PFAssembler] table[y expr={0}, x=pfAssemblyTime] {\fileTHB};
		\noexpand\addplot+[ElSolver] table[y expr={0}, x=elSolverTime] {\fileTHB};
		\noexpand\addplot+[ElAssembler] table[y expr={0}, x=elAssemblyTime] {\fileTHB};
		\noexpand\addplot+[Projection] table[y expr={0}, x=projectionTime] {\fileTHB};           
		\noexpand\addplot+[Marking] table[y expr={0}, x=markingTime] {\fileTHB};
		\noexpand\addplot+[Refinement, 
			nodes near coords={\noexpand\pgfmathprintnumber[fixed, precision=2]{\totalTHB}}
		] table[y expr={0}, x=refinementTime] {\fileTHB};
    }
    \tmp
}


\pgfplotsforeachungrouped \AT / \ORDER in {1/2,1/4,2/2,2/4}
{
    \def\MESH{0.25}
    \def\MESHl{\MESH l0}
    \edef\tmp
    {
        \noexpand\nextgroupplot[group/empty plot]
    }
    \tmp
}

\pgfplotsforeachungrouped \AT / \ORDER in {1/2,1/4,2/2,2/4}
{
    \def\MESH{0.50}
    \def\MESHl{\MESH l0}
	\def\fileTHB{Data/sum_data/implicit/shear_THB_4l0_\MESHl_AT-\AT_Order\ORDER_sum.csv}

	\pgfplotstableread[col sep=comma]{\fileTHB}\tempTableTHB

	\pgfplotstablegetelem{0}{pfSolverTime}\of\tempTableTHB \pgfmathsetmacro{\THBA}{\pgfplotsretval}
	\pgfplotstablegetelem{0}{pfAssemblyTime}\of\tempTableTHB \pgfmathsetmacro{\THBB}{\pgfplotsretval}
	\pgfplotstablegetelem{0}{elSolverTime}\of\tempTableTHB \pgfmathsetmacro{\THBC}{\pgfplotsretval}
	\pgfplotstablegetelem{0}{elAssemblyTime}\of\tempTableTHB \pgfmathsetmacro{\THBD}{\pgfplotsretval}
	\pgfplotstablegetelem{0}{projectionTime}\of\tempTableTHB \pgfmathsetmacro{\THBE}{\pgfplotsretval}
	\pgfplotstablegetelem{0}{markingTime}\of\tempTableTHB \pgfmathsetmacro{\THBF}{\pgfplotsretval}
	\pgfmathsetmacro{\totalTHB}{(\THBA + \THBB + \THBC + \THBD + \THBE + \THBF)/1e3}

	\expandafter\xdef\csname totalFine\AT\ORDER\endcsname{\totalTHB}

    \edef\tmp
    {
        \noexpand\nextgroupplot[
			xmin = 0,xmax = 0.8e3, xticklabels={},
			every node near coord/.append style={anchor=west, font=\noexpand\tiny, xshift=2pt,fill=white,text=black,}
		]
		\noexpand\addplot+[PFSolver] table[y expr={0}, x=pfSolverTime] {\fileTHB};
		\noexpand\addplot+[PFAssembler] table[y expr={0}, x=pfAssemblyTime] {\fileTHB};
		\noexpand\addplot+[ElSolver] table[y expr={0}, x=elSolverTime] {\fileTHB};
		\noexpand\addplot+[ElAssembler] table[y expr={0}, x=elAssemblyTime] {\fileTHB};
		\noexpand\addplot+[Projection] table[y expr={0}, x=projectionTime] {\fileTHB};           
		\noexpand\addplot+[Marking] table[y expr={0}, x=markingTime] {\fileTHB};
		\noexpand\addplot+[Refinement, 
			nodes near coords={\noexpand\pgfmathprintnumber[fixed, precision=2]{\totalTHB}}
		] table[y expr={0}, x=refinementTime] {\fileTHB};
    }
    \tmp
}

\pgfplotsforeachungrouped \AT / \ORDER in {1/2,1/4,2/2,2/4}
{
    \def\MESH{0.50}
    \def\MESHl{\MESH l0}
	\def\fileTHB{Data/sum_data/hybrid/shear_THB_4l0_\MESHl_AT-\AT_Order\ORDER_sum.csv}

	\pgfplotstableread[col sep=comma]{\fileTHB}\tempTableTHB

	\pgfplotstablegetelem{0}{pfSolverTime}\of\tempTableTHB \pgfmathsetmacro{\THBA}{\pgfplotsretval}
	\pgfplotstablegetelem{0}{pfAssemblyTime}\of\tempTableTHB \pgfmathsetmacro{\THBB}{\pgfplotsretval}
	\pgfplotstablegetelem{0}{elSolverTime}\of\tempTableTHB \pgfmathsetmacro{\THBC}{\pgfplotsretval}
	\pgfplotstablegetelem{0}{elAssemblyTime}\of\tempTableTHB \pgfmathsetmacro{\THBD}{\pgfplotsretval}
	\pgfplotstablegetelem{0}{projectionTime}\of\tempTableTHB \pgfmathsetmacro{\THBE}{\pgfplotsretval}
	\pgfplotstablegetelem{0}{markingTime}\of\tempTableTHB \pgfmathsetmacro{\THBF}{\pgfplotsretval}
	\pgfmathsetmacro{\totalTHB}{(\THBA + \THBB + \THBC + \THBD + \THBE + \THBF)/1e3}

	\expandafter\xdef\csname totalFine\AT\ORDER\endcsname{\totalTHB}

    \edef\tmp
    {
        \noexpand\nextgroupplot[
			xmin = 0,xmax = 0.8e3, xticklabels={},
			every node near coord/.append style={anchor=west, font=\noexpand\tiny, xshift=2pt,fill=white,text=black,}
		]
		\noexpand\addplot+[PFSolver] table[y expr={0}, x=pfSolverTime] {\fileTHB};
		\noexpand\addplot+[PFAssembler] table[y expr={0}, x=pfAssemblyTime] {\fileTHB};
		\noexpand\addplot+[ElSolver] table[y expr={0}, x=elSolverTime] {\fileTHB};
		\noexpand\addplot+[ElAssembler] table[y expr={0}, x=elAssemblyTime] {\fileTHB};
		\noexpand\addplot+[Projection] table[y expr={0}, x=projectionTime] {\fileTHB};           
		\noexpand\addplot+[Marking] table[y expr={0}, x=markingTime] {\fileTHB};
		\noexpand\addplot+[Refinement, 
			nodes near coords={\noexpand\pgfmathprintnumber[fixed, precision=2]{\totalTHB}}
		] table[y expr={0}, x=refinementTime] {\fileTHB};
    }
    \tmp
}

\pgfplotsforeachungrouped \AT / \ORDER in {1/2,1/4,2/2,2/4}
{
    \def\MESH{0.50}
    \def\MESHl{\MESH l0}
	\def\fileTHB{Data/sum_data/explicit/shear_THB_4l0_\MESHl_AT-\AT_Order\ORDER_sum.csv}

	\pgfplotstableread[col sep=comma]{\fileTHB}\tempTableTHB

	\pgfplotstablegetelem{0}{pfSolverTime}\of\tempTableTHB \pgfmathsetmacro{\THBA}{\pgfplotsretval}
	\pgfplotstablegetelem{0}{pfAssemblyTime}\of\tempTableTHB \pgfmathsetmacro{\THBB}{\pgfplotsretval}
	\pgfplotstablegetelem{0}{elSolverTime}\of\tempTableTHB \pgfmathsetmacro{\THBC}{\pgfplotsretval}
	\pgfplotstablegetelem{0}{elAssemblyTime}\of\tempTableTHB \pgfmathsetmacro{\THBD}{\pgfplotsretval}
	\pgfplotstablegetelem{0}{projectionTime}\of\tempTableTHB \pgfmathsetmacro{\THBE}{\pgfplotsretval}
	\pgfplotstablegetelem{0}{markingTime}\of\tempTableTHB \pgfmathsetmacro{\THBF}{\pgfplotsretval}
	\pgfmathsetmacro{\totalTHB}{(\THBA + \THBB + \THBC + \THBD + \THBE + \THBF)/1e3}

	\expandafter\xdef\csname totalFine\AT\ORDER\endcsname{\totalTHB}

    \edef\tmp
    {
        \noexpand\nextgroupplot[
			xmin = 0,xmax = 0.8e3,
			xlabel={CPU time},x unit={s},
			scaled x ticks=true,
			scaled x ticks=base 10:-3,
			every node near coord/.append style={anchor=west, font=\noexpand\tiny, xshift=2pt,fill=white,text=black,},
			legend to name={CommonLegend}
		]
		\noexpand\addplot+[PFSolver] table[y expr={0}, x=pfSolverTime] {\fileTHB};
		\noexpand\addlegendentry{Phase-field Solver}
		\noexpand\addplot+[PFAssembler] table[y expr={0}, x=pfAssemblyTime] {\fileTHB};
		\noexpand\addlegendentry{Phase-field Assembly}
		\noexpand\addplot+[ElSolver] table[y expr={0}, x=elSolverTime] {\fileTHB};
		\noexpand\addlegendentry{Elasticity Solver}
		\noexpand\addplot+[ElAssembler] table[y expr={0}, x=elAssemblyTime] {\fileTHB};
		\noexpand\addlegendentry{Elasticity Assembly}
		\noexpand\addplot+[Projection] table[y expr={0}, x=projectionTime] {\fileTHB};           
		\noexpand\addlegendentry{Projection}
		\noexpand\addplot+[Marking] table[y expr={0}, x=markingTime] {\fileTHB};
		\noexpand\addlegendentry{Marking}
		\noexpand\addplot+[Refinement, 
			nodes near coords={\noexpand\pgfmathprintnumber[fixed, precision=2]{\totalTHB}}
		] table[y expr={0}, x=refinementTime] {\fileTHB};
		\noexpand\addlegendentry{Refinement}
    }
    \tmp
}

\end{groupplot}

\node[anchor=east,text width=1.2cm,align=right,font=\small] (ImplicitFine)at (myplot c1r1.west) {Implicit};
\node[anchor=east,text width=1.2cm,align=right,font=\small] (HybridFine)  at (myplot c1r2.west) {Hybrid};
\node[anchor=east,text width=1.2cm,align=right,font=\small] (ExplicitFine)at (myplot c1r3.west) {Explicit};
\draw[decorate,decoration={brace,amplitude=5pt,mirror,raise=0}] (ImplicitFine.north west) -- (ExplicitFine.south west) node[midway,above,rotate=90,inner sep=10]{Fine};

\node[anchor=east,text width=1.2cm,align=right,font=\small] (ImplicitCoarse)at (myplot c1r5.west) {Implicit};
\node[anchor=east,text width=1.2cm,align=right,font=\small] (HybridCoarse)  at (myplot c1r6.west) {Hybrid};
\node[anchor=east,text width=1.2cm,align=right,font=\small] (ExplicitCoarse)at (myplot c1r7.west) {Explicit};
\draw[decorate,decoration={brace,amplitude=5pt,mirror,raise=0}] (ImplicitCoarse.north west) -- (ExplicitCoarse.south west) node[midway,above,rotate=90,inner sep=10]{Coarse};

\path (myplot c1r7) -- node[anchor=north,below=1.1cm]{\pgfplotslegendfromname{CommonLegend}} (myplot c4r7);

\end{tikzpicture}
\caption{Shear test.}
\label{fig:implicitexplicit_cpu_shear}
\end{subfigure}
\caption{Comparison of load stepping approaches for adaptive THB-spline phase-field models in terms of computational cost for the (\subref{fig:implicitexplicit_cpu_tensile}) tensile and (\subref{fig:implicitexplicit_cpu_shear}) shear tests. The bars are stacked and colored based on the solver and assembly time of the phase-field and elasticity problems. {\NEW The projection, marking, and refinement times are also included, despite their contributions are indistinguishable as these contributions are negligible. The total time for each configuration is shown as a label next to the refinement bar.}}
\label{fig:implicitexplicit_cpu}
\end{figure}


\clearpage

\section{Conclusions}
\label{sec:conclusions}
This paper presents an adaptive isogeometric method for phase-field brittle fracture simulations. The framework is based on the isogeometric analysis of higher-order phase-field simulation of the damage field based on combinations of the AT1 and AT2 formulations of second- and fourth-order phase-field models. These equations are solved on an adaptive basis composed of Truncated Hierarchical B-splines, which is refined admissibly to avoid spurious oscillations. Since brittle fracture is phenomenologically irreversible, the adaptive refinement scheme is solely based on refinement, using phase-field-based element marking. In addition, since brittle fracture is phenomenologically sudden, load-steps are optionally re-computed on a mesh obtained in the previous refinement iteration, aiming to capture significant propagations of the damage field accurately across load steps. The re-computation algorithms presented in this paper consist of so-called explicit, implicit or hyrid schemes. Lastly, we reason on the phenomenon of cross-talk in relation to the presented refinement strategy. This phenomenon is naturally eliminated by admissible meshes refined to sufficient depth.\\

The proposed adaptive isogeometric method is thoroughly tested on two common benchmarks for numerical phase-field brittle fracture: the SEN shear and tensile tests,. Based on the results obtained for the benchmark tests under different model configurations, several observations can be made, and conclusions can be drawn.\\ 
Firstly, the results show that for tensor-product meshes, the dissipated energy obtained using the fourth order AT1 model on a coarse mesh with size $h=\frac{\ell_0}{2}$ is close to the results obtained on a fine mesh with size $h=\frac{\ell_0}{4}$, while the other phase field models considered in this study show significant differences between the fine and coarse meshes. Consequently, it can be concluded that a coarse mesh is sufficient for the fourth-order AT1 model, providing a speed-up of approximately a factor 4, similar to the one observed in earlier work \cite{greco_at1}.\\
Secondly, the suddenness of fracture can only be captured accurately by using hybrid of implicit load stepping schemes. For the explicit load stepping sheme, the results show that crack propagation is limited by the propagation of the mesh. On the other hand, the implicit load stepping scheme provides a conservative counterpart: load step recomputation is performed as soon as one or more elements are marked for refinement. Consequently, load step recomputation is also performed in the pre-propagation phase, yielding high computational costs in this regime, for low accuracy gains. Finally, the hybrid load stepping scheme, provides a combination of the two methods: it aims to be explicit in the pre-propagation phase, and implicit in the propagation phase. Although this approach requires a threshold parameter to be chosen, it provides efficiency gains over the implicit scheme. For the fourth-order AT1 model, {\NEW the over-all speed-up of using the proposed adaptive meshing strategy is approximately a factor 10-13 compared to the tensor-product computation for the same model on a fine mesh}.\\

Based on the findings of this paper and on other recent advances in the simulation of brittle fracture via phase-field methods, we suggest the following directions for future work. Firstly, the present paper provides a thorough analysis of phase-field brittle fracture simulations in 2D. As volumetric fracture simulations imply a significant increase in degrees of freedom, solver costs are limiting reasonable simulation times. Therefore, we plan to investigate distributed memory parallelization methods on hierarchical meshes to reduce assembly and solver times. Secondly, the analyses in the present paper are restricted to quadratic THB-splines, because of the fact that the IPF phase-field initialization method is primarily investigated for that degree, and because quadratic basis functions provide minimal support for sufficient smoothness. Nevertheless, the methods and algorithms for refinement proposed in this paper are degree-independent, hence an investigation into the effect of higher-degree bases is on our future planning.

\section{Acknowledgments}
H.M. Verhelst acknowledges financial support the knowledge and innovation covenant Fundamental Fluid Dynamics Challenges in Inkjet Printing (Phase II) (KICH2.V4C.20.001), a joint research program of Canon Production Printing, Eindhoven University of Technology, Utrecht University, University of Twente, and the Netherlands Organization for Scientific Research (NWO).\\
A. Reali and L. Greco acknowledge the support of the Italian Ministry of University and Research (MUR) through the PRIN project COSMIC (No. 2022A79M75), funded by the European Union — Next Generation EU.\\
In addition, A. Reali is a member of the Gruppo Nazionale Calcolo Scientifico-Istituto Nazionale di Alta Matematica (GNCS-INDAM), and acknowledges the contribution of the National Recovery and Resilience Plan, Mission 4 Component 2 -- Investment 1.4 -- NATIONAL CENTER FOR HPC, BIG DATA AND QUANTUM COMPUTING, spoke 6.

\begin{appendix}
\section{Algorithmic overview}
\label{app:algorithms}
\label{appendix_algorithms}
In this section, we provide an overview of the algorithms used to solve the coupled phase-field fracture problem with adaptive refinement. \\

In \cref{alg:solve_elasticity,alg:solve_phasefield} we provide subroutines for solving the displacement and phase-field subproblems at each staggered iteration $\i\ge1$. The \textbf{solve\_elasticity} algorithm implements the Newton-Raphson scheme for solving \cref{eq:matrix_problem_u}, while the \textbf{solve\_phasefield} algorithm employs the PSOR method to solve the linear complementarity problem defined by equations (\ref{eq:matrix_problem_d1}) and (\ref{eq:matrix_problem_d2}). In both algorithms, we denote the loading step index by $n\ge1$, the staggered iteration index by $\i\ge1$, and the Newton--Raphson/PSOR iteration index by $\j\ge0$. In \cref{alg:solve_elasticity}, the subroutines \textbf{assemble\_elasticity\_LHS} and \textbf{assemble\_elasticity\_RHS} are responsible for assembling the stiffness matrix and the external force vector, respectively. In \cref{alg:solve_phasefield}, the subroutine \textbf{assemble\_phasefield\_LHS} constructs the solution-dependent part of the phase-field matrix, while \textbf{solve\_PSOR} implements the PSOR algorithm to solve the linear complementarity problem. The detailed implementation of the PSOR algorithm can be found in \cite{MARENGO2021114137}, while the assembly algorithms are based on standard procedures in isogeometric analysis.\\

\begin{algorithm}
\caption{\texttt{solve\_elasticity($\VEC{u}_n^{i-1}$,$\SCALAR{d}_n^{i-1}$,$\VEC{F}^{\text{ext}}_{n}$,$\mathcal{B}$,$\mathcal{Q}$,$\mathtt{TOL}_{\text{Pic},\VEC{u}}$)}\\
    \textbf{Description:} Newton-Raphson algorithm for the displacement field at a staggered iteration $\i\ge1$, solving \cref{eq:matrix_problem_u}.}
\label{alg:solve_elasticity}
\KwInput{Displacement field $\VEC{u}_n^{i}$ at current staggered iteration $\i$; phase-field from previous staggered iteration $\SCALAR{d}_n^{i-1}$; the basis $\mathcal{B}$, the mesh $\mathcal{Q}$ and the tolerance for the Picard iterations $\mathtt{TOL}_{\text{Pic},\VEC{u}}$.}
\For{$j=0,\dots$}
{
    \Comment{Assemble stiffness matrix}
    $\MAT{K}\gets\textbf{assemble\_elasticity\_LHS}(\VEC{u}_n^{i-1},\SCALAR{d}_n^{i-1})$\\
    \Comment{Solve \cref{eq:matrix_problem_u}}
    $\VEC{u}^{i,j}_n\gets \MAT{K}^{-1}\VEC{F}^{\text{ext}}_{n} $\\	
    \Comment{Test the residual}
    \If{$\text{Res}^{j}_{\text{Pic},\VEC{u}} < \mathtt{TOL}_{\text{Pic},\VEC{u}}$}
    {
        $\VEC{u}_n^{i}\gets \VEC{u}_n^{i,j}$\\
        \textbf{break}
    }
}
\KwOutput{Displacement field $\VEC{u}_n^{i}$ at current staggered iteration $\i$.}
\end{algorithm}

\begin{algorithm}
\caption{\texttt{solve\_phasefield}($\VEC{u}_n^{i}$,$\SCALAR{d}_n^{i-1}$,$\bm{\Phi}^{\text{II,IV}}$,$\bm{\phi}^{\text{II,IV}}$,$\mathcal{B}$,$\mathcal{Q}$,$\mathtt{TOL}_{\text{PSOR},\Delta\SCALAR{d}}$)\\
    \textbf{Description:} Solving the phase-field at staggered iteration $\i\ge1$, solving \cref{eq:matrix_problem_d1,eq:matrix_problem_d2}.}
\label{alg:solve_phasefield}
\KwInput{Displacement field $\VEC{u}_n^{i}$ at current staggered iteration $\i$; phase-field from previous staggered iteration $\SCALAR{d}_n^{i-1}$; solution-independent phase-field matrix $\bm{\Phi}^{\text{II,IV}}$ and vector $\bm{\phi}^{\text{II,IV}}$ from \cref{eq:dissipation_matrix,eq:dissipation_matrix_AT1,eq:dissipation_vector_AT1}; the basis $\mathcal{B}$, the mesh $\mathcal{Q}$ and the tolerance for the PSOR iterations $\mathtt{TOL}_{\text{PSOR},\Delta\SCALAR{d}}$.}
\Comment{Assemble solution-dependent phase-field matrix and vector from \cref{eq:SLCP_free_energy_matrix_vector}}
$\bm{\Psi}\gets\textbf{assemble\_phasefield\_LHS}(\VEC{u}_n^{i})$\\
$\bm{\psi}\gets\textbf{assemble\_phasefield\_RHS}(\VEC{u}_n^{i})$\\
\Comment{Assemble total phase-field matrix and vector from \cref{eq:Q_q}}
$\MAT{Q}^{II,IV} \gets \bm{\Psi} + G_c \,\bm{\Phi}^{\text{II,IV}}$\\
$\VEC{q}^{II,IV} \gets \MAT{Q}^{II,IV} \SCALAR{d}_n^{i-1} - \bm{\psi} - \bm{\phi}^{\text{II,IV}}$\\
\Comment{Solve equations (\ref{eq:matrix_problem_d1}) and (\ref{eq:matrix_problem_d2})}
$\Delta\SCALAR{d}_n^{i}\gets \textbf{solve\_PSOR}\qty(\MAT{Q}^{II,IV}, \VEC{q}^{II,IV},\mathtt{TOL}_{\text{PSOR},\Delta\SCALAR{d}})$\\
\Comment{Update phase-field solution}
$\SCALAR{d}_n^i \gets \SCALAR{d}_n^{i-1} + \Delta\SCALAR{d}_n^{i}$\\
\KwOutput{Phase-field solution $\SCALAR{d}_n^{i}$ at current staggered iteration $\i$.}
\end{algorithm}

\Cref{alg:load_step} summarizes the overall staggered solution scheme for a loading step $n\ge1$, utilizing the previously defined subroutines to solve the coupled problem iteratively until convergence is achieved. This load step algorithm can be embedded within a load stepping procedure to simulate the entire loading history a fixed basis $\mathcal{B}$ and mesh $\mathcal{Q}$, material parameters and boundary conditions. The subroutines \textbf{assemble\_elasticity\_RHS}, \textbf{assemble\_phasefield\_LHS} and \textbf{assemble\_phasefield\_RHS} are called once per load step, as they are independent of the solution fields. \\

\begin{algorithm}
\caption{\textbf{load\_step}($\VEC{u}^k$,$\SCALAR{d}^k$,$\mathcal{B}$,$\mathcal{Q}$,$\mathtt{TOL}_{\text{Pic},\VEC{u}}$,$\mathtt{TOL}_{\text{PSOR},\Delta\SCALAR{d}}$,$\mathtt{TOL}_{\text{stag}}$)\\
\textbf{Description:} Performs staggered solution on a fixed mesh and spline space for a load step $k$.}
\label{alg:load_step}
\KwInput{Previous load step solution $(\VEC{u}^k, \SCALAR{d}^k)$, the basis $\mathcal{B}$, the mesh $\mathcal{Q}$, the tolerance for the Picard iterations $\mathtt{TOL}_{\text{Pic},\VEC{u}}$, the tolerance for the PSOR iterations $\mathtt{TOL}_{\text{PSOR},\Delta\SCALAR{d}}$ and the tolerance for the staggered iterations $\mathtt{TOL}_{\text{stag}}$.}

\Comment{Initialize load step (assumes that the loads/displacements of load step $k+1$ are set)}
$\VEC{u}^{k+1} = \VEC{u}^k$\\
$\SCALAR{d}^{k+1} = \SCALAR{d}^k$\\

\Comment{Assemble solution-independent quantities}
$\VEC{F}^{\text{ext}}_{k+1}\gets\textbf{assemble\_elasticity\_RHS}(k+1)$\\
$\bm{\Phi}^{\text{II,IV}}\gets\textbf{assemble\_phasefield\_LHS}()$\\
$\bm{\phi}^{\text{II,IV}}\gets\textbf{assemble\_phasefield\_RHS}()$\\

\Comment{Staggered iterations}

\While{$\Vert R(\VEC{u}^{k+1},\SCALAR{d}^{k+1})\Vert>\mathtt{TOL}_{\text{stag}}$}
{
    \Comment{Solve elasticity problem}
    $\VEC{u}^{k+1} = \textbf{solve\_elasticity}(\VEC{u}^{k+1},\SCALAR{d}^{k+1},\VEC{F}^{\text{ext}}_{k+1},\mathcal{B},\mathcal{Q},\mathtt{TOL}_{\text{Pic},\VEC{u}})$\\

    \Comment{Solve phase-field problem}
    $\SCALAR{d}^{k+1} = \textbf{solve\_phasefield}(\VEC{u}^{k+1},\SCALAR{d}^{k},\SCALAR{d}^{k+1},\bm{\Phi}^{\text{II,IV}},\bm{\phi}^{\text{II,IV}},\mathcal{B},\mathcal{Q},\mathtt{TOL}_{\text{PSOR},\Delta\SCALAR{d}})$\\

    \Comment{Test the staggered residual}
    \If{$\text{Res}_{\text{stag}} < \mathtt{TOL}_{\text{stag}}$}{\textbf{break}}
}

\KwOutput{New displacement field $\VEC{u}^{k+1}$, new damage field $\SCALAR{d}^{k+1}$.}
\end{algorithm}

Contrary to \cref{alg:load_step}, we aim to define an algorithm which includes mesh adaptivity. To this end, we define \textbf{refine\_elements} which marks elements for refinement based on the damage field. The algorithms mark elements admissibly in a level-by-level fashion. In the last level, the support extension of the marked elements is also marked to eliminate cross-talk, as discussed in \cref{subsec:ElementMarking}. The elements are marked based on sampling the damage field at a set of points within the element, which can be the quadrature points or another suitable set of points. The elements are marked for refinement if the damage field exceeds a given threshold $\SCALAR{d}_{\text{min}}$ at any of the sampled points. The marked elements are then refined using the \textbf{refine\_hierarchical\_mesh} routine, and the new spline space is built using \textbf{build\_hierarchical\_space}. {\NEW The area of the marked elements is computed by the \textbf{compute\_element\_area} routine and the maximum over all levels is stored in $A$. }\\

\begin{algorithm}
\caption{\textbf{refine\_elements}($\SCALAR{d}$,$\mathcal{T}$,$\mathcal{Q}$,$m$,$\ell_{\text{max}}$,$\SCALAR{d}_{\text{min}}$)\\
\textbf{Description:} Refines elements based on the damage field while ensuring cross-talk elimination and admissibility.}
\label{alg:refine_elements}
\KwInput{Damage field $\SCALAR{d}$, spline space $\mathcal{T}$, mesh $\mathcal{Q}$, admissibility class $m$, the maximum refinement level $\ell_{\text{max}}$ and the damage threshold $\SCALAR{d}_{\text{min}}$ for marking.}

\Comment{Initialize the marked elements}
$\qty{E} = \emptyset$\\
$A = 0$\\

\Comment{Loop over all levels up to $\ell_{\text{max}}-1$}
\For{$\ell=0,\dots,\ell_{\text{max}}-1$}
{
    \Comment{Loop over all elements $Q_k$ of level $\ell$ in the mesh $\mathcal{Q}$}
    \For{$Q_k\in\mathcal{Q}^\ell$}
    {
        \Comment{Check the damage in sampled points of element $Q_k$}
        \If{$\exists \, \xi \in Q_k : \SCALAR{d}(\xi)>\SCALAR{d}_{\text{min}}$}
        {
            \Comment{If $\ell=\ell_{\text{max}}$-1, mark the support extension to eliminate cross-talk}
            \If{$\ell==\ell_{\text{max}}-1$}
            {
                \Comment{Get the support extension of element $Q_k$}
                $\qty{S} = \textbf{get\_support\_extension}\qty(\mathcal{T},\mathcal{Q},Q_k)$\\
            }
            \Else
            {
                \Comment{Otherwise, only mark the element itself}
                $\qty{S} = \qty{Q_k}$\\
            }

            \Comment{Add the element(s) to the marked set}
            $\qty{E} = \qty{E} \cup \qty{S}$\\
        }
    }

    \Comment{Mark admissibly, see \cite{Bracco2018}}
    \For{$\ell=0,\dots,\ell_{\text{max}}-1$}
    {
        $\qty{E} = \textbf{mark\_recursive}\qty(\mathcal{T},\mathcal{Q},\qty{E},m)$\\
    }

    \Comment{Refine the marked elements}
    $\mathcal{Q} = \textbf{refine\_hierarchical\_mesh}\qty(\mathcal{T},\mathcal{Q},\qty{E})$\\
    $\mathcal{T} = \textbf{build\_hierarchical\_space}\qty(\mathcal{Q})$\\
    \Comment{Compute the area of the marked elements and update $A$}
    $A = \max\qty(A,\textbf{compute\_element\_area}(\qty{E}))$\\
}

\KwOutput{Updated spline space $\mathcal{T}$, updated mesh $\mathcal{Q}$ and marked area $A$.}
\end{algorithm}

Similar to the \textbf{load\_step} algorithm, we define \textbf{adaptive\_load\_step} in \cref{alg:adaptive_load_step}, which includes the mesh refinement procedure. The algorithm performs staggered iterations until convergence, and then refines the mesh based on the damage field. The solution fields are projected to the new mesh using the \textbf{project\_field} routine, which can be based on $L_2$-projection, quasi-interpolation or another suitable method. The algorithm continues to refine the mesh until a stopping criterion based on the ratio of new elements to total elements is met with tolerance $\mathtt{TOL}_{\text{ref}}$, which allows for explicit ($\mathtt{TOL}_{\text{ref}}=1$), implicit ($\mathtt{TOL}_{\text{ref}}=0$) or hybrid mesh adaptivity. \\

\begin{algorithm}
\caption{\textbf{adaptive\_load\_step}($k$,$\SCALAR{d}^k$,$\VEC{u}^k$,$\mathcal{T}$,$\mathcal{Q}$,$m$,$\ell_{\text{max}}$,$\SCALAR{d}_{\text{min}}$,$\mathtt{TOL}_{\text{Pic},\VEC{u}}$,$\mathtt{TOL}_{\text{PSOR},\Delta\SCALAR{d}}$,$\mathtt{TOL}_{\text{stag}}$,$\mathtt{TOL}_{\text{ref}}$)\\
\textbf{Description:} Performs staggered solution and adaptive refinement for a load step $k$.}
\label{alg:adaptive_load_step}
\KwInput{Load step index $k$, damage field $d\qty(\bm{\xi})$, displacement field $\VEC{u}(\bm{\xi})$, spline space $\mathcal{T}$, mesh $\mathcal{Q}$, admissibility class $m$, $\ell_{\text{max}}$, damage threshold $\SCALAR{d}_{\text{min}}$ for marking, tolerances for the Picard iterations $\mathtt{TOL}_{\text{Pic},\VEC{u}}$, PSOR iterations $\mathtt{TOL}_{\text{PSOR},\Delta\SCALAR{d}}$, staggered iterations $\mathtt{TOL}_{\text{stag}}$ and refinement iterations $\mathtt{TOL}_{\text{ref}}$.}

\Comment{Initialize load step (assumes that the loads/displacements of load step $k+1$ are set)}
$\VEC{u}^{k+1} = \VEC{u}^k$\\
$\SCALAR{d}^{k+1} = \SCALAR{d}^k$\\

\Comment{Assemble solution-independent quantities}
$\VEC{F}^{\text{ext}}_{k+1}\gets\textbf{assemble\_elasticity\_RHS}(k+1)$\\
$\bm{\Phi}^{\text{II,IV}}\gets\textbf{assemble\_phasefield\_LHS}()$\\
$\bm{\phi}^{\text{II,IV}}\gets\textbf{assemble\_phasefield\_RHS}()$\\

\Comment{Loop over refinement iterations}
\While{\textbf{true}}
{
    \Comment{Staggered iterations}
    \While{$\Vert R(\VEC{u}^{k+1},\SCALAR{d}^{k+1})\Vert>TOL_R$}
    {
        \Comment{Solve elasticity problem}
        $\VEC{u}^{k+1} = \textbf{solve\_elasticity}(\VEC{u}^{k+1},\SCALAR{d}^{k+1},\mathcal{T},\mathcal{Q},\mathtt{TOL}_{\text{Pic},\VEC{u}})$

        \Comment{Solve phase-field problem}
        $\SCALAR{d}^{k+1} = \textbf{solve\_phase\_field}(\VEC{u}^{k+1},\SCALAR{d}^{k+1},\mathcal{T},\mathcal{Q},\mathtt{TOL}_{\text{PSOR},\Delta\SCALAR{d}})$
    }

    \Comment{Refine elements based on the damage field}
    $\mathcal{Q},\mathcal{T}, A = \textbf{refine\_elements}(\SCALAR{d}^{k+1},\mathcal{T},\mathcal{Q},m,\ell_{\text{max}},\SCALAR{d}_{\text{min}})$\\

    \Comment{Project the solution fields to the new space}
    $\VEC{u}^{k+1} = \textbf{project\_field}(\VEC{u}^{k+1},\mathcal{T})$\\
    $\SCALAR{d}^{k+1} = \textbf{project\_field}(\SCALAR{d}^{k+1},\mathcal{T})$

    \Comment{Break refinement iterations if the ratio of the refined area to the total area is below the tolerance}
    \If{$\frac{A}{\vert \mathcal{Q} \vert}<\mathtt{TOL}_{\text{ref}}$}
    {
        \textbf{Break}
    }
}
\KwOutput{New displacement field $\VEC{u}^{k+1}$, new damage field $\SCALAR{d}^{k+1}$, new spline space $\mathcal{T}$ and new mesh $\mathcal{Q}$.}
\end{algorithm}

\section{Sensitivity analysis of the marking parameters}
\label{app:sensitivity}
\input{Sections/Sensitivity}

\end{appendix}

\bibliography{bibliografia.bib}

\end{document}